\documentclass[12pt,leqno]{report} 
\usepackage {amsfonts}
\usepackage{amsthm}
\usepackage{amssymb}
\usepackage{latexsym}
\usepackage{amsmath}
\usepackage{mathrsfs}
\usepackage{pifont}

\usepackage{setspace}
\usepackage{comment}
\usepackage{fancyhdr}
\renewcommand{\theequation}{\arabic{chapter}.\arabic{equation}}

\newcommand{\tu}{\textup}
\newenvironment{alist}
{

\begin{enumerate}}
{\end{enumerate}}

\newenvironment{rlist}
{

\begin{enumerate}}
{\end{enumerate}}

\newcommand{\Hil}{\mathsf{H}}
\newcommand{\hil}{\mathsf{h}}
\newcommand{\Kil}{\mathsf{K}}
\newcommand{\kil}{\mathsf{k}}
\newcommand{\Vil}{\mathsf{V}}

\newcommand{\Wil}{\mathsf{W}}

\newcommand{\bM}{\mathsf{M}}

\newcommand{\cP}{\mathcal{P}}
\newcommand{\bulletrta}{\bullet\!\rightarrow}
\newcommand{\circrta}{\circ\!\rightarrow}

\newcommand{\up}{\upsilon}
\newcommand{\Do}{\mathcal{D}}

\newcommand{\wh}{\widehat}
\newcommand{\wt}{\widetilde}

\newcommand{\Com}{\Delta}
\newcommand{\Cou}{\epsilon}

\newcommand{\blg} {\mathsf B}
\newcommand{\alg} {\mathsf A}
\newcommand{\coalg} {\mathsf C}
\newcommand{\Xlg} {\mathsf X}

\newcommand{\Coalg}{\mathcal{C}}
\newcommand{\Alg}{\mathcal{A}}
\newcommand{\Blg}{\mathcal{B}}

\newcommand{\la} {\langle}
\newcommand{\ra} {\rangle}
\newcommand{\ot} {\otimes}
\newcommand{\vnot} {\overline{\otimes}}

\newcommand{\bc} {\Bbb C}
\newcommand{\bn}{\Bbb N}
\newcommand{\br}{\Bbb R}
\newcommand{\be}{\Bbb E}

\newcommand{\Odiamond}{\text{\ding{169}}}

\newcommand{\lla} {\left\langle}
\newcommand{\rra} {\right\rangle}

\newcommand{\ve}{\varepsilon}
\newcommand{\Fock}{\mathcal F}

\newcommand{\schur}{\boldsymbol{\cdot}}

\newcommand{\fhat}{\hat{f}}
\newcommand{\ghat}{\hat{g}}
 \newcommand{\chat}{\hat{c}}
\newcommand{\dhat}{\hat{d}}
\newcommand{\khat}{\wh{\kil}}
\newcommand{\Dhat}{\widehat{D}}
\newcommand{\kilhat}{\wh{\kil}}
\newcommand{\zerohat}{\wh{0}}

\newcommand{\QSproj}{\Delta^{QS}}

\newcommand{\bm}[1]{\mbox{\boldmath{$#1$}}}
\newcommand{\cb}{{\text{\tu{cb}}}}
\newcommand{\wb}{{\text{\tu{wb}}}}
\newcommand{\uwb}{{\text{\tu{uwb}}}}
\newcommand{\sr}{{\text{\tu{sr}}}}
\newcommand{\cbsr}{{\text{\tu{cbsr}}}}

\newcommand{\condexp}{\mathbb{E}}

\newcommand{\Geom}{\mathcal{G}}
\newcommand{\Geomdagger}{\mathcal{G}^\dag}
\newcommand{\Geominv}{\mathcal{G}^{\mathrm{inv}}}

\newcommand{\cH}{\mathcal{H}}

\newcommand{\Focktot}{\Fock_{[0,t)}}
\newcommand{\Fockaft}{\Fock_{[t, \infty)}}
\newcommand{\Focktos}{\Fock_{[0,s)}}
\newcommand{\Fockafs}{\Fock_{[s, \infty)}}
\newcommand{\Focktots}{\Fock_{[0,s+t)}}
\newcommand{\Fockafts}{\Fock_{[s+t, \infty)}}
\newcommand{\Fockstot}{\Fock_{[s,s+t)}}

\newcommand{\Exps}{\mathcal{E}}
\newcommand{\Step}{\mathbb{S}}
\newcommand{\Guich}{\mathcal{G}}
\newcommand{\Sant}{\mathcal{S}}

\newcommand{\Seq}{\mathscr{S}}
\newcommand{\Seqdagger}{\mathscr{S}^\dag}
\newcommand{\Seqinv}{\mathscr{S}^{\mathrm{inv}}}
\newcommand{\Seqstar}{\mathscr{S}^\ast}

\newcommand{\Procwc}{\mathbb{P}_{\mathrm{wc}}}

\newcommand{\Procwr}{\mathbb{P}_{\mathrm{wr}}}
\newcommand{\Proclb}{\mathbb{P}_{\mathrm{lb}}}
\newcommand{\Proctwo}{\mathbb{P}_2}
\newcommand{\Procc}{\mathbb{P}_{\mathrm{c}}}
\newcommand{\ProcHc}{\mathbb{P}_{\mathrm{Hc}}}
\newcommand{\Procb}{\mathbb{P}_{\mathrm{b}}}

\newcommand {\Op} {\mathcal{O}}
\newcommand {\Opdag} {\mathcal{O}^{\dag}}
\newcommand{\Opinv}{\mathcal{O}^{\mathrm{inv}}}
\newcommand{\Opstar}{\mathcal{O}^\ast}

\newcommand{\Proc}{\mathbb{P}}

\newcommand{\ida}{1_{\alg}}
\newcommand{\id}{\textrm{id}}

\newcommand{\idB}{\textrm{id}_{B(\khat)}}
\newcommand{\idf}{1_{\Fock}}

\newcommand{\QSC}{\mathbb{QSCC}}
\newcommand{\QSCHc}{\QSC_{\mathrm{Hc}}}
\newcommand{\QSCdag}{\QSC^\dag}

\newcommand{\QSCdagc}{\QSCdag_{\mathrm{c}}}
\newcommand{\QSCdagHc}{\QSCdag_{\mathrm{Hc}}}

\newcommand{\QSCwc}{\QSC_{\mathrm{wc}}}

\DeclareMathOperator{\Dom}{Dom}
\DeclareMathOperator{\Ran}{Ran}
\DeclareMathOperator{\Lin}{Lin}
\DeclareMathOperator{\Ker}{Ker}
\DeclareMathOperator{\re}{Re}

\theoremstyle{plain}
\newtheorem{propn}{Proposition}[section]
\newtheorem{tw}[propn]{Theorem}
\newtheorem{lemma}[propn]{Lemma}
\newtheorem{cor}[propn]{Corollary}
\newtheorem{fact}[propn]{Fact}

\theoremstyle{definition}
\newtheorem{deft}[propn]{Definition}
\newtheorem{rem}[propn]{Remark}

\renewcommand{\thepropn}{\arabic{chapter}.\arabic{section}.\arabic{propn}}
\renewcommand{\theequation}{\arabic{chapter}.\arabic{section}.\arabic{equation}}

\begin{document}
\newpage

\pagenumbering{arabic}
\setcounter{page}{1}
\pagestyle{fancy}
\renewcommand{\chaptermark}[1]{\markboth{\chaptername%
\ \thechapter:\,\ #1}{}}
\renewcommand{\sectionmark}[1]{\markright{\ #1}}


\baselineskip=0.6cm

\begin{titlepage}
\begin{center}
\vspace*{3cm}

\Large
\textbf{Quantum Stochastic}

\vspace{0.5cm}
\textbf{Convolution Cocycles}
\normalsize

\vspace{6cm}

\Large
Adam G. Skalski
\normalsize

\vfill

Thesis submitted to the University of Nottingham
\linebreak
for the degree of Doctor of Philosophy
\linebreak
December 2005
\end{center}
\end{titlepage}

\thispagestyle{empty} \vspace*{16cm} \hspace*{5 cm} \Large Moim
Rodzicom \normalsize

\onehalfspacing

\pagebreak

\chapter*{Acknowledgments}
This thesis would never have been written without the help, company
and encouragement of many people.

I would like to thank Martin Lindsay for the most friendly,
enthusiastic and at the same time challenging form of supervision.
He shared with me generously his deep knowledge and understanding of
quantum stochastic analysis. I admire his ability to see the global
picture of the problems under consideration without compromising the
mathematical rigour of the solutions; I hope to follow this approach
in my future work. The willingness to help in all aspects of
settling in Nottingham, and many friendly hours we spent together,
extended far beyond the formal duties of a supervisor. Finally, I
need also to thank Martin for his patience with my impatience on
several occasions.

I was first introduced to noncommutative mathematics at my home
university in \L \'od\'z by Ewa Hensz-Ch\c{a}dzy\'nska. My (ongoing)
fascination with functional analysis in general and operator
algebras and their applications in particular began at that time,
and throughout the period of my studies in Poland and in the United
Kingdom I felt and appreciated the support of Prof.
Hensz-Ch\c{a}dzy\'nska.

Joachim Zacharias suffered from my several visits to his office in a
good-humoured way; although my questions were probably often not
very clever, his answers were always helpful and motivating. It was
a big pleasure to work with him.

The examiners, Philippe Biane and Joel Feinstein, offered their precious time to
read this thesis, write the reports and conduct the viva. I am very grateful for their 
thoughtful remarks which suggested further possible areas of research.

On my way to this point I met many other mathematicians (and
physicists) without whose knowledge, explanations, advice, comments,
encouragement, hospitality and often practical help none of this
work would have been done. I owe special gratitude to Wojciech
Banaszczyk, Rajarama Bhat, Marek Bo\.zejko, Uwe Franz, Stanis\l aw
Goldstein, Robin Hudson, Marius Junge, Edward Kapu\'scik, Semyon
Litvinov, Andrzej \L uczak, Alexander Ushveridze, Stephen Wills,
Leonid Vainerman and Quanhua Xu.

Fellow graduate students made excellent company during my first
visit to Nottingham and later throughout my Polish and English PhD
studies - in particular I would like to mention the operator
algebraic lot: Orawan, Iulian and Nick. Apologies go to all upon
whom I inflicted numerous seminars - I hope they were at least to a
certain degree as useful to the audience as they were to the
speaker. \thispagestyle{empty}

Moi Rodzice i Siostra pokazali mi, \.ze zadawanie pyta\'n i szukanie
odpowiedzi mo\.ze by\'c przygod\c{a}, \.ze nauka, a zw\l aszcza
matematyka mo\.ze by\'c \'zr\'od\l em nieustaj\c{a}cej
przyjemno\'sci. Bardzo Wam za to dzi\c{e}kuj\c{e}.

Finally, for the daily support, the willingness to move to another
country, the relentless belief in what I was doing, and all the
other things that cannot be described, I thank Joasia.

\thispagestyle{empty}

\pagebreak \thispagestyle{empty} \vspace*{14cm} \hspace*{3 cm} Ocho
y ocho son diecis\'eis y el que cuenta.

\hspace*{4.5 cm} Julio Cortasar, \emph{Los premios}

\pagebreak \thispagestyle{empty}

\normalsize
\begin{center} \textbf{Abstract} \end{center}

A concept of quantum stochastic convolution cocycle is introduced
and studied in two different contexts -- purely algebraic and
operator space theoretic. A quantum stochastic convolution cocycle
is a quantum stochastic process $k$ on a coalgebra $\Alg$ satisfying
the convolution cocycle relation
\[ k_{s+t} = (k_s \ot (\sigma_s \circ k_t) ) \circ \Com
\]
and the initial condition $ k_0 = \iota \circ \Cou$,  where $\Com$
and $\Cou$ denote the coproduct and counit of $\Alg$, $\sigma$
denotes the time shift on operators acting on the symmetric Fock space 
over $\br_+$ and $\iota$ is the
unital embedding of $\bc$ in the algebra of bounded operators on the
Fock space. The notion generalises that of quantum L\'evy process,
which in turn is a noncommutative probability counterpart of
classical L\'evy process on a group.

       Convolution cocycles
arise as solutions of quantum stochastic differential equations. In
turn every sufficiently regular cocycle satisfies an equation of
that type. This is proved along with the corresponding existence and
uniqueness of solutions for coalgebraic quantum stochastic
differential equations. The stochastic generators of unital
$^*$-homomorphic cocycles are characterised in terms of structure
maps on a
$^*$-bialgebra. 
This yields a simple proof of the Sch\"urmann Reconstruction Theorem
for a quantum L\'evy process; it also yields a topological version
for a quantum L\'evy process on a $C^*$-bialgebra. Precise
characterisation of the stochastic generators of completely positive
and contractive quantum stochastic convolution cocycles in the
$C^*$-algebraic context is given, leading to some dilation results.
A few examples are presented and some interpretations offered for
quantum stochastic convolution cocycles and their stochastic
generators on different types of  $^*$-bialgebra.

The techniques used for the analysis of convolution cocycles in the
purely algebraic and operator space theoretic context are distinct.
In the first case the basic tool is the Fundamental Theorem on
Coalgebras. In the second the key method is the application of the
so-called $R$-transformation linearising  quantum convolution.

\newpage


\tableofcontents

\chapter{Introduction}
\newpage

This thesis is concerned with the investigation of quantum
stochastic convolution cocycles in both a purely algebraic context
and in the context of compact quantum groups. The notion has been
introduced by J.M.\,Lindsay and the author in \cite{LSqscc1}.
Quantum stochastic cocycles are quantum stochastic processes living
on a (symmetric) Fock space of a Hilbert space of vector-valued
square-integrable functions on the positive half-line and satisfying
a cocycle relation with respect to the Fock space shift. This notion
derives from the idea of the cocycle structure that may be
associated to continuous-time classical Markov processes (the
precise analogy and connections are exploited in
\cite{MartinKalyan}; the example of Brownian motion cocycles cast in
the language used in this thesis may be found in \cite{lect}).
Quantum stochastic cocycles are basic objects of interest in quantum
stochastic analysis (\cite{par}, \cite{mey}, \cite{Biane},
\cite{lect}) and in the more general theory of noncommutative white
noise (\cite{hkk}). They were first introduced in the noncommutative
probability literature in the late 1970s by L.\,Accardi and his
collaborators (\cite{Luigi}, \cite{AFL}). Later developments
unearthed connections to other areas of mathematics. As example we
can mention the application of the analysis of quantum stochastic
cocycles to resolve certain open questions in Arveson's theory of
product systems (\cite{Arv}). On the other hand there is a well
developed theory of quantum L\'evy processes, originating in the
paper \cite{asw} of L.\,Accardi, M.\,Sch\"urmann and W.\,von
Waldenfels and further extensively developed especially by the
second mentioned author (see \cite{franz}, \cite{schu} and
references therein). Quantum L\'evy processes generalise classical
stationary, independent increment stochastic processes on groups.
Close examination of the two areas described above has naturally led
to the idea of quantum stochastic convolution cocycle.

As usual in noncommutative mathematics, in order to `quantise' one
focuses on the appropriate class of functions on the underlying set
of the classical structure. When the interest lies in topological
aspects of a space, one starts with the algebra of continuous
functions, leading to the theory of $C^*$-algebras, and further, via
introducing differential structure and the $K$-theory to the
noncommutative geometry of A.\,Connes. When the underlying set is
considered as a measure space, the relevant object is the algebra of
all measurable functions, leading to the land of von Neumann
algebras. By analogy we see that quantum stochastic convolution
cocycles should `act on' a generalisation of the algebra of
complex-valued functions on a group or, to be more precise,
semigroup with identity, namely on a $^*$-bialgebra. These,
depending on whether the topology is taken into account, may be
considered in two different categories, purely algebraic and
operator-space theoretic. This dichotomy influences the structure of
the thesis, which is described in the end of this chapter.

\subsection*{Classical theory of probability on algebraic structures}

This thesis is aimed to be a contribution to the theory of
noncommutative probability on quantum/algebraic structures. In this
subsection we indicate the amazing variety and richness of the
classical theory. This is done first by presenting major ideas and
concepts of extending the standard probability theory, concerned
with real-valued random variables, to more general algebraic
structures. The choice is clearly subjective, and is not intended to
represent neither fair summary nor general introduction to the
subject. For this we refer to a beautiful early book \cite{Gren} and
to the impressive volume \cite{Heyer} (see also \cite{Heyercomp},
\cite{Heyer2}, \cite{hyp}). In the second part we concentrate on a
particular case of probability measures on compact groups, which is
of importance for the following chapters.

   Let us start with the observation that the afore-mentioned extensions of standard probability theory may be seen at
least on two basic levels.
We may change the target space of the random variables in question and consider random variables taking values
in infinite dimensional vector spaces, in abelian or nonabelian groups, or in (Banach) algebras.  We may also consider
typical distributions on the real line or complex plane, combine them using algebraic operations, and then ask questions
concerning the objects obtained, as is done for example in the theory of random matrices.

The importance of investigating properties and looking for answers
to new problems arising in the constructions described above stems
from multiple sources. Firstly, there exist concrete physical models
that can be usefully described in the algebraic language. To mention
the simplest examples, behaviour of a ball with a fixed centre
floating in a randomly moving fluid may be analysed via the theory
of random variables with values in the rotation group SO(3); a chain
of devices transforming the initial data according to certain linear
prescription and subjected to errors may be viewed as a dynamical
system whose evolution is governed by consecutively applied random
matrices. Secondly, a probabilistic approach may be helpful in
answering deterministic questions, an idea which has been
successfully pursued by the Polish school of functional analysis
(e.g.\. by constructing Banach spaces with desired properties as
suitable limits of randomly generated finite-dimensional spaces).
Thirdly, attempts to generalise allow us to understand better the
particular case of probability theory on $\br$, both its advantages
and limitations. Finally there is the most important reason
according to the author, namely pure intellectual curiosity.

Extensions of classical results usually present various technical
problems. The main tools of the standard theory, such as Fourier
analysis, need to be reformulated and in the process often lose at
least some of their well-known properties. Many theorems need
additional assumptions, automatically satisfied by real- or
complex-valued random variables. In most cases however, functional
analytic techniques remain useful; their full strength is often
revealed only in such a general context.

Let us now illustrate somewhat vague statements above, and provide a bridge with the considerations in further chapters. Assume
that $G$ is a compact semigroup. Its multiplicative structure provides the convolution action on $M(G)$, the Banach space of all
 regular complex-valued Borel measures on $G$ with total variation norm:
\[ \mu \star \nu (f) = \int_G \int_G f(st) d\mu(s) d\nu(t), \;\;\; \mu, \nu \in M(G),\; f\in C(G),\]
where the isometric isomorphism $M(G) \cong C(G)^*$ is used ($C(G)$ denotes
the Banach space of all continuous functions on $G$ with the supremum norm).  It is easy to see that
$(M(G),\star)$
becomes a Banach algebra.
Further if $G$ has a neutral element $e$, the measure $\delta_e$ (the Dirac measure concentrated
in $e$) is the unit of $(M(G),\star)$. This allows us to introduce the following notion:
a family $\{\mu_t: t \geq 0\}$ of measures in $M(G)$ is called a
\emph{convolution semigroup of measures} if
\[\mu_{s+t} = \mu_s \star \mu_t, \;\;s,t \geq 0, \;\;\; \text{ and } \;\;\; \mu_0 = \delta_e.\]

It is called \emph{weakly continuous} (or more correctly weak$^*$-continuous) if
$|\mu_t(f)-\mu_0(f)|\stackrel{t \to 0^+}{\longrightarrow} 0$ for each $f \in C(G)$,
it is called \emph{norm continuous} if $\|\mu_t-\mu_0\|\stackrel{t \to 0^+}{\longrightarrow} 0$.

Already these, apparently simple notions allow interesting, and in general highly
nontrivial, questions to be asked. For example,
\begin{alist}
\item how do (norm continuous, weakly continuous) convolution semigroups of measures look like?
      Can they be classified?
\item given a measure $\mu \in M(G)$, when does it embed in a convolution semigroup of measures (that is,
when does the convolution semigroup of measures $\{\mu_t: t \geq 0\}$ exist such that
$\mu = \mu_1$)?
 If it does embed, is the embedding unique?
\item what is the limit behaviour of the semigroup? Given a measure $\mu \in M(G)$, when
do the sequences $(\mu^{\star k})_{k=1}^{\infty}$,
$(\frac{1}{k} \sum_{n=1}^{k}\mu^{\star k})_{k=1}^{\infty}$ converge? What are the possible limits?
\end{alist}
The typical method used for approaching these problems is to transform questions concerning
convolution semigroups to more familiar questions concerning certain semigroups of operators.
Define for each $\mu \in M(G)$ the map $P_{\mu}: C(G) \to C(G)$ by
\[ P_{\mu} (f) (s) = \int_G f(st) d\mu(t), \]
$f \in C(G)$, $ s \in G$. When $\mu$ is a probability measure, the
reader will have recognised here the so-called Markov operator. Not
surprisingly, the problems above become easier when attention is
restricted to norm-continuous semigroups. These can be shown to
correspond precisely to distributions of compound Poisson processes
(\cite{Gren}). In full generality, the answers to (a) and (b) depend
on additional properties of the (semi)group $G$ (one should mention
for example the celebrated Hunt formula from \cite{Hunt},
characterising generators of weakly continuous convolution
semigroups of measures on Lie groups - see \cite{Heyer}). Here, as
in the case of $\br$, the Fourier transform, this time understood in
the language of the Peter-Weyl theory of unitary representations of
compact groups, is an indispensable tool. The limits of sequences in
(c), when they exist, are clearly idempotent measures (that is
measures satisfying $\mu \star \mu = \mu$), and one can show that if
$G$ is a group, the idempotent probability measures correspond to
Haar measures on subgroups. The corresponding property does not hold
in the quantum case (for more information on that see \cite{Pal}).

The above well-known facts are quoted here, as the same methods will
be echoed in Chapter 3, and especially in Chapter 4, in a
noncommutative context. There we will approach the stochastic
process by looking at the convolution semigroups of functionals
given by its distribution, and then proceed with the idea of
considering the corresponding semigroups of operators (formalised in
the guise of the $R$-map introduced in Section \ref{OS coalgebras}).
Further comments on relations between the topic of this subsection
and the main subject of this thesis are given in Section \ref{Exam}.

\subsection*{Summary}

Here we present a summary of the main results of the thesis,
expanding on the abstract.

A concept of quantum stochastic convolution cocycle is introduced
and studied in two different contexts, namely purely algebraic and
operator space theoretic. A quantum stochastic convolution cocycle
is a quantum stochastic process $k$ on a coalgebra $\Alg$ (with
coproduct $\Com:\Alg \to \Alg \ot \Alg$ and counit $\Cou:\Alg \to
\bc$) satisfying the convolution cocycle relation
\[ k_{s+t} = (k_s \ot (\sigma_s \circ k_t) ) \circ \Com
\]
and the initial condition
\[ k_0 = \iota_{\Fock} \circ \Cou.\]
The functional equation is for operators on a symmetric Fock space
$\Fock$ over a Hilbert space of square-integrable vector-valued
functions on $\br_+$. Here $\sigma$ denotes the time shift on the
Fock space operators and $\iota_{\Fock}$ is the unital embedding of
$\bc$ in the algebra of bounded operators on $\Fock$. The notion
generalises the notion of quantum L\'evy process, which in turn is a
noncommutative probability counterpart of classical L\'evy process
on a group. Moreover it is a convolution counterpart to a previously
studied structure, referred to here as \emph{standard} quantum
stochastic cocycle. Such a cocycle is a quantum stochastic process,
again consisting of Fock space operators, satisfying the cocycle
relation $ k_{s+t} = \widehat{k}_s \circ \sigma_s \circ k_t $ (where
$\widehat{k}_s$ denotes a  certain natural extension of $k_s$) and
the initial condition $k_0 = \iota_{\Fock}$.

      Convolution cocycles
arise as solutions of certain kinds of linear constant-coefficient
quantum stochastic differential equations, as do standard cocycles.
 We show, under suitable assumptions on the coefficients, existence and uniqueness
for solutions of coalgebraic quantum stochastic differential
equations with initial condition given by the counit; these turn out
to be quantum stochastic convolution cocycles. Conversely, we show
that every sufficiently regular convolution cocycle satisfies an
equation of that type. The stochastic generators of unital
$^*$-homomorphic cocycles on a $^*$-bialgebra may be characterised
in terms of structure maps, equivalently by Sch\"urmann triples.
This yields a simple proof of the Sch\"urmann Reconstruction Theorem
for a quantum L\'evy process; it also yields a topological version
for quantum L\'evy process on a $C^*$-bialgebra. Precise
characterisation of the stochastic generators of completely positive
and contractive quantum stochastic convolution cocycles in the
$C^*$-algebraic context is given, leading to two theorems on
dilations of such cocycles to $^*$-homomorphic convolution cocycles.
Finally we present a few examples and offer some interpretations of
quantum stochastic convolution cocycles and their stochastic
generators on different types of $^*$-bialgebra.

The techniques used for the analysis of cocycles in the purely
algebraic and operator space theoretic context are distinct. In the
first case the basic tool is the Fundamental Theorem on Coalgebras,
as exploited by Sch\"urmann in his pioneering work. In the second
the key method is the application of the so-called
$R$-transformation linearising quantum convolution to the
well-established theory of  standard quantum stochastic cocycles.

\subsection*{Description of the contents}

Having briefly described classical inspirations and summarised the
main results, the rest of this chapter consists of a detailed
description of the contents of the thesis. Quantum stochastic is
abbreviated to QS in the sequel. Chapter 2 contains preliminaries.
After introducing in Section \ref{genfactnot} the general notations,
we proceed in Section \ref{Ospaces} to recall the basic notions of
operator space theory and the definition of matrix spaces. Section
\ref{Focksection} contains standard facts concerning Fock spaces and
defines QS processes. The last part of Chapter 2, Section
\ref{QSIntegrals}, introduces the QS integration of R.L.\,Hudson and
K.R.\,Parthasarathy. Fundamental Formulae and Fundamental Estimate
are quoted, and iterated QS integrals described; the notation used
is mainly modelled on the lecture notes \cite{lect}.

In Chapter 3 QS convolution cocycles are investigated in the purely
algebraic context. The chapter starts by introducing coalgebras and
convolution structure provided by the coproduct in Section
\ref{algandconv} and proceeds with the definition of QS convolution
cocycles and their associated semigroups in Section \ref{QSccocycl}.
Section \ref{coalgQSDE} contains the existence and uniqueness
results for solutions of coalgebraic QS differential equations; it
also describes basic properties of the solutions, including a form
of H\"older continuity. The fact that the solutions of coalgebraic
QS differential equations are QS convolution cocycles is proved in
Section \ref{stochsection}; there also a converse is obtained - all
H\"older continuous QS convolution cocycles must arise in this way.
Section \ref{multalg} is devoted to multiplicative properties of the
cocycles on $^*$-bialgebras. As the operators in question may be
(and usually are) unbounded, the analysis concentrates rather on
properties of scalar products than the actual composition. The It\^o
Formula for iterated QS integrals is expressed via so-called
matrix-sum kernels, following \cite{lwhom}, and the general
structure of the stochastic generators of weakly multiplicative
cocycles is obtained. The known characterisation of generators of
$^*$-homomorphic QS convolution cocycles in terms of Sch\"urmann
triples follows, and in Section \ref{QLevy} we recall the definition
of quantum L\'evy processes due to L.\,Accardi, M.\,Sch\"urmann and
W.\,von Waldenfels (\cite{asw}) and give an alternative proof of the
Sch\"urmann Reconstruction Theorem (\cite{schu}). Section
\ref{perturb} is devoted to certain perturbation formulas (see also
\cite{franz}). The last section of Chapter 3 contains a definition
of opposite QS convolution cocycles and indicates corresponding
versions of the results of previous sections.

In Chapter 4 we continue the investigation of QS convolution
cocycles, this time in the operator space theoretic (or
$C^*$-algebraic) context. It begins with Section \ref{OS coalgebras}
introducing operator space coalgebras and convolution semigroups of
functionals. The stress is put on the so-called $R$-map which
transforms convolution semigroups into semigroups of operators,
whilst respecting their relevant continuity properties. In Section
\ref{standvsconv} QS convolution cocycles on operator space
coalgebras are defined. This time it is necessary to introduce also
the notion of a weak QS convolution cocycle.  Later the basic facts
on standard QS cocycles are recalled (following \cite{lwjfa}), and
the $R$-map is shown to yield a transformation between the two
classes of cocycles. Section \ref{QSDEOS} contains the existence and
uniqueness results for solutions of QS differential equations in the
context of operator spaces, with completely bounded coefficients and
nontrivial initial conditions. The methods extend those of
\cite{lwex}. Coalgebraic QS differential equations are the topic of
Section \ref{coalgQSDEquat}. They may be viewed as a particular case
of the equations described in the previous section. Their solutions
are weak QS convolution cocycles. From the corresponding fact in the
theory of standard QS cocycles and  properties of the $R$-map, in
Section \ref{CPCcoc} we prove that every Markov-regular completely
positive and contractive QS convolution cocycle on a
$C^*$-hyperbialgebra satisfies a coalgebraic QS differential
equation. The precise structure of stochastic generators of such
cocycles is also described. Section \ref{multCalg} is devoted to
$^*$-homomorphic, or more generally weakly multiplicative, cocycles
on $C^*$-bialgebras. Their stochastic generators may be again, as in
Chapter 3, expressed in terms of Sch\"urmann triples. Two possible
notions of quantum L\'evy process on a $C^*$-bialgebra are proposed
and a version of the Sch\"urmann Reconstruction theorem for such
processes is given. The question of dilating completely positive QS
convolution cocycles on $C^*$-bialgebras to $^*$-homomorphic ones is
addressed in Section \ref{Dilations}. Two possible forms of
dilations are described, corresponding to the results obtained for
standard QS cocycles respectively in \cite{dilate} and in
\cite{Stine}. The proofs exploit the characterisation of the
stochastic generators obtained earlier. Section \ref{Exam} contains
several examples and comments on QS convolution cocycles on specific
types of OS coalgebras. Unital $^*$-homomorphic cocycles are
discussed in three main cases: commutative (algebras of continuous
functions on compact groups), cocommutative (universal
$C^*$-algebras of discrete groups) and genuinely quantum (where the
most satisfactory results are obtained for full compact quantum
groups). The last part of this section is devoted to recalling the
standard conditional expectation construction of
$C^*$-hyperbialgebras, and describing the QS convolution cocycles on
the structures so-obtained. Finally the last short section discusses
the notion of QS convolution cocycles on multiplier $C^*$-bialgebras
(locally compact quantum semigroups) and poses some of the problems
associated with the developing of theory for these.

There are two appendices. In Appendix A the automatic continuity
properties of $(\pi_1,\pi_2)$-derivations (for representations
$\pi_1$, $\pi_2$ of a $C^*$-algebra) are discussed, based on the
celebrated work of S.\,Sakai, J.R.\,Ringrose and E.\,Christensen on
(standard) derivations on $C^*$-algebras. These results are needed
for the description of the stochastic generators of $^*$-homomorphic
cocycles in Section \ref{multCalg}. Appendix B recalls basics of the
theory of Hopf $^*$-algebras and their corepresentations
(\cite{Schm}), useful for understanding the examples in Section
\ref{Exam}.

Most of the contents of Chapter 3 have been published in
\cite{LSqscc1} and the presentation here is patterned on this paper,
with certain transpositions of the material, so that the notion of
convolution cocycles is introduced at the earliest possible stage.
The only distinct part of Chapter 3, compared to \cite{LSqscc1}, is
the section on the opposite cocycles (Section \ref{oppQSccocycl}).
Some of the main results of Chapter 4 have been announced in
\cite{LSbedlewo}; full versions including proofs will be published
in papers \cite{LSqsde} (contents of Section \ref{QSDEOS}),
\cite{LSqscc2} (Sections \ref{OS coalgebras}, \ref{standvsconv},
\ref{coalgQSDEquat}, \ref{multCalg}, \ref{Exam}, Appendix A and
parts of Section \ref{CPCcoc}) and \cite{CPQSCC} (parts of Section
\ref{CPCcoc}, Section \ref{Dilations}). The main difference between
the treatment in the thesis and the papers mentioned above lies in
that the latter are written in a more abstract language, with
abstract operator spaces replacing concrete ones as the fundamental
objects. Both approaches may be shown to be equivalent; the one
adopted here permits certain simplifications of the proofs and
remains closer to the traditional QS framework. Another important
difference between the thesis and the papers mentioned above is that
here we concentrate, whenever possible, on the class of completely
bounded, everywhere defined processes, whereas in particular in
\cite{LSqsde} the perspective is broadened to discuss processes
which are densely defined, assuming only complete boundedness of
their columns. Whilst analytically natural, this generality is not
needed for the purposes of this thesis. The paper \cite{LSqsde} also
contains a characterisation of those QS cocycles on finite
dimensional operator spaces that arise as solutions of QS
differential equations, similar in spirit to the characterisation
for QS convolution cocycles in Section \ref{stochsection}. Another
difference is that proofs of various properties of solutions of QS
differential equations in Section \ref{QSDEOS} are based on the
uniqueness theorem, rather than on the explicit checks as in
\cite{LSqsde}.

Here the introduction ends, and the proper journey begins. Everyone is invited!

\chapter{Preliminaries}
 \newpage
 This chapter contains the preliminary facts and establishes some of the terminology
 needed in the sequel.
After discussing the general notations used in the thesis,
the standard notions of the theory of operator spaces, including the concept
of the matrix spaces, are given. We proceed to describe the symmetric Fock space
setup, main Fock space operators and define quantum stochastic processes.
This is sufficient to understand the definition of QS cocycles and QS convolution cocycles
as introduced in Section \ref{standvsconv}.
The rest of the chapter is devoted to the quantum stochastic integration. The basic definitions
are given and the Fundamental Formulas and the Fundamental Estimate included; finally the preliminary
aspects of the iterated QS integrals are discussed.

\section{General facts and notations}   \label{genfactnot}
\setcounter{equation}{0}

 All vector spaces in this work are complex,
  scalar products are linear in the second variable.
In general we use $V,W,A,\ldots$ for
vector spaces (and other purely algebraic objects: algebras, coalgebras, etc.) and
$\Vil, \Wil, \alg, \ldots$ for operator spaces ($C^*$-algebras, operator space coalgebras,
 etc.); Hilbert spaces are usually denoted by letters $\hil$, $\kil$, $\Hil$, $\Kil$
 ($\kil$ will be usually reserved for the noise dimension space, cf. Section \ref{Focksection}).
The vector space of linear maps between vector spaces $V,W$ will be denoted $L(V;W)$, and
the normed space of all bounded linear maps between normed spaces $\Vil, \Wil$ will be
denoted by $B(\Vil; \Wil)$. If $S\subset V$, $\text{Lin}\,S$ denotes the linear span of $S$.

For a subset $V$ of an involutive vector space $W$, $V^{\dagger}:=\{w \in W: w^{\dagger} \in V\}$.
When both vector spaces $V,W$ are equipped with involutions, and $\phi\in L(V;W)$, the map
$\phi^{\dagger}\in L(V;W)$ is defined by $\phi^{\dagger} (v) = \phi(v^{\dagger})^{\dagger}$, $v \in V$.

For a function $f:\br_+ \to \hil$ and subinterval $I$ of $\br_+$,
$f_I$ denotes the function $\br_+\to\hil$ which agrees with $f$ on $I$
and is zero outside $I$ (cf. standard indicator-function notation).
This convention also applies to vectors, by viewing them as constant
functions---for example
\[
\xi_{[s,t[}, \text{ for } \xi\in\hil  \text{ and } 0\leq s<t.
\]

Let $\hil$ be a Hilbert space.
Ampliations are denoted
\begin{equation*}
\iota_{\hil}: B(\Hil)\to B(\Hil\ot\hil), \ T\mapsto T\ot I_{\hil},
\end{equation*}
and each vector $\xi\in\hil$ defines operators
\begin{equation}
\label{Dirac}
E_\xi : \Hil \to \Hil\ot\hil, \ v\mapsto v\ot\xi \text{ and }
E^{\xi} = (E_\xi)^*,
\end{equation}
generalising Dirac's bra-ket notation:
\begin{equation*}
E_\xi = I_{\Hil} \ot |\xi\ra \text{ and } E^\xi = I_{\Hil} \ot \la\xi |.
\end{equation*}
The particular Hilbert space $\Hil$ will always be clear from
the context. When $\xi,\eta\in \hil$, the functional
$\omega_{\xi, \eta}:B(\hil) \to \bc$ is defined by
\begin{equation} \label{omegas}
 \omega_{\xi, \eta} (T) = \la \xi, T \eta \ra,
\end{equation}
$T\in B(\hil)$. We will write simply $\omega_{\xi}$ for $\omega_{\xi, \xi}$
($\xi \in \hil$).

Now let $D$ be a subset of $\kil$.
The following notation will be employed:
\begin{equation}
\label{hats}
\wh{D} := \Lin \big\{\wh{\xi}: \xi\in D\big\},
\text{ where } \wh{\xi}:= \binom{1}{\xi}\in\wh{\kil}:=\bc \oplus \kil.
\end{equation}
Whenever $J$ is a set and $f:J \to \kil$, the function $\fhat:J \to \kilhat$ is defined by
\[ \fhat(s): = \wh{f(s)}, \;\;\; s \in J.\]
The orthogonal projection from $\wh{\kil}$ onto $\kil$ is denoted by $\QSproj$; if
$\hil_1, \hil_2$ are Hilbert spaces, $\QSproj_{\hil_1,\hil_2}:= \id_{\hil_1} \ot \QSproj \ot
\id_{\hil_2}$.

Let now $E$ be a dense subspace of some Hilbert space $\hil$.  $\Op (E)$ will denote the vector
 space of linear operators in $\hil$ with
domain $E$, $\Opdag (E)$ the vector space of these operators
in $\Op(E)$ whose adjoints' domains include $E$. The space  $\Opdag (E)$ is equipped with the
involution $\dagger$:
\begin{equation}
\label{daggerinv} T^{\dagger} := T^*|_E,  \;\;\; T \in \Opdag (E).\end{equation}
In chapter 3 some further subspaces of $\Op (E)$ will be specified. Usually when $T$ is an operator
defined on the whole of $\hil$ its restriction to $E$ will be denoted by the same letter;
this notational abuse leads to the formal inclusion $B(\hil) \subset \Op(E)$. Note that
$\Opdag(\hil)=B(\hil)$ and the formula \eqref{omegas} (for $ \xi \in \hil$, $\eta \in E$)
defines a linear functional on $\Op(E)$, denoted also by $\omega_{\xi,\eta}$.

If $\Alg$ is a $^*$-algebra, an element $a\in \Alg$ is said to be
\emph{positive} if $a = \sum_{i=1}^n a_i^* a_i$ for some $n \in
\bn$, $a_1,\ldots,a_n \in \Alg$. The set $\Alg_+$ of positive
elements is always a cone, but in general it may fail to generate
the whole of $\Alg$ if $\Alg$ is nonunital (think of the algebra of
all polynomials without constant term). By a positive map between
$^*$-algebras is understood a linear map preserving cones of
positive elements.

At some point we will need a particular notation for converting
a subset of the set $\{1,\ldots,n\}$ (denoted below by $\nu$) into  subsets of $\{1,\ldots,n+1\}$ ($n\in \bn$)
developed in \cite{lwhom}:
\begin{equation} \label{arrows}
\overset{\circrta}{\nu}
 := \{ 1+k : k \in \nu \} \;\text{ and }\;
 \overset{\bulletrta}{\nu}  := \{ 1\} \cup \overset{\circrta}{\nu}  .
\end{equation}

 Algebraic tensor product of vector spaces is denoted by $\odot$; for maps mainly
the symbol $\ot$ is used, meaning either algebraic tensor product of linear maps, or
its relevant continuous extension. When there is no danger of confusion
the symbol of tensor product between vectors is dropped (e.g.\ $u \ot \ve(f)$ is written
as $u\ve(f)$, the notation to be introduced in Section \ref{Focksection}). The set of nonnegative integers
is denoted $\bn_0:= \bn \cup \{0\}$. When $S_1,S_2$ are arbitrary sets, $S_1\subset \subset S_2$
means that $S_1$ is a \emph{finite} subset of $S_2$.
The symbol $\Delta_n [0,t]$ denotes the simplex
$\big\{ \bm{s} \in [0,t]^n \, \big| \;  s_n \geq \cdots \geq s_1 \big\}$;
for many estimates in the main text it is worth recalling that the Lebesgue measure of
$\Delta_n [0,t]$ is $\frac{t^n}{n!}$.
The symbol $\tau$ is often used for tensor flips both between Hilbert spaces
($\hil_1 \ot \hil_2 \to \hil_2 \ot \hil_1$) and algebras of bounded operators
($B(\hil_1 \ot \hil_2) \to B(\hil_2 \ot \hil_1)$) - note that the latter is a $^*$-isomorphism.
When $\alg\subset B(\hil)$ is a nondegenerate $C^*$-algebra,
$\alg''$ (the bicommutant of $\alg$) equals the strong closure of $\alg$ - this is the content
of von Neumann's Double Commutant Theorem.

By a \emph{nonnegative-definite kernel} on a non-empty set $S$ is understood a map $k:S \times S \to \bc$
such that for all $n \in \bn$, $\lambda_1, \ldots, \lambda_n \in \bc$, $ s_1, \ldots,s_n \in S$
there is $\sum_{i,j=1}^n {\overline{\lambda_i}\lambda_j k(s_i, s_j)} \geq 0$.
A pair $(\Kil,\eta)$, where $\Kil$ is a Hilbert space and $\eta:S \to \Kil$ is called
a \emph{(minimal) Kolmogorov construction} for the kernel $k$ if
\[k(s,t) = \la \eta(s), \eta(t) \ra, \;\; s,t \in S \; \text{ and  }  \; \overline{\Lin}\, \eta(S)= \Kil.\]
A minimal Kolmogorov construction
exists for each nonnegative-definite kernel  and is unique up to a unitary transformation.
The notion of nonnegative definite maps has a natural generalisation to maps on $S\times S$
taking values in a $C^*$-algebra; there is also an obvious counterpart of the notion
of Kolmogorov construction for such maps.

Quantum stochastic
is usually abbreviated to QS, operator space to OS,
completely positive to CP and completely positive and contractive
to CPC.

\section{Operator spaces}   \label{Ospaces}
\setcounter{equation}{0}

For general introduction to the theory of operator spaces we refer to \cite{ERuan};
for a functional-analytic point of view on this theory, with numerous applications to
operator algebras, we recommend \cite{Pisier}. Let us recall that
an (abstract) operator space $\Vil$ is a Banach space equipped with norms on matrix spaces
over $\Vil$ ($M_n(\Vil)$, $n \in \bn$), satisfying so-called Ruan's axioms.
If $\Vil$, $\Wil$ are operator spaces and $\phi:\Vil\to \Wil$ is a linear map, $\phi$ is
called \emph{completely bounded} if
\[ \|\phi\|_{cb} := \sup_{n\in \bn} \|\phi^{(n)}\| < \infty,\]
 where $\phi^{(n)}: M_n(\Vil) \to M_n(\Wil)$
is an obvious matricial extension of $\phi$ - $\phi$ is applied to each entry of the matrix.
Analogously we define completely contractive and completely isometric maps.
The space of completely bounded maps between operator spaces $\Vil$ and $\Wil$ is itself
an operator space, denoted further by $CB(\Vil; \Wil)$ - the matrix norms are introduced
with the help of the algebraic identification $M_n (CB(\Vil; \Wil)) \cong
CB(\Vil;M_n(\Wil))$. In particular, as every bounded linear functional is completely bounded,
we obtain an operator space structure on $\Vil^*= CB(\Vil;\bc)$.

For this work it will usually  be sufficient
to work with concrete operator spaces,
i.e.\ closed subspaces of the space $B(\hil_1;\hil_2)$ for  Hilbert spaces
$\hil_1, \hil_2$, with matricial norms induced by the identification
$M_n(B(\hil_1;\hil_2))\cong B(\hil_1^{\oplus n}; \hil_2^{\oplus n})$.
Ruan's Theorem says that each abstract operator space is completely isometric to a concrete
operator space.

The spatial (or minimal) tensor product of operator spaces $\Vil \subset B(\hil_1)$ and
$\Wil \subset B(\hil_2)$ is the norm closure of $\Vil \odot \Wil$ in $B(\hil_1 \ot \hil_2)$.
It is denoted by $\Vil\otimes\Wil$, and in fact does not depend on concrete representations
of $\Vil$ and $\Wil$ - it is completely isometric to the norm closure of the canonical
image of $\Vil \odot \Wil$ in $CB(\Vil^*;\Wil)$. The important fact for us is that
when $\Vil_1$, $\Vil_2$, $\Wil_1$, $\Wil_2$ are operator spaces, $\phi_1\in CB(\Vil_1;\Wil_1)$
and $\phi_2\in CB(\Vil_2;\Wil_2)$, then the map $\phi_1 \odot \phi_2:\Vil_1 \odot \Vil_2 \to
\Wil_1 \odot \Wil_2$  has a (unique) completely bounded extension
 $\phi_1 \ot \phi_2:\Vil_1 \ot \Vil_2 \to \Wil_1 \ot \Wil_2$, satisfying
 $\|\phi_1 \ot \phi_2\|_{cb} = \|\phi_1\|_{cb}  \|\phi_2\|_{cb}$.

Let $\alg, \blg$ be $C^*$-algebras (so in particular operator
spaces) and $\phi:\alg \to \blg$ be linear. The map $\phi$ is called
completely positive (CP) if for all $n\in \bn$ $\phi^{(n)}:M_n(\alg)
\to M_n (\blg)$ is positive. Every completely positive map
$\phi:\alg \to \blg$ is completely bounded and its complete bound
coincides with its operator bound; if $\alg$ is unital, $\|\phi\| =
\|\phi(I)\|$.

\subsection*{Matrix spaces}

 We need a concept of matrix spaces introduced by J.M.\,Lindsay and S.J.\,Wills in
\cite{lwex}.
Let $\Vil\subset B(\Kil)$ be an operator space and let $\hil$ be a supplementary Hilbert space.
The (concrete) operator space:
\[ M_{\hil} (\Vil) = \{ T \in B(\Kil \ot \hil): E^cTE_d \in \Vil \text{ for all } c,d \in \hil\}\]
is called an  \emph{$\hil$-matrix space} over $\Vil$. It is easy to see that $M_{\hil} (\Vil)$
contains the spatial tensor product $\Vil \ot B(\hil)$. If $\Vil$ is ultraweakly closed,
$M_{\hil} (\Vil)$ coincides with the ultraweak closure of $\Vil \odot B(\hil)$, denoted
by $\Vil \vnot B(\hil)$. If $\hil$ has finite dimension, $n$ say,
then $M_{\hil} (\Vil)\cong M_n(\Vil)$.
In general both inclusions
$\Vil \ot B(\hil) \subset  M_{\hil} (\Vil) \subset \Vil \vnot B(\hil)$ are proper.

Later the following functorial property of matrix spaces will be needed. If $\hil_1$,
$\hil_2$ are Hilbert spaces then there is a canonical complete isometry
\begin{equation}\label{functprop}
M_{\hil_1} \left(M_{\hil_2} (\Vil) \right)
\cong M_{\hil_2 \ot \hil_1} (\Vil). \end{equation}
The latter space is also completely isometric to $M_{\hil_1 \ot \hil_2} (\Vil)$,
with the complete isometry implemented by the tensor flip $\Kil \ot \hil_2 \ot \hil_1$
to $\Kil \ot \hil_1 \ot \hil_2$.

Whenever $\Wil \subset B(\Hil)$ is
another operator space, and $\phi \in CB(\Vil;\Wil)$, the map
$\phi \ot \id_{B(\hil)}$ extends uniquely to a completely bounded map
$\phi^{(\hil)}:  M_{\hil} (\Vil) \to   M_{\hil} (\Wil)$ satisfying
\[ E^c \big(\phi^{(\hil)}(T)\big) E_d = \phi (E^c T E_d),\]
for all $T \in  M_{\hil} (\Vil)$, $c,d \in \hil$. The map $\phi^{(\hil)}$ is called the
 \emph{$\hil$-lifting} of $\phi$ and satisfies
$\|\phi^{(\hil)}\|_{cb} = \|\phi\|_{cb}$.  Observe that when $\hil=\bc^n$ ($n \in \bn$), the map
$\phi^{(\hil)}$ is just the previously introduced lifting $\phi^{(n)}:M_n(\Vil) \to M_n(\Wil)$.

Analogously one can define \emph{$\hil$-column space} and  \emph{$\hil$-row space} over $\Vil$,
denoted respectively by $C_{\hil}(\Vil)$ and $R_{\hil} (\Vil)$ (so that
for example $C_{\hil}(\Vil)$ is a subspace of $B(\hil; \Kil \ot \hil)$). For $\Vil = \bc$
the notation $ \left|\hil\rra:= C_{\hil}(\bc)$, $\lla \hil \right|:=R_{\hil}(\bc)$ is also used.
Completely bounded maps lift to those as well, and the liftings are denoted respectively
by superscripts $^{(|\hil\rangle)}$, $^{(\langle\hil|)}$.

\section{Fock space notations and QS processes} \label{Focksection}
\setcounter{equation}{0}
In this section the definition of a symmetric Fock space is recalled and the terminology concerning
quantum stochastic processes established.

\subsection*{Symmetric Fock space}

\hspace*{0.5cm} Recall that the symmetric Fock space over a Hilbert space $\hil$ is defined as the orthogonal
Hilbert space sum $ \bigoplus_{n=0}^{\infty} (\hil^{\ot n}_{\mathrm{sym}})$, where $\hil^{\ot 0}_{\mathrm{sym}}:= \bc$,
and for each $n \in \bn$  the space $\hil^{\ot n}_{\mathrm {sym}}$ is a closed subspace
of $\hil^{\ot n}$ generated by vectors $\{u^{\ot n}: u \in \hil\}$.  It will be denoted by
$\Gamma(\hil)$. By the exponential property of Fock space is understood the existence
of the canonical
isomorphism $\Gamma(\hil_1\oplus \hil_2) \cong \Gamma(\hil_1) \ot \Gamma(\hil_2)$.

Let $\kil$ be a Hilbert space, called the \emph{noise dimension space}.
There is a natural isomorphism $L^2(\br_+) \ot \kil \cong L^2(\br_+;\kil)$, where the second
space is a space of square integrable functions on $\br^+$ with values in $\kil$. Integrability of Banach-space-valued functions is always understood as
Bochner integrability.
The stochastic arena for the action of QS stochastic cocycles, $\Gamma(L^2(\br_+;\kil))$,
will be denoted by $\Fock_{\kil}$ (or $\Fock$ if the space $\kil$ is clear from the context).
Whenever $J$ is a subinterval of $\br_+$ we will write $\Fock_{J; \kil}$ (or $\Fock_J$) for $\Gamma(L^2(J;\kil))$.
The exponential property of Fock space gives $\Fock$ a structure of a product system in the sense of Arveson
(\cite{Arv}, \cite{Raja}):
\begin{equation} \label{factor} \Fock = \Focktos \ot \Fockstot \ot \Fockafts, \;\;\; s,t \geq 0.\end{equation}

The exponential vectors in $\Fock$ are defined by
\[\ve(f):= \sum_{n=0}^{\infty} \frac{1}{\sqrt{n!}} f^{\ot n}, \;\;\;  f\in L^2(\br_+;\kil).\]
The terminology is justified by the following relation ($f,g\in L^2(\br_+;\kil)$):
\[ \la \ve(f), \ve(g) \ra = \exp (\la f,g \ra),\]
moreover the exponential vectors  behave well under the tensor decomposition of $\Fock$, given in (\ref{factor}):
\[ \ve(f) = \ve(f_{[0,s)}) \ot \ve (f_{[s,s+t)}) \ot \ve(f_{[s+t, \infty)}), \]
$ s,t \geq 0$, $f\in L^2(\br_+;\kil)$.
The vector $\ve(0)$ is denoted by $\Omega$ and called the \emph{vacuum vector};
the vector state $\omega_{\Omega}$ is called the \emph{vacuum state} on $B(\Fock)$. Observe that $\Fockstot$
is viewed as a subspace of $\Fock$ via the identification
\[ \ve(f) \longrightarrow \ve(0_{[0,s)}) \ot \ve(f) \ot \ve(0_{[s+t, \infty)}),\]
$f \in L^2([s,s+t); \kil)$.

For any dense subspace $D \subset \kil$ define
\[
\Step_D := \text{Lin}\{d_{[0,s)}: d \in D, s \in \br_+\}
\]
 ($\Step:= \Step_{\kil}$) and corresponding subspaces of $\Fock$:
\[
\Exps_D := \text{Lin} \{\ve(f): f \in \Step_D\}, \;\;\;  \Exps := \Exps_{\kil}.
\]

It may be shown
that $\Exps_D$ is a dense subspace of $\Fock$ (see for example \cite{lect}).
Elements of $\Exps_D$ will play the role of test functions for QS processes
(notice that our step functions are right-continuous).

The basic operators for the QS integration are the \emph{conservation} (known also as \emph{preservation}
or \emph{gauge}), \emph{creation} and \emph{annihilation} operators, defined on $\Exps$ respectively by
\begin{eqnarray*}
n(T) \ve(f)  &=& \frac{d}{dt} \ve(\textrm{e}^{tT}f) |_{t=0} \\
a^{\dagger} (g) \ve(f) & = & \frac{d}{dt} \ve(f + t g)|_{t=0} \\
a(g) \ve(f) & = & \la g, f \ra \ve(f),
\end{eqnarray*}
where $f,g \in \Step$ and $T \in B(L^2(\br_+;\kil))$.

Denote (for $t\geq 0$) by $s_t$ the shift isometry acting from $L^2(\br_+;\kil)$
to $L^2([t, \infty);\kil) $ defined by
\[ (s_t f )(s)  =  f (s-t),\;\; s \geq t.\]
This can be `second-quantised' to the isometry $S_t: \Fock \to  \Fockaft$, which is the continuous
extension of the map  defined on $\Exps$ by
\[ S_t(\ve(f)) = \ve(s_t f).\]

The \emph{CCR flow of index $\kil$} is the $E_0$-semigroup (\cite{Arv}) $\{\sigma_t:t \geq 0\}$ on $B(\Fock)$ defined by
\begin{equation} \label{shift} \sigma_t(X) = I_{\Focktot} \ot \left( S_t X S_t^* \right),\end{equation}
$X \in B(\Fock)$. Observe that the above formula may be extended to all operators
$X\in \Op (\Exps_D)$, in which case $\sigma_t(X)$ also belongs to $\Op (\Exps_D)$.

When $\kil_0$ is a closed subspace of a Hilbert space $\kil$, by the \emph{$\kil_0$-vacuum conditional
expectation} is understood a projection
$\condexp_{\Fock_{\kil_0}}:B(\Fock_{\kil}) \to B({\Fock_{\kil_0}})$ preserving the vacuum state. Its
action is given by 
\[T \mapsto E^{\Omega_0} T E_{\Omega_0},\]
 where $\Omega_0$ is the vacuum vector in $\Fock_{\kil_0^{\perp}}$.

It is often useful to view the symmetric Fock space from a different point, using
the description due to A.\,Guichardet (\cite{Guichardet}). For a subinterval $J\subset \br_+$ and $n \in \bn_0$
put
\[ \Gamma_J = \{ \sigma \subset J: \# \sigma < \infty\}, \text{ and }
\Gamma_J^{(n)}= \{\sigma \in \Gamma_J: \# \sigma = n \}.\]
Observe that (for $n \in \bn$) $\Gamma_J^{(n)}$ can be naturally identified with the set
$\{(t_1,\ldots,t_n) \in J^n: t_1 < \ldots < t_n\}$, and so is equipped with a measure induced by $n$-dimensional Lebesgue measure.
Using the fact that $\Gamma_J$ is a disjoint union
$\bigcup_{n=0}^{\infty} {\Gamma_J^{(n)}}$ and declaring  $\Gamma_J^{(0)}= \{\emptyset\}$
to be an atom of the unit measure allows us to introduce the \emph{Guichardet measure} (or the
\emph{symmetric measure})
on $\Gamma_J$ as the sum of the measures described above. Define
\[ \Guich_{J,\kil}:= \left\{F \in L^2 \left(\Gamma_J; \bigoplus_{n=0}^{\infty} {\kil^{\ot n}}\right)
: F(\sigma) \in \kil^{\ot \# \sigma} \text{ for a.e.\ } \sigma \right\},\]
where $\kil^{\ot 0}:= \bc$. For any $f\in L^2(J;\kil)$ a \emph{product function}
$\pi_f\in  \Guich_{J,\kil}$ is defined by
\[ \pi_f(\sigma) = \left\{\begin{matrix}
 f(s_n) \ot \cdots \ot f(s_1) & \text{ if } \sigma = \{s_1 <\ldots <s_n\} \\
  1 & \text{ if } \sigma = \{\emptyset\}. \end{matrix} \right. \]
The crucial fact here is that each product function belongs to $\Guich_{J,\kil}$ and the map
$\pi_f \to \ve(f)$ extends uniquely to a Hilbert space isomorphism $\Guich_{J,\kil} \to
\Fock_{J,\kil}$.

\subsection*{Processes}

\indent Let $\hil$ be an additional Hilbert space with a dense subspace $E$ and
fix a dense subspace $D \subset \kil$.
By an \emph{$\hil$-operator process} (with the domain $\Do:=E \odot \Exps_D$)
we understand a family
$X = \big(X_t\big)_{t \geq 0}$ of operators on $\hil \ot\Fock$,
each having the (dense) domain $\Do$, being weak-operator
measurable in $t$ and adapted to the natural Fock-space
operator-filtration.
Thus $X:\br_+ \to \Op (E \odot \Exps_D)$, $t\mapsto X_t\xi$ is weakly measurable
for all $\xi\in\Do$ and, for each $t\geq 0$, $\zeta \in E$,
$f \in \Step_D$,
\[X_t(\zeta \ot \ve(f))= X(t)\big(\zeta \ot \ve(f_{[0,t[})\big)\ot\ve(f_{[t,\infty[})\]
for some operator $X(t)\in\Op\big(E \odot \Exps_{D,[0,t)}\big)$, where
$\Step_{D,[0,t)}$ is defined as $\Step$ is, except that $\br_+$ is
replaced by $[0,t)$.
Two $\hil$-operator processes $X$, $X'$ with the domain $\Do$ are called
\emph{indistinguishable}
if for all $\xi \in \Do$, $\zeta \in \hil \ot \Fock$,
\[ \la \zeta, X_t \xi \ra = \la \zeta, X'_t \xi \ra \;\;\; \text{ for almost all } t\in \br_+.\]

Identifying indistinguishable processes allows us to consider the linear space of all
$\hil$-operator processes with the domain $E \odot \Exps_D$, denoted $\Proc (E, \Exps_D)$
(or  $\Proc (\Exps_D)$ if $\hil=\bc$).

A process $X\in \Proc (E, \Exps_D)$
is called \emph{weakly regular} if for each $f,g \in \Step_D$, $t \geq 0$
the operator $E^{\ve(f)} X_t E_{\ve(g)} :E \to \hil$ is bounded (recall the notation
(\ref{Dirac})), and, for fixed $f,g$, the respective norms are locally bounded with respect
to $t$.
It is called \emph{weakly continuous} if for each $\zeta, \xi \in \Do$
the scalar-valued function $ t \to \la \zeta, X_t \xi \ra$ is continuous.
Further $X$ as above is called \emph{square integrable} if for each $\xi \in \Do$
the vector-valued function  $t \to X_t \xi$ is locally square integrable; it is called
\emph{locally bounded} if each $t \to X_t \xi$ is locally bounded, \emph{continuous} if
each $t \to X_t \xi$ is continuous and finally \emph{H\"older continuous} if $t \to X_t \xi$ is
H\"older continuous at $0$ with exponent $\frac{1}{2}$:
\begin{equation} \limsup_{t \to 0^+} {t^{-\frac{1}{2}} \|X_t \xi - X_0\xi\|} < \infty.
\label{cHolder} \end{equation}
It is called \emph{bounded} if $X_t$ is a bounded operator for each $t\geq 0$. In such case
we will denote the continuous extension of $X_t$ to an operator in $B(\hil \ot \Fock)$ by the same
letter.

The respective subspaces of $\Proc (E, \Exps_D)$ will be denoted by
$\Procwr (E, \Exps_D)$ (weakly regular),  $\Procwc (E, \Exps_D)$
(weakly continuous), $\Proctwo (E, \Exps_D)$ (square integrable),  $\Proclb (E, \Exps_D)$
(locally bounded),  $\Procc (E, \Exps_D)$
(continuous), $\ProcHc (E, \Exps_D)$ (H\"older continuous) and $\Procb (E, \Exps_D)$
(bounded).  Note that
\begin{equation} \label{wc=wr} \Procwc (\Exps_D) \subset \Procwr (\Exps_D).
\end{equation}

Additionally let
\[ \Proc^{\dagger} (\hil, \Exps_D) := \{X \in \Proc (\hil, \Exps_D): \forall_{t \in 0} \;
 X_t \in \Opdag (\hil \odot \Exps_D) \}, \]
\[ \Proc_{\alpha}^{\dagger} (\hil, \Exps_D) := \{X \in \Proc^{\dagger} (\hil, \Exps_D):
 X, X^{\dagger} \in \Proc_{\alpha} (\hil, \Exps_D) \}, \]
 where $\alpha$ may be any of the above subscripts, and the involution $X \to X^{\dagger}$
 is induced pointwise from the involution (\ref{daggerinv}).

Let now $V$ be a vector space and $\Wil \subset B(\hil)$ an operator space.
A linear map $k$ from $V$ to
$\Proc (\hil,\Exps_D)$ is called a
\emph{process
on $V$ with values in} $\Wil$ (and domain $\hil \odot \Exps_D$)
if
for each $f,g\in \Step_D$, $t \geq 0$, $v \in V$ operator
$E^{\ve(f)} k_t(v) E_{\ve(g)}$ belongs to $\Wil$.
The vector space of all processes on $V$ with values in $\Wil$ is written
$\Proc(V;\Wil, \Exps_D)$ (or simply $\Proc(V; \Exps_D)$ if   $\Wil=\bc$).
We say that $k \in \Proc(V;\Wil, \Exps_D)$ is \emph{weakly regular} (respectively
weakly continuous, square integrable, etc.) if
each $k(v)$ ($v\in V$) satisfies the required property; similarly
$\Proc_{\alpha}^{\dagger}(V;\Wil, \Exps_D)$ denotes the space of all
$k\in \Proc(V;\Wil, \Exps_D)$ such that each $k(v)$ ($v \in V$) belongs to
$\Proc_{\alpha}^{\dagger}(\hil, \Exps_D)$ and for each $f,g\in \Step_D$, $t \geq 0$, $v \in V$ operator
$E^{\ve(f)} k_t(v)^{\dagger} E_{\ve(g)}$ belongs to $\Wil$ (again $\alpha$ is any of the available subscripts).
If $V$ is a subspace of some involutive vector space, each $k \in \Proc_{\alpha}^{\dagger}(V;\Wil, \Exps_D)$
determines a process $k^{\dagger} \in \Proc_{\alpha}^{\dagger}(V^{\dagger};\Wil^*, \Exps_D)$ by the
standard formula
$k^{\dagger}(v^{\dagger}):=(k(v))^{\dagger}$.
If $V$ is closed under the involution $\dagger$, $W$ is selfadjoint and $k=k^{\dagger}$, $k$ is called \emph{real}.

If $\Vil$ is an operator space, $k \in \Proc(\Vil; \Wil, \Exps_D)$ is called
\emph{completely bounded} if
for each $t\geq 0$, $v \in \Vil$ the process $k(v)$ is bounded and the
resulting map
$k_t:\Vil \to B(\hil\ot \Fock)$ is completely bounded. It is called
\emph{weakly bounded} if for each $t\geq 0$ $f,g \in \Step_D$ the map
\[ \Vil \ni v \mapsto E^{\ve(f)} k_t(v) E_{\ve(g)} \in \Wil\]
is bounded, and \emph{uniformly weakly bounded} if it is weakly bounded and for all
$f,g \in \Step_D$, $T >0$
\[ \sup_{t<T} \|E^{\ve(f)} k_t(\cdot) E_{\ve(g)}\| < \infty.\]
Observe that due to Banach-Steinhaus Theorem processes which are weakly regular and weakly bounded are automatically
uniformly weakly bounded.

Let $k\in \Proc(\Vil;\Wil, \Exps)$. If for all $t \geq 0$, $g\in \Step$ and $v \in \Vil$ the operator
 $k_t (v) E_{\ve(g)}$ is in $C_{\Fock}(\Wil)$ (equivalently, is bounded) and the map
\begin{equation}
\label{sreg} \Vil \ni v \mapsto k_t(v) E_{\ve(g)}\in C_{\Fock}(\Wil)\end{equation}
(denoted further by $k_{t,\ve(g)}$) is bounded then $k$ is said to have \emph{bounded columns}. If $k$ has bounded columns and for all $g\in \Step$, $T \geq 0$
\begin{equation}
\label{sregest} \sup_{t<T} \| k_{t,\ve(g)}\| < \infty,\end{equation}
then $k$ is called \emph{strongly regular}. If the map $k_{t,\ve(g)}$
 is completely bounded (i.e.\ $k$ has \emph{completely bounded columns})
and the estimate \eqref{sregest} holds for $\cb$-norms, $k$ is called \emph{$\cb$-strongly regular}.
Similarly by analogy with \eqref{cHolder} we can define \emph{$\cb$-H\"older continuous processes}.

The families
of all processes satisfying conditions introduced above are denoted respectively by $\Proc_{\cb}(\Vil; \Wil, \Exps)$ (completely
bounded), $\Proc_{\wb}(\Vil; \Wil, \Exps_D)$ (weakly bounded), $\Proc_{\uwb}(\Vil; \Wil, \Exps_D)$ (uniformly weakly bounded),
$\Proc_{\sr}(\Vil; \Wil, \Exps_D)$ (strongly regular) and $\Proc_{\cbsr}(\Vil; \Wil, \Exps_D)$ ($\cb$-strongly regular).
\emph{Completely positive, unital and contractive} processes on $\alg$ with values
in $\blg$ (where $\alg$, $\blg$ are $C^*$-algebras) are defined by analogy with completely bounded processes.

\section{Quantum stochastic integrals} \label{QSIntegrals}
\setcounter{equation}{0}

Fix again a dense subspace $D$ of the noise dimension space $\kil$, and let $E$ be a dense
subspace of another Hilbert space $\hil$. Recall the notation \eqref{hats}.

Quantum stochastic integration introduced by R.L.\,Hudson and K.R.\,Partha-sarathy gives a linear map
\[
\Proctwo (E\odot \Dhat,\Exps_D ) \rightarrow \Procc (E, \Exps_D),
\]
denoted
\[ X \rightarrow \int^{\cdot}_0 X_s \, d\Lambda_s.\]
In the matrix form, if
\[ X = \begin{bmatrix}
K & G \\ F & H
\end{bmatrix},\]
for certain processes $K, G, F$ and $H$, the QS integral is given by:
\[ \int^{\cdot}_0 X_s \, d\Lambda_s = T(K) + A^* (F) + N(H) + A(G),\]
where $T$ denotes the \emph{time integral}, $A^*$ the \emph{creation integral},
$A$ the \emph{annihilation integral} and $N$ the \emph{preservation integral}.
The latter integrals can be defined either by suitable limits of integrals of
piecewise constant processes
with respect to integrators given in terms of the basic operators on Fock space
introduced in the previous section (\cite{par}, \cite{Hud}), or by using the
gradient operator of Malliavin calculus and Hitsuda-Skorohod integral (\cite{lect}).
In the context of this work, where we are concerned with adapted, square-integrable
processes defined on exponential domains, both definitions are equivalent.
An advantage of the second (Malliavin-style) technique lies in its explicit
independence on the choice of an orthonormal basis in $\kil$.
 The square integrability assumptions on processes in the entries of the matrix above
 guranteeing the existence of the QS integral may be weakened (\cite{lect}).

Quantum stochastic integral
enjoys the \emph{Fundamental Formulae} and the \emph{Fundamental Estimate} below.

\begin{tw} [First Fundamental Formula]
Let $X \in \Proctwo (E\odot \Dhat,\Exps_D )$, $\xi,\eta \in E$, $f,g\in\Step_D$ and
$ t \geq 0$. Then
\begin{equation}  \label{FFF}
\lla \xi\ve (f) ,\left(\int_0^t X_s d \Lambda_s \right)\eta\ve (g) \rra =
\int_0^t  \lla \xi \fhat(s) \ve(f), X_s \left(\eta \ghat(s) \ve(g)\right) \rra  ds.
\end{equation}
\end{tw}

As $E \odot\Exps_D$ is dense in $\hil \ot\Fock$, the First Fundamental
Formula determines $\int^{\cdot}_0 X_s \, d\Lambda_s$  uniquely. It also implies that
for $X \in \Proctwo^{\dagger} (E\odot \Dhat,\Exps_D )$, $t\geq 0$,
\begin{equation} \left( \int_0^t X_s d \Lambda_s \right)^{\dagger} = \int_0^t X^{\dagger}_s d \Lambda_s.
\label{dagQS} \end{equation}

The next estimate
is very important for constructing solutions of QS differential equations.
We use a natural notation $\int_{r}^t X_s d \Lambda_s = \int_{0}^t X_s d \Lambda_s -
\int_{0}^r X_s d \Lambda_s$ ($0\leq r \leq t$).

\begin{tw} [Fundamental Estimate]
Let $X \in \Proctwo (E\odot \Dhat,\Exps_D )$, $\xi\in E$, $f\in\Step_D$ and
$0 \leq  r \leq t  \leq T$. Then
\begin{equation}  \label{FE}
\left\| \left(\int_{r}^t X_s d \Lambda_s \right)\xi\ve (g) \right\|^2
\leq  C(f,T) \int^t_r \| X_s  \left(\xi \fhat(s) \ve(f) \right)\|^2 \, ds,
\end{equation}
where $C(f,T)>0$ is a constant dependent only on $f$ and $T$.
\end{tw}

Finally the Second Fundamental Formula is a QS extension of the It\^o Formula.

\begin{tw} [Second Fundamental Formula]
Let $X,Y \in \Proctwo (E\odot \Dhat,\Exps_D )$, $\xi,\eta \in E$, $f,g\in\Step_D$,
$t \geq 0$ and let $L = \int^{\cdot}_0 X_s \, d\Lambda_s$,
$M = \int^{\cdot}_0 Y_s \, d\Lambda_s$.  Then
\begin{align}  \label{SFF}
\lla L_t \left(\xi \ve (f)\right), M_t \left( \eta \ve (g) \right)\rra =
 \int_0^t  ds &\left(\lla \wt{L_s} \left(\xi \fhat(s)\ve (f)\right) ,
   Y_s  \left(\eta\ghat(s) \ve (g) \right)\rra     \right.    \\
   & + \lla X_s \left(\xi \fhat(s)\ve (f)\right), \wt{M_s}
         \left(\eta\ghat(s) \ve (g) \right)\rra \notag \\
& \left. +\lla X_s  \left(\xi \fhat(s)\ve (f)\right),
\QSproj_{\hil,\Fock} \,  Y_s   \left(\eta\ghat(s) \ve (g)\right) \rra \right) \notag,
\end{align}
where $\wt{L_t}$, $\wt{M_t}\in \Op(E\odot \Dhat \odot \Exps_D)$ are given
respectively by $\wt{L_t}:= \tau \big(I_{\Dhat} \ot L_s \big)$,
 $\wt{M_t}:= \tau \big(I_{\Dhat} \ot M_s \big)$, and $\tau$ is the
tensor flip: $\tau: \Op(\Dhat \odot E \odot \Exps_D) \to  \Op(E\odot \Dhat \odot \Exps_D)$.
\end{tw}

\subsection*{Iterated QS integration}

For $L \in \Op (E \odot \Dhat^{\odot n})$, where $E$ is a dense subspace
of a Hilbert space $\hil$, define
$\Lambda^n (L) \in \Procc (E, \Exps_D)$ recursively as follows (we write
$\Lambda^n_t(L)$ for $\Lambda^n(L)_t$):
\[
\Lambda^0_t (L) = L \ot I _{\Exps_D},
\text{ and, for } n\geq 1, \
\Lambda^n_t (L) = \int^t_0
\Lambda^{n-1}_s (L) \, d\Lambda_s,
\]
by viewing $E \odot \Dhat^{\odot n}$ as $(E\odot \Dhat ) \odot \Dhat^{\odot(n-1)}$.

Letting $\xi$, $\eta$, $f$, $g$ and $T$ be as in the Fundamental Formulae
above, we obtain
\begin{equation}    \label{iterated FF1}
\la \xi\ve (f) ,\Lambda^n_t (L) \eta\ve (g)\ra =
\int_{\Delta_n[0,t]}
\lla \zeta (\bm{s}),
(L \ot I_{\Fock}) \theta (\bm{s})\rra \, d\bm{s},
\end{equation}
and
\begin{equation} \label{iterated FE}
\| \Lambda^n_t (L) \xi\ve (f) \|^2  \leq C(f,T)^n
\int_{\Delta_n[0,t]}
\| (L \ot I_{\Fock}) \zeta (\bm{s}) \|^2 \, d\bm{s},
\end{equation}
where
$\zeta (\bm{s}) := \xi \fhat^{\ot n} (\bm{s}) \ve (f)$
 ($\fhat^{\ot n} (\bm{s}) := \fhat (s_n) \cdots \fhat (s_1)$)
and similarly  $\theta(\bm{s}) := \eta \ghat^{\ot n} (\bm{s}) \ve (g)$.

Further we would like to describe the action of the sequences of iterated QS integrals.
To this end we introduce spaces of so-called \emph{matrix-sum kernels}. Let $E$ be again a dense
subspace of a Hilbert space $\hil$ and let
\[ \Seq_{E,D}:= \{F=(F_n)_{n=0}^{\infty}\,: \;\forall_{n \in \bn_0}\ F_n \in \Op (E \odot \Dhat^{\odot n})\},\]
\[ \Seqdagger_{E,D} :
= \big\{ F \in \Seq_{E,D} \,: \; \forall_{n \in \bn_0}\, F_n \in \Opdag (E \odot \Dhat^{\odot n}) \big\}, \]
(where as usual $\Dhat^{\odot 0}$ denotes $\bc$). We also write $\Seq_D:=\Seq_{\bc,D}$, $\Seqdagger_D =
\Seqdagger_{\bc,D}$.
It is obvious that $\Seqdagger_{E,D}$ is a subspace of a vector space $\Seq_{E,D}$. Moreover
$\Seqdagger_{E,D}$ is an  involutive vector space with the involution $\dagger$
induced pointwise from the involution \eqref{daggerinv}.

For QS purposes we need to impose a growth condition. Let $\Geom_{E,D}$ denote the vector
space of all $F \in \Seq_{E,D}$ satisfying
\begin{equation} \label{GsubD}
\forall_{\xi \in E, \,S\subset \subset \Dhat}\, \exists_{C_1,C_2 >0}\,
\forall_{n\in \bn_0, \chi_1 , \ldots ,\chi_n \in S}\,
\| F_n (\xi \ot \chi_1 \ot \cdots \ot \chi_n) \| \leq C_1 C^n_2 ,
\end{equation}
and define the subspace
\begin{equation} \label{Gsuper}
\Geomdagger_{E,D} :
= \left\{ F \in \Seqdagger_{E,D} \, \big| \; F,F^\dagger \in \Geom_{E,D} \right\},\end{equation}
again abbreviating to $\Geom_D$, $\Geomdagger_D$ when $E=\bc$.
The estimate \eqref{iterated FE} implies that if $F \in \Geom_{E,D}$, $\xi \in E$, $f \in \Step_D$,
then $\sum_{n \geq 0} \Lambda^n_t (F_n) \xi\ve (f)$ is absolutely convergent, 
and the convergence is locally uniform in $t$. The resulting map
\begin{equation} \label{Lambda map}
\Lambda : \Geom_{E,D} \rightarrow \Procc (E,\Exps_D)
\end{equation}
is linear. In the Guichardet notation (cf.\ Section \ref{Focksection})
\eqref{iterated FF1} yields the identity
\begin{equation} \label{Guichardet FF1}
\left\la \xi \ve (f), \Lambda_t (F) \eta \ve (g)\right\ra = \int_{\Gamma_{[0,t[}} \,
d\sigma\left\la \xi\pi_{\fhat} (\sigma), F_{\# \sigma} \eta\pi_{\ghat} (\sigma)\right\ra \,
\la \ve (f), \ve (g)\ra ,
\end{equation}
where $\xi, \eta \in E$, $f,g \in \Step_D$ and  $\int \, d \sigma$ denotes integration with respect to the Guichardet measure.
Another application of \eqref{iterated FE} shows that indeed $\Lambda(\Geom_D) \subset \ProcHc(\Exps_D)$.
Moreover \eqref{dagQS} implies that $\Lambda(\Geomdagger_D) \subset \ProcHc^{\dagger}(\Exps_D)$ and
for $F \in \Geomdagger_{E,D}$
\begin{equation} \label{Lambdadag} \Lambda(F)^{\dagger} = \Lambda(F^{\dagger}).\end{equation}
A minor modification of the proof of Proposition 2.3 of \cite{lwhom} yields injectivity
of the map $\Lambda$.

The formulation of the It\^o Formula for iterated QS integrals is postponed to respectively
Sections \ref{multalg} and \ref{multCalg} - it is more convenient to express it separately in
the purely algebraic and in the operator space theoretic context.

\chapter{Algebraic case}
\newpage

This chapter is concerned with quantum stochastic convolution
cocycles on a coalgebra. They are a natural, still purely algebraic,
generalisation of quantum L\'evy processes on a $^*$-bialgebra
defined by L.\,Accardi, W.\,von Waldenfels and M.Sch\"urmann in
\cite{asw}.

QS convolution cocycles are obtained by solving coalgebraic QS
differential equations. The original proof of the existence of
solutions for such equations (for a particular type of the
generator, guaranteeing that the solution is unital and
$^*$-homomorphic), due to  M.\,Sch\"urmann (\cite{schu}), exploited
the formulation of QS integrals in terms of integral-sum kernels.
Here the proof is
 simplified by using the language of iterated QS integration of matrix-sum kernels,
whose use has been promoted in \cite{lwex} and in \cite{lwhom}.
The uniqueness of weak solutions is also proved. The basic tool is the Fundamental
Theorem on Coalgebras, allowing essentially to reduce the considered problem to
the finite dimensional case.

The cocycles arising in this way satisfy a H\"{o}lder condition, and
it is shown that conversely every such H\"{o}lder-continuous cocycle
is governed by a QS differential equation. Algebraic structure
enjoyed by matrix-sum kernels yields a unital $^*$-algebra of
processes which allows easy deduction of homomorphic properties of
cocycles on a $^*$-bialgebra, the stochastic generators of which can
be described in terms of so-called Sch\"urmann triples. This in
particular yields a simple proof of the Sch\"urmann Reconstruction
Theorem - every quantum L\'{e}vy process may be equivalently
realised in Fock space. The perturbation of QS convolution cocycles
by Weyl cocycles is shown to be implemented on the level of their
stochastic generators by the action of the corresponding Euclidean
group on Sch\"{u}rmann triples.

There is also a corresponding concept of opposite QS convolution cocycles that
can be obtained from the usual ones either by the time-inversion or by using the
opposite comultiplication on a coalgebra. Opposite cocycles also satisfy QS differential
equations.

\section{Coalgebras and convolution semigroups} \label{algandconv}
\setcounter{equation}{0}

In this section we present the definition of a coalgebra and establish basic facts concerning
the convolution product provided by the coalgebraic structure.

\begin{deft}
A vector space $\Coalg$ is a coalgebra if there are linear
maps
 $\Com : \Coalg \rightarrow \Coalg \odot \Coalg$ and
$\Cou : \Coalg \rightarrow \bc$, called the \emph{coproduct} and
\emph{counit} respectively, enjoying the coassociativity and the counit
property, namely
\begin{align}
(\id \ot \Com ) \circ \Com &=
( \Com \ot \id ) \circ \Com,
 \label{coassoc} \\
(\id \ot \Cou )\circ \Com &= ( \Cou \ot \id ) \circ
\Com = \id.  \label{counit}
\end{align}
\end{deft}

The main examples of coalgebras of interest derive from the theory of quantum groups.
Sometimes it will be handy to think of $\bc$ as a coalgebra, with both
coproduct and counit given by the identity mapping.

M.E.\,Sweedler has introduced the notation $a_{(1)} \ot a_{(2)}$ for
$\Com a$, in which
both summation and indices are supressed (\cite{Sweedler}). With this,
\eqref{coassoc} and \eqref{counit} read
\[
a_{(1)} \ot a_{(2)(1)} \ot a_{(2)(2)} = a_{(1)(1)} \ot a_{(1)(2)} \ot a_{(2)},
\quad \text{and} \quad
a_{(1)} \Cou (a_{(2)}) = \Cou (a_{(1)}) a_{(2)} = a.
\]
Let $\Com_0 := \id_{\Coalg}$ and for $n \in \bn$ define
\begin{equation} \label{Com n}
\Com_n := \big( \id^{\ot (n-1)} \ot \Com \big)
\circ \cdots \circ (\id \ot \Com ) \circ \Com .
\end{equation}
The coassociativity implies that moving any $\Com$ to any
of the available tensor places within its bracket
(rather than the right-most, as here) has
no effect. It is easily verified that the family
$\{ \Com_n : n \in \bn_0 \}$ satisfies
\begin{equation} \label{Com ij}
\big( \Com_i \ot \Com_j \big) \circ \Com =
\Com_{i+j+1}.
\end{equation}
The Sweedler notation  extends to writing
$a_{(1)} \ot \cdots \ot a_{(n+1)}$ for $\Com_n a$ $(n \geq 1)$, $a \in\Coalg$.
Thus, for example, $a_{(1)} \ot a_{(2)} \ot a_{(3)}$ becomes a neutral
notation for the effect of \eqref{coassoc} on an element $a$.

The following important result holds (the proof may be found for example in
\cite{schu}):

\begin{tw} [Fundamental Theorem on Coalgebras]  \label{FTOC}
Let $\Coalg$ be a coalgebra, $a\in \Coalg$. The subcoalgebra of $\Coalg$
generated by $a$ is finite-dimensional.
\end{tw}

Let $U,V,W$ be vector spaces and assume there is a natural map $\cdot:( U \odot V) \to W$
(for example $V=\bc$, $U=W$).
 For linear maps $\alpha : \Coalg \rightarrow U$, $\beta : \Coalg \rightarrow
V$ define
\begin{equation} \label{convolution}
\alpha \star \beta := \cdot(\alpha \ot \beta) \Com :
\Coalg \rightarrow W.
\end{equation}

Thus, for example,  the counit property~\eqref{counit} implies that
$\Cou \star \alpha = \alpha \star \Cou = \alpha$ for any
linear map
$\alpha$ from $\Coalg$ into a vector space. In particular
$\big( L(\Coalg ; \bc ) , \star \big)$ is a unital algebra with
identity
$\Cou$.

A functional $\lambda:\Coalg \to \bc$ is called \emph{idempotent} if
$\lambda = \lambda \star \lambda$ (so for example the counit is idempotent).

\begin{deft}
A convolution semigroup of functionals
(CSF, for short) on a coalgebra $\Coalg$ is a family $\{\kappa_s:s \geq 0\}$ of
linear functionals on $\Coalg$ satisfying the following conditions::
\begin{equation} \label{CSFincr}
\kappa_{s+t} = \kappa_s \star \kappa_t , \;\;\; s,t \geq 0, \end{equation}
\begin{equation} \label{CSFinit}
\kappa_0 = \Cou.\end{equation}
If for each $a \in \Coalg$
\[ \lim_{t \to 0^+}\kappa_t (a) = \Cou (a)\]
then $\{\kappa_s:s \geq 0\}$ is said to be continuous (CCSF).
\end{deft}

The following fact appeared first in \cite{asw}, for another proof see \cite{schu}. It is
a straightforward consequence of the Fundamental Theorem on Coalgebras.

\begin{fact}
Each CCSF $\{\kappa_s:s \geq 0\}$
has a \emph{generator} $\gamma\in L(\Coalg; \bc)$ given by
\[
\gamma(a) = \lim_{t \rightarrow 0^+} \frac{\kappa_t (a) - \Cou (a)}{t}, \;\;\; a \in \Coalg.\]
The CCSF may be recovered from its generator by the formula
\begin{equation} \label{exp star}
\kappa_t (a)= \exp_{\star} t\gamma (a):=\sum_{n \geq 0} (n!)^{-1} t^n \gamma^{\star n} (a)
\end{equation}
($a\in \Coalg$) where $\gamma^{\star 0} := \Cou$.
\end{fact}

\begin{deft} \label{modcoalg}
 A coalgebra $\Coalg$ is involutive if there exists
 a conjugate linear involution $^*$ compatible with the coalgebra operations:
$\Cou (a^*)= \overline{\Cou(a)}$, $\Com (a^*) = (a_{(1)})^* \ot
(a_{(2)})^*$.  An algebra $\Alg$ is a bialgebra if it is also a
coalgebra with multiplicative coproduct and counit. It is unital if
it is unital as an algebra and $\Cou (1) =1$, $\Com (1) = 1 \ot 1$.
Finally  a $^*$-bialgebra is a bialgebra with involution which is
both algebraic and coalgebraic.
\end{deft}

An important example of an idempotent functional on a unital
$^*$-bialgebra $\Alg$ is so-called \emph{Haar state},  i.e.\ a state
(positive unital functional) $h: \Alg \to \bc$ satisfying the (left
and right) invariance properties:
\begin{equation}  h(a_{(1)}) a_{(2)}  = a_{(1)} h(a_{(2)}) = h(a) 1 \label{Haar} \end{equation}
for all $a \in \Alg$. It is easy to see that the Haar state on $\Alg$ is unique (if it exists).

Note the following traffic between properties of a CCSF
$(\kappa_t)_{t\geq0}$ and its generator $\gamma$ when the coalgebra
has more structure. The semigroup consists of real (respectively,
unital) functionals if and only if the generator is real
(resp.\,vanishes at the unit). Moreover, on a $^*$-bialgebra $\Alg$,
if the functionals are positive then the generator is
\emph{conditionally positive}:
\begin{equation} \label{conditionally positive}
\gamma (a) \geq 0 \text{ for } a \in \Alg_+ \cap \Ker \Cou,
\end{equation}
as for all $a\in \Alg_+ \cap \Ker \Cou$, $t \geq 0$,
\[
\frac{\kappa_t (a) - \Cou (a)}{t} = \frac{\kappa_t (a)}{t} \geq 0.\]

\section{QS convolution cocycles} \label{QSccocycl}
\setcounter{equation}{0}

Here we introduce for the first time the main object analysed in the thesis.

\begin{deft}
A quantum stochastic convolution cocycle
(on $\Coalg$ with domain $\Exps_D$) is a process
$l\in\Proc (\Coalg; \Exps_D)$ such that, for all $s,t \geq 0$,
\begin{equation}
\label{convinc}
l_{s+t} = l_s \star (\sigma_s \circ l_t)
\end{equation}
 and for all $a\in \Coalg$
\begin{equation} \label{couinit}
l_0 (a) = \Cou(a) I_{\Exps_D}.
\end{equation}
\end{deft}

Family of all QS convolution cocycles is denoted by $\QSC (\Coalg; \Exps_D)$; the
notation adorned with subscripts and superscript $^\dagger$ according to our convention.

\begin{rem} \label{cocprop}
The formula \eqref{convinc} makes use of the assumed adaptedness of $l$ and the identification
$\Focktots \cong \Focktos \ot \Fockstot$.
It is referred to as the
\emph{convolution increment property}, and any process satisfying it is
called a \emph{convolution increment process}. The initial condition for a
convolution increment process $k\in \Proc(\Coalg; \Exps_D)$ must have the form
$k_0(a) = \lambda(a) I_{\Exps_D}$ ($a \in \Coalg$), where $\lambda:\Coalg\to \bc$ is
an idempotent functional.
\end{rem}

The convolution increment property has an equivalent `weak' description:

\begin{lemma} \label{M}
Let $l \in \Proc (\Coalg ; \Exps_D)$. Then $l$  is
a convolution increment process if and only if for each $f,g \in \Step_D$,
\begin{equation} \label{semigroup decomp}
\mathrm{e}^{-\la f_{[0,t)}, g_{[0,t)} \ra } \lla \ve (f_{[0,t)}), l_t (a)
\ve (g_{[0,t)}) \rra = \prod^{n-1}_{i=0} \lambda^{f(t_i),g(t_i)}_{t_{i+1}-t_i}
(a_{(i+1)})
\end{equation}
where $0 = t_0 \leq t_1 \leq \cdots \leq t_n =t$ contains the discontinuities
of $f_{[0,t)}$ and $g_{[0,t)}$,
$a_{(1)} \ot \cdots \ot a_{(n)} = \Com_{n-1} (a)$
(Sweedler-style), and
\begin{equation} \label{assoc semigroups}
\lambda^{c,d}_t := \mathrm{e}^{-t\la c,d\ra} \lla \ve (c_{[0,t[}), l_t (\,\cdot\,)
\ve (d_{[0,t[} ) \rra.
\end{equation}
\end{lemma}

\begin{proof}
The identity \eqref{semigroup decomp} results from repeated application
of the cocycle relation \eqref{convinc}. The other direction is trivial.
\end{proof}

\begin{cor} \label{N}
Let $l \in \QSCwc (\Coalg ; \Exps_D)$, $c,d \in D$. Then
\eqref{assoc semigroups} defines a continuous convolution semigroup of functionals
on $\Coalg$ : $ (\lambda^{c,d}_t )_{t \geq 0}$.
\end{cor}

\begin{deft}    \label{defsemig}
Semigroups defined by $\eqref{assoc semigroups}$ (for all pairs $c,d \in D$) are called
\emph{associated convolution semigroups} of the cocycle $l$.
\end{deft}

It is clear that two QS convolution cocycles with identical convolution semigroups are equal.
The above `semigroup decomposition' of a QS convolution cocycle is crucial for the results in Section
\ref{stochsection}.

Observe that if $\Coalg$ is a trivial coalgebra $\bc$, the QS convolution cocycle $l$
on $\bc$ is  in fact determined by the operator process $l(1)$, and the latter is
an \emph{operator Markovian cocycle}:

\begin{deft}
A process $X\in \Proc(\Exps_D)$ is called an operator Markovian cocycle if $X_0 = I_{\Exps_D}$ and
\[ X_{s+t} = X_s \sigma_s (X_t), \;\;\; s,t \geq 0.\]
\end{deft}

\section{Coalgebraic QS differential equations}    \label{coalgQSDE}
\setcounter{equation}{0}

In this section we describe the conditions assuring the existence and uniqueness
of solutions of coalgebraic QS differential equations and present basic properties of the
solutions.

Let $\varphi \in L\big(\Coalg ; \Op (\Dhat )\big)$. A
\emph{coalgebraic QS differential equation}
(with the coefficient $\varphi$) is the equation of the form
\begin{equation} \label{qgQSDE}
dk_t = k_t \star_{\tau} \varphi \,d\Lambda_t , \; \;
k_0 = \iota\circ\epsilon
\end{equation}
($\tau$ denoting a tensor flip exchanging the order of $\khat$ and $\Fock$,
$\iota$ indicating an ampliation).
A process $k \in \Proc (\Coalg , \Exps_D)$ is a \emph{weak solution} of the
equation \eqref{qgQSDE}
if for all $f,g \in \Step_D$, $t \geq 0, a \in \Coalg$
\begin{equation} \label{wQSDEsweedler}
\left\la \ve (f),\big( k_t (a)-\Cou (a)I_{\Fock}\big) \ve(g)\right\ra
=
\int^t_0
\lla \fhat (s), \varphi (a_{(2)}) \ghat (s) \rra
\lla \ve (f) ,k_s (a_{(1)}) \ve (g)\rra \, ds.
\end{equation}
Thus a weak solution $k$ is necessarily weakly continuous,
and therefore also weakly regular by \eqref{wc=wr}.
If $k\in\Proc_2 (\Coalg; \Exps_D )$ then the First Fundamental Formula implies that
\[
k_t (a) =
\epsilon (a) I +
\int^t_0 \varphi(a_{(2)}) \ot k_s(a_{(1)}) \, d\Lambda_s, \]
and $k$ is called a \emph{strong solution} of \eqref{qgQSDE}.

\subsection*{Uniqueness of the solution}

Before establishing the uniqueness of the solution of \eqref{qgQSDE}
we need a following lemma:

\begin{lemma} \label{D}
Let $k \in \Procwr (\Coalg ; \Exps_D)$  and let $V$ be a
finite dimensional subspace of $\Coalg$ equipped with some norm
$\|\cdot\|$. Then, for
each $f,g \in \Step_D$ and $T \geq 0$,
\begin{equation} \label{constants}
C_{f,g,T,V} :=
\sup \Big\{  \big| \lla \ve (f), k_t (a) \ve (g)\rra\big|\,
: \;  a \in V, \| a\| \leq 1 ,\, 0 \leq t \leq T
\Big\} < \infty .
\end{equation}
\end{lemma}

\begin{proof}
Let $e_1, \ldots e_N$ be a basis for $V$ and
$\| \cdot \|'$ be the $l^1$-norm with respect to this basis.
Then, for $a \in V$,
\[
\big|\left\la \ve (f), k_t (a) \ve (g)\right\ra \big| \leq \| a\|' \max_i
\big|\left\la \ve (f), k_t (e_i) \ve (g)\right\ra \big|.
\]
Since all
norms on $V$ are equivalent, the result follows.
\end{proof}

Define $\phi:\Coalg \to \Coalg \odot \Op(\Dhat)$ by
\begin{equation} \phi = (\id_{\Coalg} \ot \varphi) \circ \Com,
\label{bigphi} \end{equation}

\begin{propn} \label{E}
The coalgebraic QS differential equation \eqref{qgQSDE} has at most one  weak solution.
\end{propn}

\begin{proof}
Let $k \in \Proc (\Coalg ; \Exps_D )$ be the difference of two
weak solutions of \eqref{qgQSDE}, and let $a \in \Coalg$,
$f,g \in \Step_D$ and $T \geq 0$. Then
$k\in\Procwc (\Coalg;\Exps_D )\subset\Procwr (\Coalg ; \Exps_D )$,
by \eqref{wc=wr}.
By iteration
\[
\lla \ve (f), k_t (a) \ve (g)\rra = \int_{\Delta_n[0,t]}
\lla \ve (f), k_{s_1} \Big( \phi^{\fhat (s_1)}_{\ghat (s_1)} \circ
\cdots
\circ\phi^{\fhat (s_n)}_{\ghat (s_n)}(a)\Big) \ve (g)\rra\, d\bm{s}
\]
for each $n\in\bn$ and $t\in [0,T]$, where for $\xi, \eta \in \Dhat$
\[
\phi^{\xi}_{\eta}: =
\big(\id_{\Coalg} \ot (\omega_{\xi ,\eta} \circ \varphi)  \big) \circ \Com \]
(recall the notation \eqref{omegas}).
The coalgebra $\Coalg_a$, defined as the  subcoalgebra of $\Coalg$
generated by $a$ is finite-dimensional by Theorem \ref{FTOC}.
As each $\phi^{\xi}_{\eta}$ leaves $\Coalg_a$
invariant, fixing a norm for $\Coalg_a$
and appealing to Lemma \ref{D} allows us to claim that the integrand is bounded by
\[
\Big( \max \big\{ \| \phi^{\fhat (s)}_{\ghat(s)} \| : 0 \leq s \leq T  \big\}
\Big)^n C_{f,g,T,\Coalg_a }\| a\|.
\]
The result follows.
\end{proof}

\subsection*{Existence of the solution}

Let $\up^{\varphi}$ be the linear
map $\Coalg \rightarrow \Seq_D$ defined by
$\up^\varphi (a)_n = \up^\varphi_n (a)$ ($n \in \bn_0$) where
$\upsilon_0^{\varphi}:=\Cou$ and,
in the notation \eqref{Com n},
\begin{equation} \label{u phi n}
\up^\varphi_n = \varphi^{\ot n} \circ \Com_{n-1}
: \Coalg \rightarrow
\Op (\Dhat)^{\odot n} \subset \Op(\Dhat^{\odot n})\
\text{ for } n\geq 1.
\end{equation}
Note the recursive identity
\begin{equation} \label{u n+1}
\up^\varphi_{n+1} = (\up^\varphi_n \ot \varphi) \circ \Com.
\end{equation}
If $\Coalg$ is involutive and $\varphi\in L(\Coalg; \Opdag(\Dhat))$ then
\begin{equation} \label{u phi dagger}
(\up^\varphi )^\dagger =
\up^{\varphi^{\dagger}}.
\end{equation}

In terms of the associated map
$\phi$ introduced in \eqref{bigphi}
\[
\upsilon^{\varphi}_n = \epsilon_n \circ \phi_n,
\]
where the maps
$\phi_n : \Coalg \rightarrow \Coalg \odot \Op (\Dhat)^{\odot n}$
and
$\Cou_n :
\Coalg \odot \Op (\Dhat)^{\odot n} \rightarrow \Op (\Dhat )^{\odot n}$
are defined (recursively) by
$\phi_0 = \id_{\Coalg}$, $\Cou_0 = \Cou$ and, for $n\geq 1$,
$\phi_n = (\phi \ot \id_{\Op (\Dhat )^{\odot n}} ) \circ
\phi_{n-1}$, $\Cou_n = \Cou \ot \id_{\Op (\Dhat )^{\odot n}}.$

Define additionally $\wt{\up}^{\varphi}:\Coalg \to \Seq_D$ by
\begin{equation} \label{upflipped}
\wt{\up}^{\varphi}_n = \tau_n \circ \up^{\varphi}_n, \;\;\; n \in \bn_0,
\end{equation}
where $\tau_n: \Op (\Dhat^{\odot n}) \to \Op (\Dhat^{\odot n})$
denotes the flip reversing the order of the $n$ copies of $\Dhat$.

\begin{lemma} \label{F}
For any
$\varphi \in L(\Coalg ; \Op (\Dhat))$, $\up^\varphi\in L(\Coalg ; \Geom_D)$ and
$\wt{\up}^\varphi\in L(\Coalg ; \Geom_D)$.
Moreover, if
$\varphi$ is $\Opdag (\Dhat)$-valued then both $\up^\varphi$ and $\wt{\up}^\varphi$
take values in $\Geomdagger_D$.
\end{lemma}

\begin{proof}
Fix an element $a \in \Coalg \setminus \{ 0\}$
and let again $\Coalg_a$ denote the coalgebra generated by $a$.
By Theorem \ref{FTOC} $\Coalg_a$ is finite dimensional;
let $a^1 ,\ldots ,a^N$ be its linear basis in which $a^1 =a$. Let
$(\nu^i_{jk})$ be the coefficients of $\Com$ (viewed as a map
$\Coalg_a \rightarrow \Coalg_a \ot \Coalg_a$)
 with respect to this basis, i.e.\ for all $i=1,\ldots,N$
 \[ \Com(a^i) =\sum_{j,k=1}^N \nu_{j,k}^i a^j \ot a^k.\]
Set ($i,j=1,\ldots,N$)
\[
T^i_j = \sum_k \nu^i_{jk} \varphi (a^k) \in \Op (\Dhat).
\]
Then $\phi (a^i) = \sum_{j=1}^N a^j \ot  T^i_j  $ and
\begin{equation} \label{upinHD}
\up^\varphi_n (a) = \sum_{\mathbf{k}} \Cou (a^{k_n}) T_{k_n}^{k_{n-1}} \ot
T_{k_{n-1}}^{k_{n-2}} \ot \cdots \ot T_{k_1}^1 ,
\end{equation}
a sum of $N^n$ terms of the form $X_1 \ot \cdots \ot X_n$ in which
\[
X_i \in \{ T^j_k : 1 \leq j,k \leq N\} \cup \{ \Cou (a^j) T^l_j
: 1 \leq j,l \leq N \},
\]
so $\up^\varphi (a) \in \Geom_D$.
The rest is easily verified.
\end{proof}

\begin{tw} \label{G}
Let $\varphi \in L\left(\Coalg ; \Op (\Dhat )\right)$ and set
$\wt{\up} = \wt{\up}^\varphi$. Then the process $k := \Lambda \circ \wt{\up}$
strongly satisfies the coalgebraic quantum stochastic differential equation
\eqref{qgQSDE}.
\end{tw}

\begin{proof}
By \eqref{Lambda map} and Lemma \ref{F} the process $k$ is
continuous. It therefore suffices to
show that it satisfies the equation weakly. In Guichardet notation,
\eqref{Guichardet FF1}, \eqref{u n+1} and \eqref{upflipped} imply
\begin{align*}
&\int^t_0 ds \lla \fhat (s) , \varphi (a_{(2)}) \ghat (s) \rra
\lla \ve (f), k_s (a_{(1)}) \ve (g) \rra \\
= &\int^t_0 \, ds \int_{\Gamma_{[0,s]}} d \tau
\lla \pi_{\fhat} (\tau \cup s), \big(  \varphi (a_{(2)})  \ot
\wt{\up}_{\# \tau} (a_{(1)}) \big) \pi_{\ghat} (\tau \cup s) \rra
\lla \ve (f), \ve (g)\rra \\
=& \int_{\Gamma_{[0,t]}} d\sigma \big(1-\delta_\emptyset (\sigma) \big)
\lla \pi_{\fhat} (\sigma) ,\wt{\up}_{\# \sigma} (a) \pi_{\ghat} (\sigma) \rra
\lla \ve (f) , \ve (g)\rra  \\
= & \left\la \ve (f), k_t (a)\ve (g) \rra - \epsilon (a)
\lla \ve (f) , \ve(g) \right\ra,
\end{align*}
and so, by \eqref{wQSDEsweedler}, $k$ is the weak solution of \eqref{qgQSDE}.
\end{proof}

Thus the coalgebra QS differential equation \eqref{qgQSDE} has a unique  weak
solution; it
is a strong solution and is given by  $\Lambda \circ \wt{\upsilon}^\varphi$ ---
we denote it $l^\varphi$.

\subsection*{Properties of the solution}

\begin{lemma} \label{traffic}
Let $\varphi \in L\left(\Coalg ; \Op (\Dhat) \right)$. Then
$l^\varphi \in \ProcHc (\Coalg ; \Exps_D)$ and the following hold.
\begin{alist}
\item
The map $\varphi \mapsto l^\varphi$ is injective.
\item
If $\varphi \in L\left( \Coalg ; \Opdag (D) \right)$ then
$l^\varphi \in \ProcHc^{\dagger} (\Coalg ; \Exps_D)$.
\item
If $\Coalg$ is unital then $l^\varphi$ is unital if and only if
$\varphi (1) =0$.
\item
If $\Coalg$ is involutive and $\varphi \in L \left( \Coalg ; \Opdag
(\Dhat)\right)$ then
$(l^\varphi )^\dagger = l^{\varphi^{\dagger}}$.
In particular,
$l^\varphi$ is real if and only if $\varphi $ is real.
\end{alist}

\end{lemma}

\begin{proof}
In view of the inclusion $\Lambda(\Geom_D) \subset\ProcHc (\Exps_D)$ and Lemma
\ref{F},
$l^\varphi$ belongs to $\ProcHc (\Coalg ; \Exps_D)$.
\newline
(a) follows from the identity
\begin{equation} \label{generator}
\left\la \chat ,\varphi (a) \dhat \, \right\ra =
\lim_{t \rightarrow 0^+} t^{-1}
\Big(\left\la \ve (c_{[0,t)}), l^\varphi_t (a) \ve (d_{[0,t)}) \right\ra
- \epsilon (a) \mathrm{e}^{t\la c,d \ra} \Big) ,
\end{equation}
and the totality of $\{ \chat : c \in D\}$ in $\kilhat$.
\newline
(b) follows from Lemma \ref{F} and \eqref{Lambdadag}.
\newline
(c) follows from \eqref{wQSDEsweedler}, \eqref{generator} and
the unitality of $\Cou$ and $\Com$.
\newline
(d) By part (b) $l^\varphi\in\ProcHc^{\dagger}(\Coalg,\Exps_D)$ and
by \eqref{u phi dagger} and \eqref{Lambdadag}
\[
(l^\varphi)^\dagger = \big(\Lambda \circ (\wt{\up}^\varphi)\big)^\dagger
= \Lambda \circ \wt{\up}^{\varphi^{\dagger}} = l^{\varphi^{\dagger}} .
\]
The last part follows by injectivity.
\end{proof}

\begin{rem}
If we cast $\varphi \in L\big( \Coalg ; \Op (\Dhat ) \big)$ in block
matrix form:
\[
\varphi =
\begin{bmatrix}
\gamma & \alpha \\ \chi & \nu - \iota\circ\Cou
\end{bmatrix},
\]
where $\iota (z)=zI_{D}$,
then $\gamma \in L(\Coalg ; \bc)$, $ \chi \in L\big(\Coalg ; L (\bc;
D)\big)$,
$\alpha \in L\big(\Coalg ; L(D;\bc)\big)$ and $\nu \in L\big(\Coalg ;
\Op(D)\big)$,
so that
\begin{equation}\label{block}
\varphi(a) =
\begin{bmatrix}
\gamma (a) &\alpha (a) \\
\chi (a) & \nu (a)-\Cou (a) I \end{bmatrix}.
\end{equation}
Moreover
$\varphi \in L\big(\Coalg ; \Opdag (\Dhat )\big)$ if and only if
$\alpha (\Coalg )\subset \big\la D \big| :=
\big\{ \la d| \,\big|\; d \in D\big\}$ and
$\nu (\Coalg ) \subset \Opdag (D)$.
Thus if $\Coalg$ is involutive then
\[
\varphi^\dagger =
\begin{bmatrix}
\gamma^\dagger & \chi^{\dagger} \\
\alpha^{\dagger} & \nu^\dagger - \iota\circ\Cou
\end{bmatrix}.
\]
In particular,
\[
\varphi =
\varphi^\dagger \text{ if and only if }
\gamma = \gamma^\dagger, \,
\nu = \nu^\dagger \text{ and }
\alpha = \chi^{\dagger}.
\]
\end{rem}

\section{Stochastically generated convolution cocycles} \label{stochsection}
\setcounter{equation}{0}

In this section we first note that solutions of coalgebraic QS
differential equations are QS convolution cocycles.
Heeding the fact that solutions are
$\frac{1}{2}$-H\"{o}lder continuous
we then establish the converse: every QS convolution cocycle in
$\ProcHc^{\dagger} (\Coalg ; \Exps_D)$ necessarily satisfies a coalgebraic QS
differential equation.

\begin{lemma} \label{J}
Let $\varphi \in L(\Coalg ; \Op (\Dhat ))$. Then
$\up := \up^\varphi$  defined in \eqref{u phi n}  satisfies
\[
\up_{n+m} = (\up_n \ot \up_m) \circ \Com
\]
for all $n,m \in \bn_0$.
\end{lemma}

\begin{proof}
Since $\up_k = \varphi^{\ot k} \circ \Com_{k-1}$, this reduces to the identity
\eqref{Com ij}.
\end{proof}

\begin{propn} \label{K}
Let $l=l^{\varphi}$ where $\varphi \in L(\Coalg ; \Op (\Dhat))$.
Then $l\in \QSCHc(\Coalg; \Exps_D)$. If $\Coalg$ is involutive and $\varphi$
takes values in $\Opdag(\Dhat)$ then $l\in \QSCdagHc(\Coalg; \Exps_D)$

\end{propn}

\begin{proof}
Set $\wt{\up} = \wt{\up}^\varphi$ as in Theorem \ref{G} and let $a \in \Coalg$.
Using Guichardet notation again,  we obtain for all $f,g \in \Step_D$, $s,t \geq 0$
\begin{align*}
\lefteqn{\lla \ve(f), (l_s \star ( \sigma_s \circ l_t) )(a) \ve(g) \rra =
\lla \ve(f), \left(l_s (a_{(1)}) \ot (\sigma_s \circ l_t) (a_{(2)}) \right) \ve(g) \rra} \\
&=\left\la \ve (f_{[0,s)} ) ,l_s (a_{(1)}) \ve (g_{[0,s)} ) \right\ra
\left\la \ve(f_{[s,\infty)}), \sigma_s (l_t (a_{(2)})) \ve (g_{[s,\infty)})\right\ra  \\
& =\int_{\Gamma_{[0,s[}} \, d\sigma
\left\la \pi_{\fhat} (\sigma), \wt{\up}_{\# \sigma} (a_{(1)})
\pi_{\ghat} (\sigma)\right\ra \int_{\Gamma_{[s,s+t[}} \, d \tau
\left\la \pi_{\fhat} (\tau ), \wt{\up}_{\# \tau} (a_{(2)}) \pi_{\ghat} (\tau )\right\ra
\mathrm{e}^{\la f,g \ra} \\
& =\int_{\Gamma_{[0,s+t[}} \, d \omega\left\la \pi_{\fhat}
(\omega),   \wt{\up}_{\# \omega \cap [s,s +t)} (a_{(2)}) \ot \wt{\up}_{\# \omega \cap [0,s)}
(a_{(1)}) \pi_{\ghat}
(\omega)\right\ra \mathrm{e}^{\la f,g\ra} \\
& =\int_{\Gamma_{[0,s+t[}} \, d \omega\left\la \pi_{\fhat}
(\omega), \tau_{\# \omega} \left( \up_{\# \omega \cap [0,s)}
(a_{(1)}) \ot \up_{\# \omega \cap [s,s +t)} (a_{(2)}) \right) \pi_{\ghat}
(\omega)\right\ra \mathrm{e}^{\la f,g\ra} .
\end{align*}
Applying Lemma \ref{J}, and identity \eqref{Guichardet FF1} once more, we see that this is
equal to $\la \ve (f), l_{s+t} (a) \ve (g) \ra$, which proves that $l$ is a convolution
increment process. As its initial condition is given by the counit, Lemma \ref{traffic}
ends the proof.
\end{proof}

Thus
$\varphi \mapsto l^\varphi$ gives  maps
\[
L\big(\Coalg ; \Op (\Dhat) \big) \rightarrow \QSCHc (\Coalg;\Exps_D) \text{
and
}
L\big(\Coalg ; \Opdag (\Dhat) \big) \rightarrow \QSCdagHc (\Coalg;
\Exps_D) .
\]
These maps are injective, by Lemma \ref{traffic};
our aim now is to establish bijectivity of the second map.
When $l = l^\varphi$ we refer to $\varphi$ as the \emph{stochastic generator} of
the QS convolution cocycle $l$.
Note that if $l$ is a cocycle in $\Proc^{\dagger} (\Coalg ; \Exps_D)$, and
$\Coalg$
is involutive, then $l^\dagger$ is a cocycle too.

Recall Definition \ref{defsemig} and Corollary \ref{N}.

\begin{lemma}  \label{O}
Let $l = l^\varphi $ for $\varphi \in L\big(\Coalg ; \Op (\Dhat)\big)$.
Then its
associated CCSFs $\lambda^{c,d}$ have generators
$ \la \chat, \varphi (\,\cdot\,) \dhat \ra$   ($ c,d\in D$).
\end{lemma}

\begin{proof}
Since $l_t (a) = \Lambda_t \big( \wt{\up}^\varphi (a) \big)$ this is an
immediate consequence of
\eqref{Guichardet FF1}
and the definition of $\upsilon^{\varphi}$ by formula \eqref{u phi n} - see also
the proof of Lemma \ref{traffic}.
\end{proof}

\begin{propn}  \label{sesquilinear}
Let $l \in \QSCdagc (\Coalg ; \Exps_D)$. Define a map
\[
q :\Dhat \times\Dhat \rightarrow L (\Coalg ; \bc), \quad
\left( \binom{z}{c},\binom{w}{d} \right) \mapsto
\begin{bmatrix}\overline{z-1}&1\end{bmatrix}\,
\begin{bmatrix}
\gamma_{0,0}&\gamma_{0,d} \\
\gamma_{c,0}&\gamma_{c,d}
\end{bmatrix} \, \begin{bmatrix} w-1\\1 \end{bmatrix} ,
\]
where $\{ \gamma_{c,d} :c, d\in D\}$ are the generators of the CCSFs
associated with $l$. Then $q$ is sesquilinear.
\end{propn}

\begin{proof}
The proposition amounts to the sesquilinearity of each form
$q_a := q(\cdot , \cdot)(a)$. Thus let $a \in \Coalg$.
First note the identity
\[
q_a \left( \chi , \eta \right)  = \lim_{t \rightarrow 0^+}
t^{-1} \lla \xi (t), \big( w,d_{[0,t)} , (2!)^{-1/2} (d_{[0,t)})^{\ot 2},
\cdots \big) \rra ,
\]
for $\chi = \binom{z}{c}$ and $\eta = \binom{w}{d}  $ in $\Dhat$,
where
\[
\xi (t) = \big[  l_t (a)^\dagger - \overline{\epsilon (a)} I \big]
\big( (z-1) \ve (0) + \ve (c_{[0,t[}) \big) .
\]
Thus if
$\eta = \eta_1 + \alpha \eta_2$, $\eta_1 = \binom{w_1}{d_1}$, $\eta_2 = \binom{w_2}{d_2}  $ then
\[
q_a \left( \chi ,\eta \right)- q_a \left( \chi,
\eta_1\right) - \alpha q_a \left(\chi, \eta_2 \right)
= \lim_{t \rightarrow 0^+} \la \xi (t) , \zeta (t) \ra
\]
where
\[
\zeta (t) = t^{-1} \Big( (n!)^{-1/2}
\big[
d^{\ot n} - (d_1)^{\ot n} - \alpha (d_2)^{\ot n}
\big]
\ot 1_{[0,t[^n} \Big)_{n \geq 2} .
\]
Since $\zeta$ is locally bounded and $\xi (t) \rightarrow 0$ as $t \rightarrow 0^+$,
by the  continuity of the process $\big(l_t(a)^\dagger\big)_{t \geq 0}$, this shows that
$q_a$ is linear in its second argument.
A very similar argument, this time using the continuity of the process
$\big(l_t (a)\big)_{t \geq 0}$, shows that $q_a$ is conjugate linear in
its first argument. The result follows.
\end{proof}

\begin{propn} \label{P}
Let $l \in\QSCdagHc (\Coalg ; \Exps_D)$ and let  $q$ be defined as in
\tu{Proposition \ref{sesquilinear}}.
Then, for each $a \in \Coalg$, the sesquilinear form $q( \cdot , \cdot )
(a)$ is
separately continuous in each argument.
\end{propn}
\begin{proof}

Let $\{ \gamma_{c,d} : c, d \in D \}$ be the generators
of CCSFs associated with $l$, let $a \in \Coalg$, and $\chi = \binom{z}{c}, \, \eta =
\binom{w}{d} \in \Dhat$.
Then
\begin{align*}
q(\chi , \eta)(a) =
& \overline{z} \big( (w-1) \gamma_{0,0} (a) +
\gamma_{0,d} (a) \big) + \\
& (w-1) \big( \gamma_{c,0} (a) - \gamma_{0,0} (a) \big)
  + \big( \gamma_{c,d} (a) - \gamma_{0,d} (a) \big)
\end{align*}
and, for $e \in D$ and $T>0$,
\begin{align*}
\lefteqn{\big| \gamma_{c,e} (a) - \gamma_{0,e} (a) \big|} \\
&=
\lim_{t\to 0^+} t^{-1}
\bigg|
\Big\la
e^{-t\la e,c \ra}
\ve (c_{[0,t)})
- \ve (0), \big( l_t (a) - \epsilon (a) \big) \ve (e_{[0,t)}) \Big\ra \bigg| \\
&=
\lim_{t\to 0^+} t^{-1} \bigg| e^{-t\la c,e \ra} \Big\la \ve (c_{[0,t)})
- \ve (0), \big( l_t (a) - \epsilon (a) \big) \ve (e_{[0,t)}) \Big\ra \bigg| \\
&\leq \limsup_{t\to 0^+}  t^{-1/2} \big\| \ve (c_{[0,t)}) - \ve (0) \big\|
t^{-1/2} \Big\| \big[ l_t (a) - \epsilon (a) \big] \ve (e_{[0,T)}) \Big\|
\mathrm{e}^{ -\| e\|^2 (T-t)/2} \\
& \leq \| c \| C(a,e,T)
\end{align*}
for some constant $C$ depending only on $a$, $e$ and $T$.
Thus, setting $T=1$,
\begin{align*}
\big| q(\chi ,\eta )(a) \big|
& \leq |z| \big| (w-1) \gamma_{0,0} (a) + \gamma_{0,d} (a) \big|
 + \| c \| \big( |w-1| C(a,0,1) + C(a,d,1) \big) \\
& \leq M\| \chi \|
\end{align*}
for a constant $M$ depending only on $a$ and $\eta$. This establishes continuity
in the first argument.
Continuity in the second argument is proved by a very similar argument,
this time using
the H\"{o}lder-continuity of the process $\big( l_t (a)^\dagger \big)_{t
\geq 0}$.
\end{proof}

\begin{rem}
Proposition \ref{P} may be proved in the same way under the following weaker hypothesis:
estimates of the type
\[
\limsup_{t \to 0^+} t^{-1} \Big| \big\la \ve (c_{[0,t)}) - \ve (0) ,
l_t (a) \ve (d_{[0,t)}) \big\ra \Big| \leq M_{a,d}\| c\|,
\]
for constants $M_{a,d}$ depending only on $a$ and $d$ should be satisfied
by both $l$ and $l^{\dagger}$.
\end{rem}

We are ready for the main result of this section.

\begin{tw} \label{Q}
Let $k \in \Proc (\Coalg ;\Exps_D )$. Then the following are
equivalent \tu{:}
\begin{rlist}
\item $k \in \QSCdagHc (\Coalg ; \Exps_D )$;   \label{Qb}
\item
$k = l^\varphi$ for some $\varphi \in L\big(\Coalg ; \Opdag
(\Dhat)\big)$.
\label{Qa}
\end{rlist}
\end{tw}

\begin{proof}
Let $k \in \QSCdagHc (\Coalg ; \Exps_D)$ and let $\{ \gamma_{c,d} :c,d
\in D \}$ denote
the generators of CCSFs associated with $k$. By Propositions \ref{sesquilinear}
and \ref{P}, there is a
map $\varphi \in L\big(\Coalg ;\Opdag (\Dhat) \big)$ such that
\[
\lla \binom{z}{c} , \varphi (a) \binom{w}{d} \rra =
\begin{bmatrix} \overline{z-1} & 1 \end{bmatrix} \,
\begin{bmatrix} \gamma_{0,0} (a) & \gamma_{0,d} (a) \\
\gamma_{c,0} (a) & \gamma_{c,d} (a) \end{bmatrix} \,
\begin{bmatrix} w-1 \\ 1 \end{bmatrix} ,
\]
in particular,
\[
\big\la \chat, \varphi (\,\cdot\,) \dhat \, \big\ra = \gamma_{c,d} .
\]
Thus, by Lemma \ref{O}, the QS convolution cocycles $l^\varphi$
and $k$
have the same CCSFs and so coincide. Thus (\ref{Qb}) implies
(\ref{Qa}). The converse has already been established in Proposition \ref{K} and
Lemma \ref{traffic}(c).
\end{proof}

\begin{rem}
The transformation between the family $\{ \gamma_{c,d} :c,d \in D\}$
and $\varphi$ is a familiar one in the
analysis of stochastic cocycle generators (cf.\ \cite{lwjfa}).
\end{rem}

As a special case of theorem \ref{Q} for $\Coalg = \bc$, we obtain
the following result.

\begin{cor} \label{L}
Let $X \in \ProcHc^{\dagger} (\Exps_D)$. Then the following are
equivalent\tu{:}
\begin{rlist}
\item $X$ is an operator Markovian cocycle;
\item $X$ satisfies a QS differential equation of the form $dX_t = (L \ot X_t )d\Lambda_t$,
$X_0 =I$, for some $L \in \Opdag (\Dhat)$.
\end{rlist}
\end{cor}
\noindent
This type of cocycle is used in Section \ref{perturb} for perturbing
QS convolution cocycles on a general coalgebra.

\section{Multiplicativity}   \label{multalg}
\setcounter{equation}{0}

The section is devoted to characterising the stochastic generators
of weakly multiplicative (and $^*$-homomorphic) QS cocycles on a
$^*$-bialgebra $\Alg$.

\subsection*{Matrix-sum kernels revisited}

First we introduce certain further spaces of unbounded operators
useful for formulating the It\^o Formula for iterated QS integrals in the
purely algebraic context.
For a dense subspace $E$ of a Hilbert space $\hil$ let
\begin{align*}
&\Opinv (E) :
= \big\{ T \in \Op (E) \, : \; \Ran T \subset E \big\},\\
&\Opstar (E) :
= \big\{ T \in \Opdag (E) \, : \; T,T^\dagger \in \Opinv (E) \big\}
\end{align*}
(``inv'' for invariant).

Operator composition
$\Opdag (E) \times \Opinv (E) \rightarrow \Op (E)$ extends to pairs
$(S,T)$ in $\Opdag (E) \times \Op (E)$ for which
$\Dom (S^\dagger)^* \supset \Ran T$, as follows:
\begin{equation} \label{bullet}
S \schur T : = (S^\dagger)^* T.
\end{equation}
This partially defined product is bilinear in an obvious sense.
Associativity relations however have to be justified (if needed) in every considered case separately.

The definitions above enable us to introduce the following subspaces of $\Seq_D$:
\begin{align*}
&\Seqinv_D :
=
\big\{ F \in \Seq_D \, : \;
\forall_{n \in \bn}\, F_n \in \Opinv (\Dhat^{\odot n}) \big\}, \\
&\Seqstar_D :
= \big\{ F \in \Seqdagger_D \, : \; F,F^\dagger \in \Seqinv_D \big\}.
\end{align*}

To introduce the \emph{matrix-sum convolution product}
$*: \Seq_D \times \Seqinv_D \rightarrow \Seq_D$, reflecting the multiplicative nature of iterated QS
integrals, we need some more notations.
For  $n \in \bn$ and $\alpha \subset \{ 1,\ldots ,n \}$ write $\alpha = \{ \alpha_1 < \cdots < \alpha_k \}$ and
$\{ 1, \ldots ,n\} \setminus \alpha = \{ \overline{\alpha}_1 < \cdots < \overline{\alpha}_{n-k} \}$ and define
$\Pi_{\alpha ;n} \in \Opstar(\Dhat^{\odot n} )$ as a certain tensor flip - the linear extension of the map
\[
\chi_1 \ot \cdots \ot \chi_n \mapsto \chi_{\alpha_1} \ot \cdots \ot \chi_{\alpha_k}
\ot \chi_{\overline{\alpha}_1} \ot \cdots \ot \chi_{\overline{\alpha}_{n-k}}.\]
If now $F \in \Seq_D$, define
\begin{equation} \label{Pi}
F(\alpha ;n) := \Pi^*_{\alpha ;n}
\big( F_k \ot I_{n-k} \big)  \Pi_{\alpha ;n}
 \in \Op (\Dhat^{\odot n}),
\end{equation}
where to lighten the notation we write $I_l$ for $I_{\Dhat^{\odot l}} $ ($l\in \bn_0$). Observe that in particular
\[ F(\emptyset;n) = F_0 I_n.\]
Additionally let
\[\QSproj  [\alpha ;n ]:= \Pi^*_{\alpha ;n} \left( (\QSproj)^{\ot k} \ot I_{n-k} \right)\Pi_{\alpha ;n}.
\]
Thus if $F_k$ is a simple tensor $ T_1 \ot \cdots \ot T_k$,
$\alpha = \{\alpha_1<\ldots < \alpha_k\}$, then
$F(\alpha ;n)=S_1 \ot \cdots \ot S_n$, where
\[  S_i = \begin{cases}  T_j & \text{if } i \in \alpha,
 i = \alpha_j \\ I &\text{if } i \notin \alpha \end{cases}.\]
The announced above product $*$ is defined for all pairs $(F,G) \in \Seq_D \times \Seqinv_D$ by
\begin{equation} \label{convolution product}
(F \ast G)_n = \sum_{| \boldsymbol{\alpha} | = \{ 1,\ldots,n\} }
F(\alpha_1 \cup \alpha_2 ;n) \QSproj [\alpha_2 ; n] G(\alpha_2 \cup \alpha_3 ; n)
\end{equation}
where $n \in \bn_0$ and the sum is over all $3^n$ disjoint partitions
 $\alpha_1 \cup \alpha_2 \cup \alpha_3$ of
$\{ 1, \ldots ,n\}$.
Heuristically, the product $*$ may be equivalently thought of as a sum over all possible
fillings of $n$ places with operators coming from $F$ and $G$ with $\QSproj$ intervening whenever
chosen operators `overlap'. To clarify this statement note that for $n=1$ \eqref{convolution product} gives
\[ (F \ast G)_1 = F_0 G_1 + G_0 F_1 + F_1 \QSproj G_1.\]

\begin{lemma} \label{B}
The product enjoys the following properties\tu{:}
\begin{alist}
\item if $F,G \in \Seqinv_D$ then $F \ast G \in \Seqinv_D$\tu{;}
\item if $F \in \Seq_D$ and $G, H \in \Seqinv_D$ then
$F \ast (G \ast H) = (F * G) *H$\tu{;}  \label{Bb}
\item if $E = 1 \delta_0$ then $E \in \Seqstar_D$ and
$E *F=F*E=F$ for all $F \in \Seq_D$\tu{;}
\item if $F,G \in \Seqdagger_D$ with $ G,F^\dagger \in \Seqinv_D$ then
$F*G \in \Seqdagger_D$ and
$(F *G)^\dagger = G^\dagger * F^\dagger$.
\end{alist}
In particular, $(\Seqstar_D , *)$ is a unital $^*$-algebra.

\end{lemma}
\begin{proof}
To see (\ref{Bb}) note that
\[
\sum_{| \bm{\alpha} |=\{ 1,\ldots ,n\} } F(\alpha_1 \cup \alpha_2 ;n)
\QSproj [\alpha_2 ;n] G(\alpha_2 \cup \alpha_3 ;n) \QSproj [\alpha_3 ;n] H(\alpha_3 \cup \alpha_4 ;n)
\]
is a common expression for $\big( F*(G*H) \big)_n$ and
$\big( (F*G)*H\big)_n$. The rest is easily verified.
\end{proof}

We have already seen in Section \ref{QSIntegrals} that for QS purposes a growth condition needs to be imposed
on elements of $\Seq_D$.
To obtain algebras of processes we need to restrict further.
The choice of restriction here is directly motivated by the Fundamental
Theorem on Coalgebras and its consequences expressed in the proof of Lemma \ref{F}.
Thus let $\cH_D$ denote the set of $F \in \Seq_D$ satisfying
\begin{multline}
\exists_{p,q \in \bn ,R \subset\subset \Op(\Dhat)}
\, \forall_{n \in \bn_0}\;
F_n \text{ may be expressed as a sum of }
\\
pq^n \text{ terms of the form }
X_1 \ot \cdots \ot X_n \text{ with } X_1 , \ldots , X_n \in R,
\label{HsubD}
\end{multline}
with $\cH^\dagger_D ,\cH^{\mathrm{inv}}_D$ and $\cH^*_D$ defined as for $\Seq$
(the same remark applies to $\Geominv_D$, $\Geom^*_D$, see \eqref{GsubD}, \eqref{Gsuper}).
It is elementary to check that all of these are subspaces of $\Geom_D$.

\begin{propn}\label{2.2*}
Let $F \in \Geom_D$ and $G \in \cH^{\mathrm{inv}}_D$. Then $F * G \in \Geom_D$,
moreover if $F \in \cH_D$ then $F*G \in \cH_D$ too.
\end{propn}

\begin{proof}
Let $H = F *G$ and
choose $p,q$ and $R$ for $G$
according to \eqref{HsubD}. Let $S \subset\subset \Dhat$ and let
$n \in \bn$ and $\chi_1, \ldots, \chi_n \in S$. Then, for any partition
$\alpha \cup \beta \cup \gamma$ of $\{ 1, \ldots , n \}$,
\[
F(\alpha \cup \beta ; n) \QSproj [\beta ;n] G(\beta \cup \gamma ; n)
(\chi_1 \ot \cdots \ot \chi_n )
\]
is a sum of $pq^{\# (\beta \cup \gamma)}$ terms of the form
\[
F(\alpha \cup \beta ; n) (\eta_1 \ot \cdots \ot \eta_n)
\]
where each $\eta_i$ belongs to the finite set
$S' :=RS \cup \QSproj RS \cup S$. Thus, choosing $C_1$ and $C_2 \geq 1$ for the pair
$(F, S')$ according to \eqref{GsubD},
and setting $M= \max \{ \| \eta \| : \eta \in S' \}$,
\[
\| H_n (\chi_1 \ot \cdots \ot \chi_n) \| \leq
\sum pq^{\# (\beta \cup \gamma)} C_1 C_2^{\# (\alpha \cup \beta )}
M^{\# \gamma} = C'_1 (C'_2)^n ,
\]
where $C'_1 = p C_1$ and $C'_2 = 3(C_2 + qC_2 +qM)$. Thus
$H \in \Geom_D$.

If $F \in \cH_D$ then, choosing $p',q'$ and $R'$ for
$F$ (and assuming without loss that $I \in R \cap R'$),
$F(\alpha \cup \beta ;n) \QSproj [\beta ;n] G(\beta \cup \gamma ;n)$
is a sum of $p' (q')^{\# (\alpha \cup \beta )} pq^{\# (\beta \cup \gamma)}$
terms of the form $Z_1 \ot \cdots \ot Z_n$ where
$Z_i \in R' R \cup R' \QSproj R$. Thus $H_n$ is a sum of
$pp'3^n(q+qq'+q')^n$ terms of this form and $H \in \cH_D$.
\end{proof}

As an immediate consequence we have the following result.

\begin{tw} \label{cH as algebra}
$(\cH^*_D , *)$ is a unital $*$-subalgebra of $(\Seqstar_D , *)$.
\end{tw}

\subsection*{Multiplicativity for iterated QS integrals}

In this subsection we show how the product $\ast$ introduced above
corresponds to the multiplicative structure of iterated QS integrals.

Elementary inductive argument using \eqref{SFF} yields the following fact
(for a more sophisticated version,
relevant for $C^*$-algebraic processes see Theorem 2.2 of \cite{lwhom}; it is
an essential tool leading to Theorems \ref{0Corolla} and \ref{Corolla}).

\begin{fact}   \label{C}
Let $F \in \Geomdagger_D$ and $G \in \Geominv_D$. Then, for each
$n \in \bn$, $f,g \in \Step_D$,
\[
\sum^n_{i,j=0} \left\la \Lambda^i_t (F_i) \ve (f), \Lambda^j_t (G_j) \ve
(g) \right\ra
= \sum^{2n}_{k=0} \lla \ve (f), \Lambda^k_t (H_k) \ve (g) \rra
\]
where $H = F^{\dagger}*G$.
\end{fact}

Here is the relevant consequence of the above fact for the
integral $\Lambda$ encompassing all elements of the sequence at once. Recall the
partially defined product \eqref{bullet}.

\begin{propn} \label{multiplicativity}
Let $F \in \Geomdagger_D$ and $G \in \Geominv_D$ be such that $F*G \in \Geom_D$.
Then, for each $t \geq 0$,
$\Dom \Lambda_t (F)^{\dagger *} \supset \Ran \Lambda_t (G)$
and
\[
\Lambda_t (F *G) = \Lambda_t (F) \schur \Lambda_t (G) .
\]
\end{propn}

\begin{proof}
For $H \in \Seq_D$ and $n \in  \bn_0$ write
$\Lambda^{[n]} (H)$ for $\sum^n_{i=0} \Lambda^i (H_i) \in \Procc (\Exps_D)$,
so that if $H \in \Geom_D$ then for each $f \in \Step_D$
\[\Lambda^{[n]}_t (H) \ve (f)\stackrel{n \to \infty}{\longrightarrow} \Lambda_t (H) \ve (f).\]
By the fact above and \eqref{Lambdadag}, for all $f,g \in \Step_D$
\begin{align*}
\la \Lambda_t (F)^\dagger \ve (f), \Lambda_t (G) \ve (g) \ra
&= \lim_{n \rightarrow \infty}\lla \ve (f), \Lambda^{[n]}_t (F*G) \ve (g)\rra \\
&=\lla \ve (f), \Lambda_t (F*G) \ve (g)\rra ,
\end{align*}
and so the result follows.
\end{proof}

The following theorem summarises the content of
this subsection and the previous one; it is an immediate consequence of Proposition
\ref{multiplicativity} and Theorem \ref{cH as algebra}.

\begin{tw} \label{Lambda isom}
Let $\cP = \Lambda (\cH^*_D) \subset \ProcHc^{\dagger}(\Exps_D)$.
Then, with respect to the product defined in \eqref{bullet} ---
extended pointwise, the map $\Lambda$ restricts to a unital
$^*$-algebra isomorphism of unital $^*$-algebras:
\[
(\cH^*_D ,*) \rightarrow (\cP , \schur ).
\]
\end{tw}

\subsection*{Multiplicative QS convolution cocycles}

As mentioned before, we fix a bialgebra $\Alg$.

\begin{deft}  \label{multipalg}
A process $k\in \Proc^{\dagger}(\Alg;\Exps_D)$ is called weakly multiplicative if
for all $a, b \in \Alg, t \geq 0$
\[ k_t (a b) =k_t(a) \schur k_t(b).\]
\end{deft}

It is easy to see that by the arguments above
 the question of multiplicativity for a stochastically generated QS convolution
cocycle $l^{\varphi}\in \QSC(\Alg;\Exps_D)$
is equivalent to that of the multiplicativity of the map
$\wt{\up}^\varphi : \Alg \rightarrow\cH^*_D$
derived from $\varphi \in L(\Alg;\Op(D))$
(as described in Section \ref{coalgQSDE}).
The following proposition shows that the latter question may be reduced to
a simple statement concerning properties of $\varphi$.

\begin{propn} \label{Adam's Lemma 5.4}
Let $\up = \up^\varphi$, $\wt{\up} = \wt{\up}^\varphi$ for
$\varphi \in L\big(\Alg ; \Opinv (\Dhat) \big)$.
Then the following are equivalent\tu{:}
\begin{rlist}
\item $\varphi (ab) = \varphi (a) \epsilon (b) + \epsilon (a) \varphi (b)
+ \varphi(a) \QSproj \varphi (b)$ for all $a, b \in \Alg$.  \label{Adam
a}
\item
$\up (ab) = \up (a) * \up (b)$ for all
$a,b \in \Alg$.  \label{Adam b}
\item
$\wt{\up} (ab) = \wt{\up} (a) * \wt{\up} (b)$ for all
$a,b \in \Alg$.  \label{Adam c}
\end{rlist}
\end{propn}

\begin{proof}
(\ref{Adam a}) is contained in (\ref{Adam b}) since
$\up_1 (ab) = \varphi (ab)$ and
$\big( \up (a) *\up (b) \big)_1$ is the right hand side of
(\ref{Adam a}). Thus (\ref{Adam b}) implies (\ref{Adam a}). Conversely,
if (\ref{Adam a}) holds then $\cP (1)$ holds where $\cP (n)$ is the
proposition
\[
\forall_{a,b \in \Alg} \quad
\up_n (ab) = \big( \up (a) * \up (b) \big)_n .
\]
By the multiplicativity of $\Cou = \up_0$, $\cP (0)$ holds.
Assume therefore that $\cP(k)$ holds for $k \leq n$, and fix $a,b \in
\Alg$.
Employing Sweedler notation and using~\eqref{u n+1}, $\cP (1)$ and then
$\cP (n)$,
\begin{align*}
\lefteqn{\up_{n+1} (ab) } \\
& = \varphi (a_{(1)} b_{(1)}) \ot \up_n (a_{(2)} b_{(2)} ) \\
& = \big[ \varphi (a_{(1)}) \Cou (b_{(1)}) + \Cou (a_{(1)})
\varphi (b_{(1)}) + \varphi (a_{(1)}) \QSproj \varphi (b_{(1)}) \big] \ot \\
& \hspace*{3cm}
\sum_{| \bm{\alpha} | = \{ 1, \ldots , n\} }
\up (a_{(2)}) (\alpha_1 \cup \alpha_2 ; n) \QSproj [\alpha_2 ; n] \up (b_{(2)})
(\alpha_2 \cup \alpha_3 ; n ) ,
\end{align*}
where the sum is over all partitions of the set $\{1,\ldots,n\}$ into
three disjoint subsets.
The identity
\[
\varphi (c_{(1)}) \ot \up (c_{(2)}) (\lambda ; n) = \up (c)
(\overset{\bulletrta}{\lambda}; n +1 ),
\]
(for $c \in \Alg$, $\lambda\subset \{1,\ldots,n\}$) gives the following equalities
\begin{align*}
&\varphi (c_{(1)}) \ot \up (c_{(2)})(\lambda ;n) \QSproj [ \lambda \cap \mu ;n]
\up (d) (\mu ;n) \\
& \hspace*{2cm} = \up (c) (\overset{\bulletrta}{\lambda}; n+1) \QSproj
[\overset{\circrta}{\nu} ;n+1 ] \up (d)
(\overset{\circrta}{\mu} ;n+1),
\end{align*}
and
\begin{align*}
&\varphi (c_{(1)}) \QSproj \varphi (d_{(1)}) \ot \up (c_{(2)}) (\lambda ; n)
\QSproj [ \lambda \cap \mu ; n] \up (d_{(2)}) (\mu ; n) \\
& \hspace*{2cm} = \up (c) (\overset{\bulletrta}{\lambda}; n+1) \QSproj
[\overset{\bulletrta}{\nu} ; n+1 ] \up (d)
(\overset{\bulletrta}{\mu} ; n+1 ),
\end{align*}
in which $\nu = \lambda \cap \mu$. Thus
\begin{align*}
\lefteqn{\up_{n+1} (ab)} \\
&= \sum_{| \bm{\alpha} |=\{ 1,\ldots ,n\} }
\Big( \up (a) (\overset{\bulletrta}{\alpha_1} \cup
\overset{\circrta}{\alpha_2} ; n+1 ) \QSproj
[ \overset{\circrta}{\alpha_2} ;n+1 ] \up (b)
(\overset{\circrta}{\alpha_2} \cup \overset{\circrta}{\alpha_3} ; n+1) \\
& \qquad \qquad \qquad +
\up (a) (\overset{\circrta}{\alpha_1} \cup \overset{\circrta}{\alpha_2} ;
n+1) \QSproj [ \overset{\circrta}{\alpha_2} ; n+1 ] \up (b)
(\overset{\circrta}{\alpha_2} \cup \overset{\bulletrta}{\alpha_3} ; n+1 ) \\
& \qquad \qquad \qquad +
\up (a) (\overset{\circrta}{\alpha_1} \cup \overset{\bulletrta}{\alpha_2} ;
n+1) \QSproj [ \overset{\bulletrta}{\alpha_2} ; n+1 ] \up (b)
(\overset{\bulletrta}{\alpha_2} \cup \overset{\circrta}{\alpha_3} ; n+1 ) \Big) \\
&= \big( \up (a) * \up (b) \big)_{n+1} .
\end{align*}
The implication \eqref{Adam a}$\Rightarrow$\eqref{Adam b} therefore follows by induction.

For \eqref{Adam b}$\Leftrightarrow$\eqref{Adam c} note the following
general fact. Let $F,G \in \Seqinv_{\Dhat}$, and assume
$\wt{F}, \wt{G}\in \Seqinv_{\Dhat}$ are constructed from $F,G$ by composing with the tensor flips
$\tau_n$ reversing the order of the copies of $\Dhat$
in $\Op(\Dhat^{\odot n})$ for each $n\in \bn$. For any $k \in \bn$
\[ \tau_k \circ \ (F \ast G)_k = \tau_k \circ \left( \sum_{| \boldsymbol{\alpha} | = \{ 1,\ldots,k\} }
F(\alpha_1 \cup \alpha_2 ;k) \QSproj [\alpha_2 ; k] G(\alpha_2 \cup \alpha_3 ; k) \right)\]
\[  = \sum_{| \boldsymbol{\alpha} | = \{ 1,\ldots,k\} } \wt{F}(\wt{\alpha}_1 \cup \wt{\alpha}_2 ;k)
 \QSproj [\wt{\alpha}_2 ; k] \wt{G}(\wt{\alpha}_2 \cup \wt{\alpha}_3 ; k) =
 (\wt{F} \ast \wt{G} )_k,\]
where $\wt{\alpha}_1$ denotes the set $\{k-i: i \in \alpha_1\}$ and $\wt{\alpha}_2, \wt{\alpha}_3$
are defined in an analogous way. The formula above applied to $\up(a)$ and $\up (b)$
(or respectively $\wt{\up}(a)$ and $\wt{\up} (b)$) in place of $F$ and $G$ ends the proof.
\end{proof}

\begin{propn} \label{H}
Let $\varphi \in L\big( \Alg ; \Opstar (\Dhat ) \big)$ and set
$k=l^\varphi$.
\begin{alist}
\item
If
$k$ is
weakly multiplicative then, for all $a,b \in \Alg$,  \label{H a}
$\Dom \big( \varphi (a)^\dagger \big)^* \supset \Ran \QSproj \varphi (b)$ and
\begin{equation} \label{structure}
\varphi (ab) = \varphi (a) \Cou (b) + \Cou (a) \varphi (b) + \varphi (a)
\schur \QSproj \varphi (b) .
\end{equation}
\item
Conversely, if
$\varphi$ satisfies \eqref{structure} then $k$ is weakly multiplicative.
\label{H b}
\end{alist}
\end{propn}

\begin{proof}
Let $\wt{\up} = \wt{\up}^\varphi$.
If
$k$ is weakly multiplicative then, for any $a,b \in \Alg$, using \eqref{generator} and the
Second Fundamental Formula,
\[
\lla \chat, \varphi (ab) \dhat\, \rra =\lla \varphi (a)^\dagger \chat ,
\Cou (b) \dhat \, \rra
+\lla \overline{\epsilon (a)} \chat , \varphi (b) \dhat \, \rra
+\lla \varphi (a)^\dagger \chat ,\QSproj \varphi (b) \dhat \,\rra
\]
for all $c, d \in D$. By sesquilinearity, this implies that
\[
\la \varphi (a)^\dagger \chi, \QSproj \varphi (b) \eta \ra =
\lla \chi ,\big[ \varphi (ab) - \varphi (a) \epsilon (b)
- \epsilon (a) \varphi (b) \big] \eta \rra
\]
for all $\chi, \eta \in \Dhat$, and so (\ref{H a}) holds.

Conversely, if $\varphi$
satisfies
\eqref{structure} then, by Lemma \ref{qgQSDE}, $\wt{\up}$ is $\cH^*_D$-valued
and,  by Proposition \ref{Adam's Lemma 5.4},
$\wt{\up} (a) * \wt{\up} (b) = \wt{\up} (ab)$.
Thus, by Theorem \ref{Lambda isom}, $k$ satisfies
$\Dom (k_t (a)^\dagger)^* \supset \Ran k_t (b)$ and
\[
k_t (a) \schur k_t (b) = \Lambda _t \big( \wt{\up} (a) * \wt{\up} (b)\big)
= \Lambda_t \big( \wt{\up} (ab) \big) = k_t (ab) ,
\]
so $k$ is weakly multiplicative.
\end{proof}

In view of Theorem \ref{Lambda isom}, the following characterisation is obtained
from Lemma \ref{traffic} and Proposition \ref{H}. Its origins date back
to the paper \cite{glo} and the book \cite{schu}. Recall the algebra of
processes defined in Theorem \ref{Lambda isom}, and the remark on block
matrix forms after Lemma \ref{traffic}.

\begin{tw} \label{HH}
Let $k=l^{\varphi}$, where $\varphi \in L(\Alg ; \Opstar (\Dhat ))$,
and suppose that $\Alg$ is a unital $^*$-bialgebra. Then the
following are equivalent\tu{:}
\begin{rlist}
\item
$k$ is unital and $^*$-homomorphic as a map $\Alg \rightarrow (\cP
,\schur\,)$;
\item
$\varphi$ vanishes at $1_{\Alg}$ and satisfies
\begin{equation} \label{alghomstruct}
\varphi (a^*b) = \varphi (a)^* \Cou (b) + \overline{\epsilon (a)}
\varphi (b) + \varphi (a)^* \QSproj \varphi (b) ;
\end{equation}
\item
$\varphi$ has block matrix form
\begin{equation} \label{varphi as block}
\begin{bmatrix}
\gamma & \delta^\dagger \\
\delta & \rho - \iota\circ\epsilon
\end{bmatrix}
\end{equation}
in which $\iota$ is the ampliation $z\mapsto zI_{D}$\tu{;}
\begin{align}
& \rho : \Alg \rightarrow \Opstar (D) \text{ is a unital
*-homomorphism};
\label{Sch T1} \\
& \delta : \Alg \rightarrow |D\ra \text{ is a }
(\rho,\Cou)\text{-derivation\tu{:} }  \notag \\
& \hspace{3cm}
\delta (ab) = \delta (a) \epsilon (b) + \rho (a) \delta (b) ;
\label{Sch T2} \\
& \gamma : \Alg \rightarrow \bc \text{ is linear and satisfies }
\notag \\
&\hspace{2cm}
\gamma(a^*b) =
\overline{\gamma (a)} \Cou (b) + \overline{\Cou (a)} \gamma (b)
+ \delta (a)^*\delta (b).
\label{Sch T3}
\end{align}
\end{rlist}
\end{tw}
\noindent
Following P.-A.\,Meyer (\cite{mey}) we shall refer to such
$(\gamma , \delta , \rho)$
as a $D$-\emph{Sch\"{u}rmann triple} on $\Alg$.

\section{Quantum L\'evy processes}  \label{QLevy}
\setcounter{equation}{0}

In this section Sch\"{u}rmann's theorem on the
reconstruction of a quantum L\'{e}vy process from its `generator' is described.
A new simple proof of the existence of an equivalent realisation of each quantum L\'evy process on a Fock space is given,
based on the results established in previous sections.

Let $\Alg$ be a unital $^*$-bialgebra.

\begin{deft}[\cite{asw}, \cite{schu}]
By a \emph{quantum L\'{e}vy process} on $\Alg$ over a unital
$^*$-algebra-with-state $(\Blg , \omega)$ is meant a family $\big\{
j_{s,t}: \Alg \rightarrow \Blg \, : \; 0 \leq s \leq t \big\}$ of
unital $^*$-homomorphisms satisfying
\begin{enumerate}
\item[(QL1)]
$j_{r,t} = j_{r,s} \star j_{s,t}$  for
$0 \leq r \leq s \leq t$;
\item[(QL2)]
$j_{t,t}(a) = \Cou (a) 1_{\Blg}$ for $t\geq 0$, $ a\in \Alg$;
\item[(QL3)]
$\{ j_{s_i,t_i} (\Alg) : i=1, \ldots n \}$ commute and
\[
\omega \left( \prod^n_{i=1} j_{s_i,t_i} (a_i) \right) = \prod^n_{i=1}
\omega \big( j_{s_i,t_i} (a_i)\big),
\]
whenever $n \in \bn$, the intervals $[s_i,t_i ),\ldots ,[s_n, t_n)$ are disjoint and $a_1, \ldots, a_n \in
\Alg$;
\item[(QL4)]
$\omega \circ j_{s,t} = \omega \circ j_{0,t-s}$ for
$0 \leq s \leq t$;
\item[(QL5)]
$\omega \circ j_{0,t} (a) \stackrel{t \to 0^+}{\longrightarrow}\Cou (a)$
for all $a \in \Alg$.
\end{enumerate}
\end{deft}

Condition (QL1) is known as the \emph{increment property};
the others respectively as the initial condition,
(\emph{tensor}) \emph{independence of
increments},
\emph{time-homogeneity} and \emph{continuity}.
It is immediately verified that
\[
\kappa_t := \omega \circ j_{0,t}
\]
defines a continuous convolution semigroup of states on $\Alg$,
called the \emph{one-dimensional distribution} of the quantum
L\'{e}vy process; its generator is also referred to as the
\emph{generator} of the quantum L\'{e}vy process. For more
information on quantum L\'evy processes on $^*$-bialgebras,
generalisations to free, boolean and monotone case, connections with
Lie algebras and many examples we refer to the book \cite{schu} and
the lecture notes \cite{franz}.

Quantum L\'{e}vy processes $j^i$ on $\Alg$ over $(\Blg^i ,\omega^i)$
$(i=1,2)$ are said to be \emph{equivalent} if they satisfy
\[
\omega^1 \left( \prod^n_{k=1} j^1_{s_k,t_k} (a_k) \right) =
\omega^2 \left( \prod^n_{k=1} j^2_{s_k,t_k} (a_k) \right)
\]
for all $n \in \bn$, disjoint intervals $[s_k,t_k)$ and elements
$a_k$ $(k=1, \ldots, n)$. In view of the axioms (QL1-4) it is clear that
two quantum L\'{e}vy processes are equivalent if and only if their one-dimensional
distributions coincide --- equivalently, if their generators are equal.

Let $\varphi \in L\left(\Alg;\Opstar(\Dhat)\right)$, and
$k=l^{\varphi} \in \QSCdagHc (\Alg ; \Exps_D)$
be unital, real and weakly multiplicative.
Then, setting
\begin{align*}
\Alg^\varphi &= \Lin \Big\{
k_{s_1} (a_1) \schur\sigma_{s_1} \big(k_{s_2-s_1} (a_2) \big)
\schur \cdots \schur\sigma_{s_{n-1}} \big(k_{s_n -s_{n-1}} (a_n ) \big) :
\\
& \hspace{4cm} n \in \bn, 0 \leq s_1 \leq \cdots \leq s_n,
a_1 , \ldots , a_n \in \Alg \Big\},\\
j^\varphi_{s,t} &=
\sigma_s \circ k_{t-s}: \Alg\rightarrow\Alg^\varphi , \text{
and }
\\
\omega^\varphi &= \omega_\Omega |_{\Alg^\varphi} ,
\end{align*}
$\Alg^\varphi$ is a unital $^*$-algebra in the involutive linear
space $\Opdag (\Exps_D)$ with product given by \eqref{bullet},
$\omega^\varphi$ is a state on $\Alg^\varphi$ and it is easily
checked that $j^\varphi$ is a quantum L\'{e}vy process over
$(\Alg^\varphi , \omega^\varphi)$ with generator $\gamma$, where
$\gamma=\varphi_0^0$ (the top-left component of the block matrix
form of $\varphi$). Let us call this type a \emph{Fock space quantum
L\'{e}vy process}.

Note that since a quantum L\'{e}vy process is unital (real) and positive,
its generator $\gamma$ vanishes on $1_{\Alg}$, is real and conditionally
positive (see \eqref{conditionally positive}).

\begin{tw}[\cite{schu}] \label{recon}
Let $\gamma$ be a real,
conditionally
positive linear functional on $\Alg$ vanishing at $1_{\Alg}$.
Then there is a Fock space quantum L\'{e}vy process with generator
$\gamma$.
\end{tw}

\begin{proof}
The proof follows by a GNS-type construction applied
to $\gamma$ viewed as a positive functional on $\Ker \Cou$. Set $D= \Ker \Cou \big/ N$ where
\[
N = \big\{ a \in \Ker \Cou \, : \; \gamma (a^* a) =0 \big\}.
\]
Then $\big([a],[b]\big) \mapsto \gamma (a^*b)$ defines an inner product
on $D$; let $\kil$ be its completion.
Then $\rho (a): [c] \mapsto [ac]$ defines for each $a \in \Alg$ an operator on $D$.
It is obvious that $\rho$ is a unital representation of $\Alg$ on $D$ satisfying
\[
\lla \rho (a)[b],[c] \rra = \lla [b], \rho(a^*) [c] \rra .
\]
Thus $\rho$ is a unital $^*$-homomorphism
$\Alg \rightarrow \Opstar (D)$.
Moreover the linear map
$\delta : a \mapsto  |d(a)\ra $, where $d(a)=[a-\Cou (a) I]$,
is easily seen to be a
$(\rho,\Cou)$-derivation $\Alg\rightarrow |D\ra$ satisfying
\[
\delta (a)^*\delta (b)
= \gamma (a^*b) - \overline{\gamma (a)}
\Cou (b) - \overline{\Cou (a)} \gamma (b) .
\]
Set $k = l^\varphi$, where $\varphi$ is the map $\Alg \rightarrow
\Opstar (\Dhat)$ with block matrix form given by the prescription
\eqref{varphi as block}. Then Theorem \ref{HH} implies that $k$ is
$^*$-homomorphic (i.e.\,real and weakly multiplicative) and unital.
Since $\varphi^0_0 =\gamma$ the result follows.
\end{proof}

\begin{cor}
Every quantum L\'{e}vy process is equivalent to a Fock space quantum
L\'{e}vy process.
\end{cor}

\section{Perturbation}      \label{perturb}
\setcounter{equation}{0} The section is concerned with the
perturbation of $^*$-homomorphic QS convolution cocycles by unitary
Weyl cocycles.

 Consider first the case of the trivial bialgebra $\bc$, and let
 $\varphi \in L(\bc; \Op (\Dhat ))$. Then $\varphi$ and $l^\varphi$ are
determined
by the operator $L := \varphi (1) \in \Op (\Dhat )$ and
the operator Markovian cocycle $X^L := l^\varphi (1) \in \ProcHc (\Exps_D )$ which satisfies
the operator
QS differential equation
\begin{equation} \label{XL QSDE}
\textrm{d}X_t = (L \ot X_t) \, \textrm{d}\Lambda_t , \quad X_0 = I.
\end{equation}
These processes have explicit action on exponential vectors: for any
$f \in \Step_D$
\begin{equation} \label{explicit}
X^L_t \ve (f) = \exp \Big(tz + \int^t_0 \beta (f(s)) \, ds\Big)
\ve \Big( (Rf)_{[0,t[} +d_{[0,t[} \Big)
\end{equation}
where
\[
\begin{bmatrix} z&\beta \\ |d\ra & R-I \end{bmatrix},
\text{ with } z \in \bc, \, d \in \kil , \, \beta \in L(D; \bc ),
\text{ and } R \in \Op (D),
\]
is the block matrix form of $L$.
From either of the above descriptions it is clear that the map
$L \mapsto X^L$ is injective $\Op (\Dhat) \rightarrow \ProcHc (\Exps_D)$.
Moreover if $L \in \Opdag (\Dhat)$ (equivalently,
$R \in \Opdag (D)$ and $\beta = \la c|$ for some
$c \in \kil$) then $X^L \in \ProcHc^{\dagger} (\Exps_D)$ and
$(X^L)^\dagger =X^{L^{\dagger}}$. Similarly, if
 $L \in \Opinv (\Dhat)$ (equivalently, $R \in \Opinv (D)$
and $d \in D$) then $X^L_t \in \Opinv (\Exps_D)$ for each $t$. If
$L\in \Op (\Dhat)$ and $M \in \Opinv (\Dhat)$ then,  by the explicit action
\eqref{explicit},
\begin{equation} \label{XLXM}
X^L X^M = X^{L \Odiamond M}
\end{equation}
where
\begin{equation} \label{diamond}
L \Odiamond M := L + M+ L \QSproj M
\end{equation}
By the above injectivity $\big( \Opinv (\Dhat), \Odiamond \big)$
is a semigroup with identity $0$; clearly $\big( \Opstar (\Dhat),\Odiamond \big)$
is an involutive semigroup: $(L \Odiamond M)^\dagger = M^\dagger \Odiamond L^\dagger$.
Note that these identities contain the Weyl commutation relations.

The above formula implies that for $L\in \Opstar (\Dhat)$
\[
X^L \text{ is isometric} \iff L^{\dag} \Odiamond L =0, \text{ whereas }
X^L \text{ is coisometric} \iff L \Odiamond L^{\dag} =0,
\]
cf.\ analogous characterisations described in~\cite{lwptrf}.

In the next proposition, \eqref{diamond} is extended by left and right actions
of (parts of) $\Opstar (\Dhat )$  on $L\big(\Coalg ; \Op (\Dhat)\big)$,
for a coalgebra $\Coalg$.

\begin{propn}
Let $\varphi \in L(\Coalg ; \Op (\Dhat ))$ and let
$L,M \in \Op (\Dhat )$.
\begin{alist}
\item
If $\varphi \in L \left( \Coalg ; \Opdag (\Dhat) \right)$ and
$ M \in \Opinv (\Dhat)$ then
\[
 l^\varphi (\, \cdot \, ) X^M = l^{\varphi \Odiamond M} ,
\]
\item
If
$ L \in \Opdag (\Dhat)$ and $ \varphi \in L\left(\Coalg ; \Opinv (\Dhat)
\right)$
then
\[
 l^\varphi_t (\Coalg )\subset \Dom (X^L_t)^{\dagger *} \text{ and }
X^L \schur l^\varphi (\, \cdot \, )= l^{L\Odiamond \varphi},
\]
\end{alist}
where  for each $a \in \Coalg$
\begin{equation*}
(\varphi \Odiamond M)(a) := \varphi (a) (I + \QSproj M) + \Cou(a) M
\end{equation*}
and
\begin{equation*}
(L \Odiamond \varphi ) (a) := (I+ L\QSproj ) \varphi (a)
+ \Cou (a) L.
\end{equation*}

\end{propn}

\begin{proof}
These follow easily from the two Fundamental Formulae \eqref{FFF} and \eqref{SFF}.
\end{proof}

\noindent
The above formulae extend \eqref{diamond} -- this may be checked by setting
$\Coalg = \bc$ and $a=1$.

Let $\varphi \in L \big(\Coalg ; \Opstar (\Dhat) \big)$ and
$L_1, L_2  \in \Opstar (\Dhat)$. Then their block matrix forms
(see \eqref{block}) are respectively given by
\[
\begin{bmatrix}
\gamma & \alpha \\
\chi & \nu - \iota\circ\epsilon
\end{bmatrix}
\quad \text{and} \quad
\begin{bmatrix} z_i&\la c_i| \\ |d_i\ra & R_i-I \end{bmatrix},
\]
where
$z_i \in \bc$, $c_i \in D$, $R_i \in \Opstar (D)$,
and
\begin{align*}
\lefteqn{(L^\dagger_1 \Odiamond \varphi \Odiamond L_2 )(a)} \\
&= (I+\QSproj L_1)^\dagger \varphi (a) (I+\QSproj L_2 )
+ \Cou (a) L^\dagger_1 \Odiamond L_2 \\
&= \begin{bmatrix}
\wt{\gamma} (a) &
\big( \alpha (a) + \la d_1|\nu (a)\big) R_2
+ \Cou (a) \la c_2|
\\[0.2cm]
R^\dagger_1 \big( \chi (a) +\nu (a) |d_2\ra\big)
+ \Cou (a) |c_1\ra  &
R^\dagger_1  \nu (a) R_2- \Cou (a) I \end{bmatrix}
\end{align*}
where
\[
\wt{\gamma} (a) =
\gamma (a)
+ (z^*_1+z_2)\epsilon (a)
+  \alpha (a) |d_2 \ra + \la d_1| \chi (a)
+ \big\la d_1, \nu (a) d_2 \big\ra .
\]

Let now $\Alg$ be a unital $^*$-bialgebra and consider conjugation
by a single map $L \in \Opstar (\Dhat)$:
\[
\wt{\varphi} = L^{\dag} \Odiamond \varphi \Odiamond L,
\]
where $\varphi\in L\big(\Alg ; \Opstar (\Dhat)\big)$.
It is easily checked that if
$l^{\varphi}$ is real then $l^{\wt{\varphi}}$ is real;
if $l^{\varphi}$ is unital
then
\[
l^{\wt{\varphi}} \text{ is unital } \iff X^L \text{ is isometric};
\]
and if $l^{\varphi}$ is weakly multiplicative then
$l^{\wt{\varphi}}$ is weakly multiplicative if and only if  for all $a \in \Alg$
\[
  (\QSproj L +I)^*
\big(\QSproj \varphi(a) + \epsilon(a)I\big)^\dagger L \Odiamond L^\dagger
\big(\QSproj \varphi(a) + \epsilon(a) I\big) (\QSproj L +I)
=0.
\]

Therefore, considering perturbations by unitary (Weyl) cocycles, one
obtains the action of the Euclidean group of $D$ on Sch\"urmann
triples associated with unital $^*$-homomorphic QS convolution
cocycles on $\Alg$ (cf.\ \cite{franz}). This action has a simple
matricial description: if
$$
L=
\begin{bmatrix}
i\mu - \frac{1}{2} \|v\|^2 & - \langle v |V \\
    |v \rangle & V - I
\end{bmatrix},
$$
where $\mu \in \br, v \in D$ and $V \in \Opstar (D)$ is unitary,
then   for all $a \in \Alg$
$$
\wt{\varphi} (a)
 =
\begin{bmatrix} 1  & \langle v|
 \\[0.2cm]
0 & V^* \end{bmatrix}
\varphi (a)
\begin{bmatrix} 1  & 0
 \\[0.2cm]
|v\rangle & V \end{bmatrix}  .
$$
Thus if
$$
\varphi =
\begin{bmatrix}
\gamma  & \delta^{\dag} \\
\delta & \rho - \iota\circ\Cou
\end{bmatrix}
$$
then  for all $a \in \Alg$
\begin{align*}
\wt{\gamma} (a)& =
\gamma (a) + \delta^{\dag}(a)|v\ra
+ \la v|\delta (a) + \big\la v , \big( \rho (a) - \Cou (a) I
\big) v \big\ra , \\
\wt{\delta} (a) &=
V^* \Big( \delta (a) + \big( \rho (a) - \Cou (a) I\big) |v\ra \Big)
\text{ and } \\
\wt{\rho} (a) &= V^* \rho (a) V.
\end{align*}

Notice that the part of the action determined by $V$ is trivial
in the sense that only a unitary transformation of the Sch\"urmann triple
$(\lambda,\delta,\rho)$ leaving $\lambda$ invariant is effected, so that
the
perturbed quantum L\'{e}vy process $l^{\wt{\varphi}}$ is equivalent to the
unperturbed one $l^\varphi$. For nonzero $v$
the perturbation still does not change the characteristics of
the quantum L\'{e}vy process, in the sense that  \emph{Gaussian} processes remain Gaussian
and the same is
true for \emph{Poisson} and \emph{drift} processes (for relevant definitions see \cite{franz}).

\section{Opposite QS convolution cocycles} \label{oppQSccocycl}
\setcounter{equation}{0}

This section is concerned with the opposite counterpart of the notion of a QS convolution cocycle introduced in Section \ref{QSccocycl}.

\begin{deft}
An opposite quantum stochastic convolution cocycle
(on $\Coalg$ with domain $\Exps_D$) is a process
$l\in\Proc (\Coalg; \Exps_D)$ such that, for all $s,t \geq 0$,
\begin{equation}
\label{oppconvinc}
l_{s+t} = (\sigma_s \circ l_t) \star l_s
\end{equation}
 and for all $a\in \Coalg$
\[
l_0 (a) = \Cou(a) I_{\Exps_D}.\]
\end{deft}

Family of all opposite QS convolution cocycles is denoted by $\QSC^{\tu{opp}} (\Coalg; \Exps_D)$; the
notation adorned with subscripts and superscript $^\dagger$ according to our convention.

\begin{rem}
The opposite convolution increment formula makes use of the identification
$\Focktots \cong \Fockstot \ot \Focktos$.
The weak description of the opposite convolution increment property takes the form:
\begin{equation}
\mathrm{e}^{-\la f_{[0,t)}, g_{[0,t)} \ra } \lla \ve (f_{[0,t)}), l_t (a)
\ve (g_{[0,t)}) \rra = \prod^{n-1}_{i=0} \lambda^{f(t_i),g(t_i)}_{t_{i+1}-t_i}
(a_{(n-i)})
\label{oppsemigroup decomp}
\end{equation}
(compare with the formula \eqref{semigroup decomp} of Lemma \ref{M}).

Opposite cocycles also have associated convolution semigroups, given by
formula \eqref{assoc semigroups}; two opposite QS convolution cocycles with
identical convolution semigroups are equal.
\end{rem}

Opposite QS convolution cocycles arise as solutions of the opposite coalgebraic QS
differential equations of the form:
\begin{equation} \label{opqgQSDE}
dk_t = \varphi \star k_t  \,d\Lambda_t , \; \;
k_0 = \iota\circ\epsilon,
\end{equation}
where $\varphi \in L\big(\Coalg ; \Op (\Dhat )\big)$.

Below we formulate the relevant opposite versions of basic theorems
proved earlier for QS convolution cocycles.

\begin{tw} \label{opG}
Let $\varphi \in L\left(\Coalg ; \Op (\Dhat )\right)$ and set
$\up = \up^\varphi$. Then the process $^{\varphi}l := \Lambda \circ \up$
strongly satisfies the opposite coalgebra quantum stochastic differential equation
\eqref{opqgQSDE}. It is a unique weak solution of \eqref{opqgQSDE}.
\end{tw}

\begin{tw} \label{opQ}
Let $k \in \Proc (\Coalg ;\Exps_D )$. Then the following are
equivalent \tu{:}
\begin{rlist}
\item $k \in {\QSCdagHc}^{\tu{opp}} (\Coalg ; \Exps_D )$;
\item
$k =\, ^{\varphi}l $ for some $\varphi \in L\big(\Coalg ; \Opdag
(\Dhat)\big)$.
\end{rlist}
\end{tw}

\begin{tw} \label{opHH}
Let $k=\,^{\varphi}l$, where $\varphi \in L(\Alg ; \Opstar (\Dhat
))$, and suppose that $\Alg$ is a unital $^*$-bialgebra. Then the
following are equivalent\tu{:}
\begin{rlist}
\item
$k$ is unital and $^*$-homomorphic as a map $\Alg \rightarrow (\cP
,\schur\,)$.
\item
$\varphi$ vanishes at $1_{\Alg}$ and satisfies
\[
\varphi (a^*b) = \varphi (a)^* \Cou (b) + \overline{\epsilon (a)}
\varphi (b) + \varphi (a)^* \QSproj \varphi (b) . \]
\end{rlist}
\end{tw}

The proofs of the above may be conducted exactly along the same
lines as for usual QS convolution cocycles - whenever it was
convenient, the arguments in the proofs were formulated to make this
transparent (for example by working with both $\up$ and $\wt{\up}$).
Alternatively, one may exploit the correspondence described below.

There is a bijective correspondence between QS convolution cocycles
and their opposite counterparts. It may be implemented either by the
time-reversal operation or by passage to the opposite coalgebra. To
formulate this precisely we need to introduce some further notation.
For each $t\geq 0$ define the following time reversal operator on
$L^2(\br_+; \kil)$:
\[( r_t (f)) (u)=  \begin{cases} f(u-t) & \text{ if } u \leq t\\ f(u) & \text{ if } u > t \end{cases}, \;\;\;\;\; f \in L^2(\br_+; \kil),\,  u \geq 0.\]
This in turn may be second-quantised to the unitary selfadjoint
operator $R_t \in B(\Fock)$ whose action on exponential vectors is
given by
\[ R_t (\ve(f)) = \ve (r_t f).\]
The latter provides, by conjugation, time reversal on the level of
operators in $\Op(\Exps_D)$:
\[ \rho_t(Z) = R_t Z R_t, \; \; Z \in \Op(\Exps_D).\]

Recall also that if $\Coalg$ is a coalgebra with coproduct $\Com$, the \emph{opposite coalgebra} $\Coalg^{\text{opp}}$ is the same vector space as $\Coalg$
equipped with the same counit and with the coproduct $\tau \circ \Com$, where $\tau$ denotes the tensor flip on $\Coalg \odot \Coalg$.

\begin{propn}
Let $l \in \Proc(\Coalg; \Exps_D )$. Then the following are
equivalent \tu{:}
\begin{rlist}
\item  $l$ is an opposite QS convolution cocycle on $\Coalg$;
\item  the process $\wt{l} \in \Proc(\Coalg^{\text{opp}}; \Exps_D )$  given by
 \[ \wt{l}_t = l_t, \;\;\; t \geq 0,\]
   is a QS convolution cocycle on $\Coalg^{\text{opp}}$;
\item the process $\wt{l} \in \Proc(\Coalg; \Exps_D )$ defined by
     \[ \wt{l}_t = \rho_t \circ l_t, \;\;\; t \geq 0,\]
  is a QS convolution cocycle on $\Coalg$.
\end{rlist}
\end{propn}

\begin{proof}
All the equivalences follow directly from the relevant semigroup decomposition formulas \eqref{semigroup decomp} and \eqref{oppsemigroup decomp}.
\end{proof}

In particular time-reversal (or passage to the opposite coalgebra)
exchanges $l^{\varphi}$ with $^{\varphi}l$.

Finally note that in \cite{LSqscc1} we actually worked with
opposite cocycles (although it was not clarified there).

 \chapter{$C^*$-algebraic case}

\newpage

This chapter is concerned with quantum stochastic convolution
cocycles on an operator space coalgebra, in the second part
specialised to a $C^*$-hyperbialgebra or a $C^*$-bialgebra. Many
results in this part mirror the ones obtained in the purely
algebraic case. Technical conditions under which theorems hold are
however usually different, and although some proofs are similar, it
is not possible to conduct them precisely so that the argument is
valid for both cases.

The operator space coalgebraic structure, defined by analogy with the purely algebraic case, allows for
defining the convolution
product. In the topological context it becomes important that the convolution
may be transformed via the so-called $R$-map into a composition operation
 preserving all the relevant continuity properties.
This corresponds to the procedure, familiar from classical
probability, of transforming the convolution of measures into the
composition of Markov operators. On the level of cocycles the
$R$-map transforms a QS convolution cocycle, with its associated
convolution semigroups, into a standard QS cocycle, with its
associated semigroups.

Operator-space theoretic QS convolution cocycles, as in the purely
algebraic case, are obtained by solving coalgebraic QS differential
equations. This time however to prove the existence and uniqueness
theorems for solutions we need to assume complete boundedeness of
coefficients (in \cite{LSqsde} it is shown that in fact it is enough
to assume that the coefficient has completely bounded columns). The
methods used here extend those of \cite{lwex}, allowing nontrivial
initial conditions. As in general there is no guarantee that the
solutions of a given coalgebraic QS differential equation will be
completely bounded, in this topological context we need to
distinguish the class of weak QS convolution cocycles. These satisfy
the cocycle relation only on the level of relevant matrix elements.

A converse of the fact that each solution of a coalgebraic QS
differential equation is a weak QS convolution cocycle may be
obtained when the initial object is a $C^*$-hyperbialgebra, so that
one can exploit the order structure. Using the known results from
the theory of standard QS cocycles (\cite{lp}, \cite{lwptrf}) and
the $R$-map described in the previous paragraph we prove that every
Markov-regular, completely positive and contractive QS convolution
cocycle on a $C^*$-hyperbialgebra satisfies a coalgebraic QS
differential equation. We also characterise the general form of
stochastic generators yielding such cocycles. This leads further to
the characterisation of stochastic generators of $^*$-homomorphic QS
convolution cocycles on a $C^*$-bialgebra $\alg$ in terms of
structure maps on $\alg$ (which can again be equivalently described
via Sch\"urmann triples). The interesting feature here is that `the
algebra determines the analysis' - every structure map on $\alg$ is
necessarily completely bounded, and even inner. The proof of this
fact is based on versions of well-known results of S.\,Sakai,
J.R.\,Ringrose and E.\,Christensen on continuity and innerness
properties of derivations for the case of
$(\pi_1,\pi_2)$-derivations (the proofs are provided in Appendix A).
The innerness of structure maps on $\alg$ may be viewed as a
noncommutative counterpart to the fact that every classical L\'evy
process on a topological group which has a bounded generator must be
a compound Poisson process.

It is possible to axiomatise quantum L\'evy processes on
$C^*$-bialgebras, either using solely the concept of distributions,
or exploiting the language of Arveson's theory of product systems.
Accordingly, we propose definitions of a weak quantum L\'evy process
and a product system quantum L\'evy process. For each of them a
topological version of the Sch\"urmann Reconstruction Theorem
remains valid.

The explicit characterisations of the stochastic generators enable us to prove two types of
dilation theorems for completely positive and contractive QS convolution cocycles
on a $C^*$-bialgebra, corresponding to the dilations obtained for standard QS cocycles
in \cite{dilate} and \cite{Stine}.

After the general theory has been presented, we describe the basic examples of commutative, cocommutative
and \emph{genuinely quantum} $C^*$-bialgebras.
Different perspectives on QS convolution cocycles or their stochastic generators
in each of the cases is offered. In particular in the context of full compact quantum groups
purely algebraic and operator space theoretic cocycles are shown to coexist on the same underlying space.
This supports the view (explicitly described in the
expository paper \cite{LSbedlewo}) that despite the obvious differences both categories of QS convolution cocycles
may nevertheless usefully be seen from a common vantage point.

It would be interesting to obtain extensions of the results of this
chapter to multiplier $C^*$-bialgebras (locally compact quantum
semigroups) - a short explanation of the technical difficulties
related to such a project is provided in the last section.

Finally we would like to stress that, in contrast to \cite{LSqscc2}, all operator space coalgebras,
$C^*$-hyperbialgebras and $C^*$-bialgebras are assumed to be concrete. For the latter two we mean
by that explicitly that they are represented in a faithful and nondegenerate way on a Hilbert space.

\section[OS coalgebras and
convolution semigroups]
{Operator space coalgebras and convolution semigroups}    \label{OS coalgebras}
\setcounter{equation}{0}

In this section we present the definition of an operator space coalgebra, establish basic facts concerning
the convolution product provided by the coalgebraic structure and define the $R$-map facilitating the traffic
between the convolution and composition operations.

\begin{deft}
An operator space $\coalg$ is an operator space (OS) coalgebra if there are completely
contractive maps
 $\Com : \coalg \rightarrow \coalg \ot \coalg$ and
$\Cou : \coalg \rightarrow \bc$, called the \emph{coproduct} and
\emph{counit} respectively, enjoying the coassociativity and the counit
property, namely
\begin{align}
(\id \ot \Com ) \circ \Com &=
( \Com \ot \id ) \circ \Com,
 \label{Ccoassoc} \\
(\id \ot \Cou )\circ \Com &= ( \Cou \ot \id ) \circ
\Com = \id.  \label{Ccounit}
\end{align}
\end{deft}

The formula \eqref{Com n} defines again a map $\Com_n$, and this time
$\Com_n : \coalg^{\ot n} \to \coalg^{\ot n+1}$ - recall that
$\Com_0:= \id_{\coalg}$. As we no longer can expect that
$\Com (a)$ ($a\in \coalg$) is a finite sum of simple tensors, we generally avoid Sweedler
notation in this chapter.

\begin{deft} \label{Cbialg} A unital $C^*$-algebra $\alg$ is called a $C^*$-hyperbialgebra, if
it is an OS coalgebra, the counit is a character (unital
multiplicative functional) and the coproduct is unital and
completely positive. A $C^*$-hyperbialgebra is called a
$C^*$-bialgebra if its coproduct is multiplicative (so
$^*$-homomorphic).
\end{deft}

\begin{rem}
One can also introduce an apparently natural category of operator
system coalgebras, but it will not be of much use here. Some authors
reserve the name $C^*$-bialgebra for a not-necessarily-unital
$C^*$-algebra $\alg$ with the coproduct taking values in $M(\alg \ot
\alg)$ (the multiplier algebra of $\alg \ot \alg$); this is relevant
for the considerations in Section \ref{multiplier}.
\end{rem}

The motivating examples of $C^*$-bialgebras and $C^*$-hyperbialgebras come
respectively from the theory of compact quantum groups (\cite{CMP}, \cite{wor2})
and hypergroups (\cite{ChaV}).

Let us introduce a modification of the convolution  product $\star$ of Chapter 2.
If $\coalg$ is an OS coalgebra, $\Vil_1, \Vil_2$ - operator spaces and
$\varphi_1 \in CB (\coalg;\Vil_1)$,  $\varphi_2 \in CB (\coalg;\Vil_2)$, define
\[ \varphi_1 \star \varphi_2 = (\varphi_1 \ot \varphi_2) \circ \Com : \coalg \to \Vil_1 \ot \Vil_2.\]

It is easily seen that $\star$ is associative and enjoys submultiplicativity and unital properties:
\begin{align}
&\varphi_1 \star \varphi_2 \star \varphi_3 = (\varphi_1 \ot \varphi_2 \ot \varphi_3 )
   \circ \Com_2 \label{star associativity} \\
& \| \varphi_1 \star \varphi_2 \|_{\cb} \leq
\| \varphi_1 \|_{\cb} \| \varphi_2 \|_{\cb}, \text{ and } \notag \\
& \Cou \star \varphi = \varphi = \varphi \star \Cou . \label{conv unit}
\end{align}

In particular, $(\coalg^* , \star)$ is a unital Banach algebra.
The following observation will be used further: if $\Vil$ is an operator space,
$\varphi \in CB(\coalg ; \Vil)$ and $n \in \bn$, then
\begin{equation} \label{phi star n}
\varphi^{\star n} = \varphi^{\ot n} \circ \Com_{n-1} .
\end{equation}
Define also $\varphi^{\star 0} := \Cou$, which is consistent with \eqref{conv unit}.

\begin{deft}
A convolution semigroup of functionals
(CSF, for short) on an OS coalgebra $\coalg$ is a family $\{\kappa_s:s \geq 0\}$ of
bounded linear functionals on $\coalg$ satisfying the following conditions:
\begin{equation} \label{CCSFincr}
\kappa_{s+t} = \kappa_s \star \kappa_t , \;\;\; s,t \geq 0, \end{equation}
\begin{equation} \label{CCSFinit}
\kappa_0 = \Cou.\end{equation}
If $\| \kappa_s - \Cou \| \to 0$ as $s \to 0^+$
then $\{\kappa_s:s \geq 0\}$ is said to be norm continuous (it is then abbreviated to CCSF).
\end{deft}

The notions
of cocommutative coalgebras and idempotent functionals introduced in Section \ref{algandconv}
have obvious counterparts in this context.

\subsection*{Convolution and the $R$-map}
To analyse existence and properties of generators of CCSFs, and also for the
investigations in Section \ref{standvsconv} we need to describe the properties of
a certain map `linearising' the convolution.

Given an operator space coalgebra $\coalg$, each operator space $\Vil$
determines maps
\begin{align*}
&R_{\Vil} : CB(\coalg;\Vil) \to CB(\coalg ;\coalg \ot \Vil),
     \quad \varphi \mapsto (\id_{\coalg} \ot \varphi) \circ \Com = \id_{\coalg} \star \varphi; \\
&E_{\Vil} :CB(\coalg;\coalg \ot \Vil) \to CB(\coalg ;\Vil),
     \quad \phi \mapsto (\Cou \ot \id_{\Vil} )\circ \phi .
\end{align*}
Thus the action of $R_{\Vil}$ is given by the convolution with the identity map on $\coalg$,
putting the argument on the right, and the action of $E_{\Vil}$ is given by the composition
with the tensor-extension of the counit.
Symbols $R_{\bc}$ and $E_{\bc}$ will be respectively abbreviated to $R$ and $E$.

The following proposition is a straightforward consequence of
the definitions of $R$ and $E$.

\begin{propn} \label{thm: R and E}
Let $\coalg$ be an OS coalgebra, and let $\Vil_1$, $\Vil_2$,
$\Vil$ be operator spaces.
\begin{alist}
\item
$R_{\Vil}$ and $E_{\Vil}$ are complete isometries satisfying
\[
E_{\Vil} \circ R_{\Vil} = \id_{CB(\coalg ;\Vil)} .
\]
\item
If $\varphi_1 \in CB(\coalg ;\Vil_1)$ and $\varphi_2 \in CB(\coalg ; \Vil_2)$
then
\[
R_{\Vil_1 \ot \Vil_2} (\varphi_1 \star \varphi_2)
= \left( R_{\Vil_1} \varphi_1 \ot \id_{\Vil_2} \right) \circ R_{\Vil_2} \varphi_2 .
\]
\end{alist}
\end{propn}

Write $CB^{\Com}(\coalg ;\coalg \ot \Vil)$ for $\Ran R_{\Vil}$.

\begin{cor} \label{cbisom}
For each operator space $\Vil$, $R_{\Vil}$ is a complete isometry of operator spaces
\[
CB(\coalg ; \Vil) \cong CB^{\Com}(\coalg ;\coalg \ot \Vil) .
\]
Moreover, $R$ is an isometric isomorphism of unital Banach algebras
\[
(\coalg^* , \star) \cong \big( CB^{\Com} (\coalg ), \circ \big) .
\]
\end{cor}

A further interesting consequence is the following

\begin{cor} \label{Rtraffic}
In $CB^{\Com} (\coalg ; \coalg \ot M_n )$,
\[
\| \phi\|_{\cb} = \| \phi^{(n)} \| .
\]
\end{cor}

\begin{proof}
Let $\phi \in CB^{\Com} (\coalg ; \coalg \ot M_n )$, say $\phi = R_{M_n} \varphi$
for some $\varphi\in CB(\coalg; M_n)$. Then
\[
\| \phi \|_{\cb} = \| \varphi \|_{\cb} = \| \varphi^{(n)} \| = \| (\Cou \ot \id_{M_n(\coalg)}) \phi^{(n)} \|
\leq \| \phi^{(n)} \| ;
\]
the result follows.
\end{proof}

In particular, in $CB^{\Com} (\coalg)$ the completely bounded norm coincides with the bounded
operator norm. As a result $CB^{\Com} (\coalg)$ is a closed subspace of $B(\coalg)$.
The next proposition collects the structure-preserving properties of the map $R_{\Vil}$ and
$E_{\Vil}$ under a number of relevant assumptions on $\coalg$ and $\Vil$.

\begin{propn} \label{Requiv}
Let $\coalg$ be an OS coalgebra and $\Vil$ an operator space, let
$\varphi \in CB(\coalg ;\Vil)$ and
$\phi = R_{\Vil}\varphi \in CB (\coalg ; \coalg \ot \Vil)$.
\begin{alist}
\item
$\phi$ is completely contractive if and only if $\varphi$ is.
\item
If $\coalg$ is a $C^*$-hyperbialgebra and $\Vil$ is an operator system
then $\phi$ is real (respectively, completely positive, or unital) if and
only if $\varphi$ is.
\item
If $\coalg$ is a $C^*$-bialgebra and $\Vil$ is a $C^*$-algebra then $\phi$
is multiplicative if and only if $\varphi$ is.
\end{alist}
\end{propn}

The following proposition paves the way for analysing generators of convolution
semigroups of functionals on OS coalgebras.

\begin{propn}
Let $\coalg$ be an operator space coalgebra. The map
$\kappa:= \{\kappa_t: t \geq 0\}  \mapsto P := \{(R \kappa_t)_{t \geq 0}\}$
is a bijection from the set of CSFs on $\coalg$ to the set of one-parameter semigroups in
$CB^{\Com} (\coalg)$. Moreover, the conditions in (a) below are equivalent, and so are the conditions in (b):
\begin{alist}
\item
\begin{rlist}
\item
\label{lim} $\lim_{t \to 0} \kappa_t (a) = \Cou (a)$ for all $a \in \coalg$;
\item
\label{c0 semigroup} $P$ is a $c_0$-semigroup on $\coalg$.
\end{rlist}
\item  \label{part bb}
\begin{rlist}
\item $\kappa$ is norm continuous in $t$;
\item $P$ is norm continuous in $t$;
\item $P$ is $\cb$-norm continuous in $t$;
\item the generator of $P$ is completely bounded.
\end{rlist}
\end{alist}
\end{propn}

\begin{proof}
The first part follows from Corollary \ref{cbisom}.

(a) Since $\Cou \circ P_t = \kappa_t$, (ii) implies (i). Suppose
therefore that (i) holds. Then, for any $\lambda \in \coalg^*$,
\[
\lambda \circ P_t = \kappa_t \circ (\lambda \ot \id_{\coalg}) \circ \Com\]
 and
\[ \Cou \circ (\lambda \ot \id_{\coalg} ) \circ  \Com = \lambda ,\]
so  $P_t a$ converges to $a$ weakly as $t \to 0^+$ for all $a \in \coalg$.
But this implies that $P$ is
strongly continuous (\cite{Davies}, Proposition 1.23) and thus a
$c_0$-semigroup.

(b) By Corollary \ref{Rtraffic}
\[
\| P_t -\id_{\coalg} \|_{\cb} = \| \kappa_t - \Cou \| = \| P_t -
\id_{\coalg} \| ,
\]
and so (\ref{part bb}) follows.
\end{proof}

Thus each norm-continuous CSF $\{\kappa_t: t \geq 0\}$ on $\coalg$ has a \emph{generator}:
\[
\gamma := \lim_{t \to 0^+} \frac{(\kappa_t - \Cou)}{t}\]
from which the CSF may be recovered:
\[
\kappa_t = \exp_{\star} t \gamma:= \sum_{n \geq 0} (n!)^{-1} t^n \gamma^{\star n} \]
($t \geq 0$), where $\gamma_0:=\Cou$.

The corresponding one-parameter semigroup on $\coalg$ has a completely bounded generator:
$R \kappa_t = \mathrm{e}^{t\tau}$, $t \geq 0$ , where $\tau = R \gamma \in CB(\coalg)$.
Remarks after Definition \ref{modcoalg} apply also here.

\section[QS convolution cocycles]{QS convolution cocycles and standard QS cocycles}
 \label{standvsconv}
\setcounter{equation}{0}
Here the main object of the thesis is defined in the operator space theoretic context, and the connection
with the standard theory of QS cocycles is made.

\subsection*{QS convolution cocycles on OS coalgebras}

\begin{deft}
A \emph{quantum stochastic convolution cocycle}
(on $\coalg$ with domain $\Exps$) is a process
$l\in\Proc_{\cb} (\coalg; \Exps)$ such that, for all $s,t \geq 0$,
\begin{equation}
\label{Cconvinc}
l_{s+t} = l_s \star (\sigma_s \circ l_t)
\end{equation}
 and for all $a\in \coalg$
\begin{equation} \label{Ccouinit}
l_0 (a) = \Cou(a) I_{\Fock}.
\end{equation}
\end{deft}

Family of all such QS convolution cocycles is denoted by $\QSC (\coalg; \Exps)$; the
notation adorned with subscripts and superscript $^\dagger$ according to our convention.
Remark \ref{cocprop} remains valid, and QS convolution cocycles on $\coalg$ have
associated CSFs, as the following lemma shows.

\begin{lemma} \label{CM}
Let $l \in \Proc_{\cb} (\coalg ; \Exps)$ be a
 convolution increment process. For each $f,g \in \Step$,
\begin{equation} \label{Csemigroup decomp}
\mathrm{e}^{-\la f_{[0,t)}, g_{[0,t)} \ra } \lla \ve (f_{[0,t)}), l_t (a)
\ve (g_{[0,t)}) \rra = \left( \lambda^{f(t_0),g(t_0)}_{t_{1}-t_0} \star \cdots \star
 \lambda^{f(t_{n-1}),g(t_{n-1})}_{t_{n}-t_{n-1}} \right) (a)
\end{equation}
where $0 = t_0 \leq t_1 \leq \cdots \leq t_n =t$ contains the discontinuities
of $f_{[0,t)}$ and $g_{[0,t)}$, and
\begin{equation} \label{Cassoc semigroups}
\lambda^{c,d}_t := \mathrm{e}^{-t\la c,d\ra} \lla \ve (c_{[0,t)}), l_t (\,\cdot\,)
\ve (d_{[0,t)} ) \rra, \;\; c,d \in \kil.
\end{equation}
\end{lemma}

\begin{proof}
As for purely algebraic cocycles, the identity \eqref{Csemigroup decomp} results from
repeated application
of the cocycle relation \eqref{Cconvinc}.
\end{proof}

In the topological case there is a need to distinguish QS convolution cocycles
from the class of processes satisfying \eqref{Csemigroup decomp}. Recall the notion of
weakly bounded processes introduced in Section \eqref{Focksection}.

\begin{deft}
A process $l \in \Proc_{\wb} (\coalg ; \Exps)$ is a weak convolution
increment process if it satisfies formula \eqref{Csemigroup decomp}; it is
a weak QS convolution cocycle if in addition $l_0 (a) = \Cou(a) I_{\Fock}$ for all
$a \in \coalg$.
\end{deft}

\begin{rem}
Note that weak boundedness of $l$ guarantees boundedness (so also complete
boundedness) of all functionals $\lambda^{c,d}$ and the right hand side
of formula \eqref{Csemigroup decomp} makes sense. Lemma \ref{CM} says that
each QS convolution cocycle is a weak QS convolution cocycle. It is easy to observe that
if a weak QS convolution cocycle is completely bounded, it is a QS convolution cocycle.
\end{rem}

\begin{cor} \label{CN}
Let $l \in \Proc_{\wb} (\coalg ; \Exps)$ be a weak QS convolution cocycle, $c,d \in \kil$. Then
\eqref{Cassoc semigroups} defines a convolution semigroup of functionals
on $\coalg$ : $ (\lambda^{c,d}_t )_{t \geq 0}$.
\end{cor}

\begin{deft}    \label{Cdefsemig}
Semigroups defined by $\eqref{Cassoc semigroups}$ (for all pairs $c,d \in \kil$) are called
associated convolution semigroups of the cocycle $l$. The semigroup
$\lambda^{0,0}$ is called the Markov semigroup of the cocycle $l$.
\end{deft}

It is clear that two weak QS convolution cocycles with identical
associated convolution semigroups are equal. For further analysis we need one more definition.

\begin{deft}    \label{cMarkreg}
A weak QS convolution cocycle is Markov-regular if its Markov semigroup is norm-continuous.
\end{deft}

\subsection*{Standard QS cocycles}

Standard QS cocycles have been analysed in numerous papers. Here we follow the
point of view of \cite{lect} (see also \cite{lwjfa} and references therein); the
terminology is however different, as we intend to treat completely bounded processes
as the objects of primary interest. None of the results in this section is new,
they are recalled here as they will be used to prove many important properties
of QS convolution cocycles.

Let $\Vil\subset B(\hil)$ be a concrete operator space. Recall the matrix spaces introduced
in Section \ref{Ospaces}. Whenever $k \in \Proc_{\cb} (\Vil;\Vil, \Exps)$,
$s,t \geq 0$, we obtain the map $\widehat{\sigma}_s \circ k_t:\Vil \mapsto M_{\Fockafs} (\Vil)$,
where $\widehat{\sigma}_s$ is a natural extension of the shift endomorphism to $M_{\Fock}(\Vil)$ and
$k_s$ may be equivalently viewed as a map from $\Vil$ to $M_{\Focktos} (\Vil)$. Denote by
$\widehat{k}_s$ the matricial extension of $k_s$ to $M_{\Fockafs} (\Vil)$.
The functorial property \eqref{functprop} and the tensor factorisation \eqref{factor}
of the Fock space imply that the following definition is consistent.

\begin{deft}
A process $k \in \Proc_{\cb} (\Vil;\Vil, \Exps)$ is a standard QS cocycle (on $\Vil$ with
domain $\Exps$) if
\begin{equation}
\label{standinc}
k_{s+t} = \widehat{k}_s \circ \widehat{\sigma}_s \circ k_t,
\end{equation}
 and for all $v \in \Vil$
\begin{equation} \label{standinit}
k_0 (v)   = v \ot  I_{\Fock}.
\end{equation}
\end{deft}

The associated semigroups of a standard QS cocycle are semigroups acting on $\Vil$, as the next lemma shows:

\begin{lemma} \label{standCM}
Let $k \in \Proc_{\cb} (\Vil ;\Vil, \Fock)$ be a standard QS cocycle. For each $f,g \in \Step$,
\begin{equation} \label{stan Csemig decomp}
\mathrm{e}^{-\la f_{[0,t)}, g_{[0,t)} \ra } \lla \ve (f_{[0,t)}), k_t (a)
\ve (g_{[0,t)}) \rra = \left( P^{f(t_0),g(t_0)}_{t_{1}-t_0} \circ \cdots \circ
 P^{f(t_{n-1}),g(t_{n-1})}_{t_{n}-t_{n-1}} \right) (a)
\end{equation}
where $0 = t_0 \leq t_1 \leq \cdots \leq t_n =t$ contains the discontinuities
of $f_{[0,t)}$ and $g_{[0,t)}$, and
\begin{equation} \label{standard Cassoc semigroups}
P^{c,d}_t := \mathrm{e}^{-t\la c,d\ra} \lla \ve (c_{[0,t)}), k_t (\,\cdot\,)
\ve (d_{[0,t)} ) \rra \;\; c,d \in \kil.
\end{equation}
\end{lemma}

\begin{deft}
A process $k \in \Proc_{\wb} (\Vil ; \Vil, \Exps)$ is a weak standard QS cocycle
if it satisfies formulas \eqref{stan Csemig decomp} and \eqref{standinit}.
\end{deft}

\begin{rem}
 Lemma \ref{standCM} says that
each standard QS cocycle is a weak standard QS cocycle; again, it is easy to observe that
if a weak standard QS cocycle is completely bounded, it is a standard QS cocycle.
\end{rem}

\begin{cor} \label{stand CN}
Let $k \in \Proc_{\wb} (\Vil ; \Vil, \Exps)$ be a weak QS  cocycle, $c,d \in \kil$. Then
\eqref{standard Cassoc semigroups} defines a semigroup acting on $\Vil$:
$ (P^{c,d}_t )_{t \geq 0}$.
\end{cor}

\begin{deft}    \label{standCdefsemig}
Semigroups defined by $\eqref{standard Cassoc semigroups}$ (for all pairs $c,d \in \kil$) are called
associated semigroups of the cocycle $k$. The semigroup
$P^{0,0}$ is called the Markov semigroup of the cocycle $k$.
\end{deft}

It is clear that two weak standard QS cocycles with identical
associated  semigroups are equal.

\begin{deft}    \label{standcMarkreg}
A weak standard QS cocycle is Markov-regular if its Markov semigroup is norm-continuous.
\end{deft}

The crucial result here, connecting the two types of cocycles introduced above is the following consequence of
Corollary \ref{cbisom}.

\begin{fact}  \label{Rcocycl0}
Let $l\in \Proc_{\wb} (\coalg; \Exps)$, $k\in \Proc_{\wb} (\coalg;\coalg, \Exps)$ and assume that for all
$f,g \in \Step$, $t \geq 0$
\[ E^{\ve(f)} k_{t,\ve(g)} = R \left(E^{\ve(f)} l_{t,\ve(g)}\right).\]
Then $l$ is completely bounded if and only if $k$ is, $l$ is a weak QS convolution cocycle
if and only if $k$ is a weak standard QS cocycle, and in the latter case
$l$ is Markov-regular if and only if $k$ is.
\end{fact}

\begin{proof}
The first fact follows directly from the definition of map $R$; observe that if either $l$ or $k$ (so in fact both)
is  completely bounded, then, for each $t \geq 0$,  $k_t=R_{B(\Fock)} l_t$. The second fact follows from Corollary
\ref{cbisom}, which gives for all $f_1, f_2, g_1, g_2 \in \Step$, $t \geq 0$
\begin{align*} R \left( \left(E^{\ve(f_1)} l_{t,\ve(g_1)}\right) \star \right.&
\left. \left(E^{\ve(f_2)} l_{t,\ve(g_2)}\right) \right) \\
 &=  \left(E^{\ve(f_1)} k_{t,\ve(g_1)}\right) \circ \left(E^{\ve(f_2)} k_{t,\ve(g_2)}\right) .\end{align*}
\end{proof}

\section[QSDEs in operator space context]{QS differential equations in operator space context}
\label{QSDEOS}
\setcounter{equation}{0}

As we have already seen in Section \ref{coalgQSDE}, to obtain QS convolution cocycles as solutions
of QS differential equations it is necessary to consider equations with nontrivial initial conditions.
Therefore in this section we extend the results of \cite{lwex} to such situations. As usual in this chapter,
the main stress will be put on completely bounded coefficients (called also stochastic generators).

Let $\Vil, \Wil$ ($\Wil \subset B(\hil)$) be operator spaces and assume that $\phi \in CB(\Vil; M_{\khat} (\Vil))$, $\theta \in CB(\Vil; \Wil)$.
By a QS differential equation with the coefficient $\phi$ and the initial condition $\theta$ is understood the equation
\begin{equation} \label{genQSDE}
dk_t = \wh{k}_t \circ \phi \,d \Lambda_t, \;\;\; k_0 = \iota \circ \theta.
\end{equation}
A process $k\in \Proc(\Vil; \Wil, \Exps)$ is a \emph{weak solution} of the equation \eqref{genQSDE} if for all
$\xi, \eta \in \hil$, $f,g \in \Step$, $t \geq 0$, $ v \in \Vil$
\begin{equation} \label{wQSDE}
\lla \xi \ve (f),\big( k_t (v)-\theta (v)\ot I_{\Fock}\big) \eta \ve(g)\rra
= \int^t_0 \lla \xi \ve (f) ,k_s \left(E^{\fhat(s)} \phi(v) E_{\ghat(s)}\right) \eta\ve (g)\rra \, ds.
\end{equation}
The definition of a strong solution requires more care - we need to explain how to define the process
that is actually QS integrated and whose matrix elements are represented by the right hand side of \eqref{wQSDE}.
Assume that a process $k \in \Proc(\Vil; \Wil, \Exps)$ has completely bounded columns.
For each $t \geq 0$, $ f\in \Step$, define
\[K_{t, \ve(f)} = \tau \circ (k_{t,\ve(f)})^{(\khat)} \circ \phi : \Vil \to C_{\Fock} \left(M_{\khat}(\Wil)\right)\]
($\tau$ denotes here a map implemented by the tensor flip $\hil \ot \Fock \ot \khat \mapsto \hil \ot \khat \ot \Fock$ - in fact this is the same
flip which is explicitly present in coalgebraic QS differential equations \eqref{qgQSDE} and
\eqref{CqgQSDE}).
Let the process $K\in \Proc (\Vil;M_{\kilhat}(\Wil), \Exps)$ be defined
by
\begin{equation} \label{integratedproc} K_t(v)  (\xi \ot \zeta \ot \ve(f))= K_{t, \ve(f)}(v) (\xi \ot \zeta)\end{equation}
($t \geq 0, \xi \in \hil, \zeta \in \khat, f \in \Step$, $v \in \Vil$).

We say that a process $k \in \Proc(\Vil; \Wil, \Exps)$ is a \emph{strong solution} of the equation \eqref{genQSDE}
if it is a weak solution, it has completely bounded columns and the process $K$ introduced above is locally square
integrable. In this case the First Fundamental Formula implies that for each $v \in \Vil$
\[ k_t(v)  = \theta(v) \ot I_{\Fock} + \int^t_0 K_s(v) \, d\Lambda_s. \]

\subsection*{Uniqueness of the solution}

\begin{propn} \label{uniq}
The QS differential equation \eqref{genQSDE} has at most one  weakly regular weak solution.
\end{propn}

\begin{proof}
Let $k\in \Procwr(\Vil; \Wil, \Exps)$ be the difference of two
weak solutions of \eqref{genQSDE}, and let $v \in \Vil$, $\xi, \eta \in \hil$,
$f,g \in \Step_D$ and $T \geq 0$. Then
by iteration
\[
\lla \xi \ve (f), k_t (v) \eta \ve (g)\rra = \int_{\Delta_n[0,t]}
\lla \xi \ve (f), k_{s_1} \Big( \phi^{\fhat (s_1)}_{\ghat (s_1)} \circ
\cdots
\circ\phi^{\fhat (s_n)}_{\ghat (s_n)}(v)\Big) \eta \ve (g)\rra\, d\bm{s}
\]
for each $n\in\bn$ and $t\in [0,T]$, where for $\zeta_1, \zeta_2 \in \khat$, $w \in \Vil$
\[
\phi^{\zeta_1}_{\zeta_2}(w): = E^{\zeta_1} \phi(w)  E_{\zeta_2}. \]
Now  weak regularity of $k$ allows us to claim that the integrand is bounded by
\[
C_{f,g,T} \left( \max \{ \|\dhat\|^2  : d \in \Ran f \cup \Ran g \}\| \phi\|_{\cb}\right) ^n
\| v\| \|\xi \| \|\eta\|,\]
where $C_{f,g,T}>0$ is a certain constant. The result follows.
\end{proof}

It is clear from the proof above that for the uniqueness of the solution one does not need to assume any
continuity properties of the initial condition.

\subsection*{Existence of the solution}

The solution of the equation \eqref{genQSDE} is constructed by means of iterated QS integrals.
Recall the notations of Section \ref{QSIntegrals}
and the formal inclusion $B(\hil) \subset \Op(E)$
(for $E$ - dense subspace of $\hil$).
Let $\up^{\theta, \phi}$ be the linear
map $\Vil \rightarrow \Seq_{\hil, \kil}$ defined by
$\up^{\theta, \phi} (v)_n = \up^{\theta,\phi}_n (v)$ ($n \in \bn_0$) where
\[
\upsilon^{\theta, \phi}_n = \theta_n \circ \phi_n,
\]
the maps
$\phi_n : \Vil \rightarrow M_{\khat^{\ot n}} (\Vil)$
and
$\theta_n : M_{\khat^{\ot n}} (\Vil) \to M_{\khat^{\ot n}} (\Wil)$
are defined (recursively) by
$\phi_0 = \id_{\Vil}$, $\theta_0 = \theta$,
$\phi_1 = \phi$, $\theta_1 = \theta^{(\khat)}$, and, for $n\geq 2$,
$\phi_n = \phi^{(\khat^{\ot (n-1)})}\circ \phi_{n-1}$,
$\theta_n = \theta^{(\khat^{\ot n})}$. Let
$\wt{\up}^{\theta, \phi}:\Vil \rightarrow \Seq_{\hil, \kil}$ be given by ($n \in \bn_0$)
\[\wt{\upsilon}^{\theta, \phi}_n = \tau_n \circ \upsilon^{\theta, \phi}_n,\]
 where $\tau_n: B(\hil \ot\kilhat^{\ot n}) \to B(\hil \ot \kilhat^{\ot n})$ is the tensor flip
 reversing the order of the copies of $\kilhat$.

The reader will recognise here relevant continuous extensions of the maps
introduced in Section \ref{coalgQSDE}.

\begin{lemma} \label{Funiq}
For any
$\phi \in CB(\Vil ; M_{\khat}(\Vil))$ and $\theta \in CB(\Vil; \Wil)$, the maps
$\up^{\theta, \phi},$ $ \wt{\up}^{\theta, \phi}$ take values in
$\Geomdagger_{\hil, \kil}$.
\end{lemma}

\begin{proof}
The result follows from properties of the liftings of completely bounded maps to matrix spaces, implying the estimates
($n \in \bn_0$)
\[ \|\phi_n \| \leq \|\phi\|_{\cb}^n, \;\;\; \|\theta_n \| \leq \|\theta\|_{\cb}, \;\;
\|\up^{\theta, \phi}_n \| \leq  \|\phi\|_{\cb}^n \|\theta\|_{\cb}.\]
\end{proof}

\begin{tw} \label{Guniq}
Let $\phi \in CB(\Vil ; M_{\khat}(\Vil))$ and $\theta \in CB(\Vil; \Wil)$. Then the process
$k := \Lambda \circ \wt{\up}^{\theta, \phi}$
is a strong solution of the QS differential equation
\eqref{genQSDE}.
\end{tw}

\begin{proof}
Write $\up:= \up^{\theta, \phi}, \wt{\up}:= \wt{\up}^{\theta, \phi}$
and prove first the following equality:
\begin{equation}
\label{upn} \up_n (E^{\zeta} \phi(v) E_{\chi} ) = E^{\zeta} \up_{n+1}(v) E_{\chi},
\end{equation}
where $\zeta, \chi \in \khat$, $v \in \Vil$, $n \in \bn$ (note that
on the right hand side Dirac operators are applied to the last copy of
$\khat$ in the tensor product). It is enough to show that
\[    \phi_n (E^{\zeta} \phi(v) E_{\chi} ) = E^{\zeta} \phi_{n+1}(v) E_{\chi},\]
and for that it is enough to compare matrix elements
of both sides, as all the maps involved are (completely) bounded.
Equality of the matrix elements follows from the inductively proved formula:
\[ E^{\zeta_1} \ldots E^{\zeta_n} \phi_n(v) E_{\chi_n} \ldots E_{\chi_1}
= E^{\zeta_1} \left(\phi \left(\cdots \left( E^{\zeta_n} \phi(v) E_{\chi_n} \right) \cdots
\right)\right) E_{\chi_1}, \]
valid for all $n \in \bn$, $ v \in \Vil$, $\zeta_1, \ldots, \zeta_n, \chi_1, \ldots, \chi_n
\in \khat$.

Note that there is a version of the formula \eqref{upn} for the map $\wt{\up}$:
\begin{equation}
\label{wtupn} \wt{\up}_n (E^{\zeta} \phi(v) E_{\chi} ) = E^{\zeta} \wt{\up}_{n+1}(v) E_{\chi};
\end{equation}
the only difference being that Dirac operators
on the right hand side are now applied to the first copy of
$\khat$ in the tensor product. With \eqref{wtupn} in hand we are ready to prove that
 $k$ satisfies \eqref{genQSDE}
weakly. Choose $f,g \in \Step$, $v \in \Vil$, $\xi, \eta \in \hil$.
In Guichardet notation,
\eqref{Guichardet FF1}  implies
\begin{align*}
&\int^t_0 ds \lla \xi \ve (f), k_s \left(E^{\fhat(s)} \phi(v) E_{\ghat(s)}\right) \eta \ve (g) \rra \\
=&\int^t_0 \, ds \int_{\Gamma_{[0,s]}} d \tau
\lla \xi \pi_{\fhat} (\tau), \big( \wt{\up}_{\# \tau}
(E^{\fhat(s)} \phi(v) E_{\ghat(s)})  \big)  \eta \pi_{\ghat} (\tau)  \rra \lla \ve (f), \ve (g)\rra \\
=&\int^t_0 \, ds \int_{\Gamma_{[0,s]}} d \tau
\lla \xi \pi_{\fhat} (\tau \cup s), \big( \wt{\up}_{\# (\tau \cup s)}
(v)  \big) \eta \pi_{\ghat} (\tau \cup s) \rra \lla \ve (f), \ve (g)\rra \\
=& \int_{\Gamma_{[0,t]}} d\sigma \big(1-\delta_\emptyset (\sigma) \big)
\lla \xi \pi_{\fhat} (\sigma) ,\wt{\up}_{\# \sigma} (v) \eta \pi_{\ghat} (\sigma) \rra
\lla \ve (f) , \ve (g)\rra  \\
= & \left\la \xi \ve (f), k_t (v)\eta\ve (g) \rra - \lla \xi, \theta (v) \eta \rra
\lla\ve (f) , \ve(g) \rra, \end{align*}
so $k$ satisfies the equation weakly.

Fix $T>0$ and $f \in \Step$.
The estimate \eqref{iterated FE} yields (for any
$0< t <T$, $ v \in \Vil$, $\xi \in \hil$)
\[
\|k_t (v) \xi \ve(f)\|^2 \leq \sum_{n=0}^{\infty} (C_{f,T})^n
\int_{\Delta_n[0,t]}
\|\phi \|_{\cb}^{2n} \|\theta\|_{\cb}^2 M_f^{2n} \|v\|^2 \|\xi\|^2 \, d\bm{s},
\]
where $M_f:= \max\{\|\dhat\|: d \in \Ran f\}$. This implies (recall the notation
introduced after \eqref{sreg})
\begin{equation} \label{sregsol}
\|k_{t,\ve(f)} \| \leq C_{f,T, \|\phi\|_{\cb}} \|\theta\|_{\cb},
\end{equation}
for some constant $C_{f,T, \|\phi\|_{\cb}}>0$. Thus
$ k\in \Proc_{\sr} (\Vil;\Wil, \Exps) \subset  \Procwr (\Vil;\Wil, \Exps) $.
This fact will be used below in conjunction with certain functorial property
of the processes constructed with the help of the map $\wt{\up}^{\theta, \phi}$
to prove the analog of \eqref{sregsol} for the completely bounded norm.

To this end, choose $n \in \bn$ and consider the maps
\[ \theta^{(n)}: M_n(\Vil) \to M_n(\Wil)\]
and
\[\wt{\phi^{(n)}} := \tau \circ \phi^{(n)}:M_n(\Vil) \to M_{\khat}(M_n (\Vil)),\]
where now $\tau$ denotes a tensor flip $B(\hil \ot \khat \ot \bc^n) \to
B(\hil \ot  \bc^n \ot \khat)$.
Denote $\wt{k} = \Lambda \circ \wt{\up}^{\theta^{(n)}, \wt{\phi^{(n)}}}$. Then
by the first part of the proof $\wt{k} \in \Proc_{\sr}(M_n(\Vil); M_n(\Wil), \Exps)$
is a weak solution of \eqref{genQSDE} with $\theta$ and $\phi$ replaced respectively by
$\theta^{(n)}$ and $\wt{\phi^{(n)}}$. Moreover by \eqref{sregsol} for $0<t < T$
\begin{equation} \|\wt{k}_{t,\ve(f)} \| \leq C_{f,T, \|\wt{\phi^{(n)}}\|_{\cb}} \|\theta^{(n)}\|_{\cb} =
 C_{f,T, \|\phi\|_{\cb} }\|\theta\|_{\cb}. \label{cbsregsol}\end{equation}
Further consider the equality
\begin{equation} \label{nlift}
\tau' \circ k_{t, \ve(f)}^{(n)} =  \wt{k}_{t, \ve(f)},\end{equation}
where this time $\tau'$ denotes the flip yielding the canonical complete isometry
$ M_n \left(C_{\Fock} (\Vil)\right) \cong C_{\Fock} \left(M_n (\Vil)\right) $.
To prove \eqref{nlift}, define a process $k' \in \Proc_{\sr}(M_n(\Vil); M_n(\Wil), \Exps) $
by the formula
\[ k' (\wt{v}) \wt{\xi} \ve(f) = \left( \tau' \circ k_{t, \ve(f)}^{(n)} (\wt{v}) \right) \wt{\xi},\]
for all $\wt{v} \in M_n (\Vil),$ $ f \in \Step,$ $ t \geq 0,$ $ \wt{\xi} \in \hil^{\oplus n}$.
As both processes $\wt{k}$ and $k'$ are weakly regular, to obtain \eqref{nlift}
it is enough to prove that $k'$ satisfies weakly  \eqref{genQSDE} with $\theta$ and $\phi$ replaced respectively by
$\theta^{(n)}$ and $\wt{\phi^{(n)}}$. Choose then any
$\wt{v}= [v_{i,j}]_{i,j=1}^n \in M_n (\Vil)$, $f,g  \in \Step$, $t \geq 0$,
$\wt{\xi} = \xi_1 \oplus \cdots \oplus \xi_n \in \hil^{\oplus n}$,
$\wt{\eta} = \eta_1 \oplus \cdots \oplus \eta_n \in \hil^{\oplus n}$ and compute:
\[
\lla \wt{\xi} \ve(f), \left( k_t'(\wt{v}) - \theta^{(n)} (\wt{v}) \right) \wt{\eta} \ve(g) \rra
=  \sum_{i,j=1}^n \lla \xi_i \ve(f), \left( k_t(v_{i,j}) - \theta(v_{i,j}) \right) \eta_j \ve(g) \rra \]
\[ =  \sum_{i,j=1}^n  \int_0^t
    \lla  \xi_i \ve(f), k_s \left(E^{\fhat(s)} \phi(v_{i,j}) E_{\ghat(s)}\right) \eta_j \ve(g) \rra ds \]
\[ = \sum_{i,j=1}^n  \int_0^t
    \lla  \xi_i \ve(f), k_{s,\ve(g)} \left(E^{\fhat(s)} \phi(v_{i,j}) E_{\ghat(s)}\right) \eta_j \rra ds \]
\[=  \int_0^t \lla
 \wt{\xi} \ve(f), \tau' \circ k_{s,\ve(g)}^{(n)} \left(E^{\fhat(s)} \wt{\phi^{(n)}} (\wt{v})  E_{\ghat(s)}\right)
 \wt{\eta} \rra ds \]
\[  = \int_0^t \lla
 \wt{\xi} \ve(f),  k'_s \left(E^{\fhat(s)} \wt{\phi^{(n)}} (\wt{v})  E_{\ghat(s)}\right)
 \wt{\eta} \ve(g) \rra ds.\]
The equality \eqref{nlift}, estimates \eqref{sregsol} and \eqref{cbsregsol} and the fact
that $\tau'$ is a (complete) isometry imply the estimate
\begin{equation} \label{cbs}
\|k_{t,\ve(f)} \|_{\cb} \leq C_{f,T, \|\phi\|_{\cb}} \|\theta\|_{\cb}.
\end{equation}
Thus $k$ has completely bounded columns and is $\cb$-strongly regular.
The estimate \eqref{cbs} shows that
the process $K$ introduced in \eqref{integratedproc} is also $\cb$-strongly regular (in fact even
$\cb$-H\"older continuous),
so in particular locally square integrable. This ends the proof.
\end{proof}

Thus the QS differential equation \eqref{genQSDE} has a unique  weak
solution; it
is also a strong solution, it is $\cb$-H\"older continuous and is given by  $\Lambda \circ \wt{\upsilon}^{\theta, \phi}$.
In the sequel we denote it by  $k^{\theta, \phi}$, simplified to $k^{\phi}$ if $\Vil = \Wil$,
$ \theta = \id_{\Vil}$.

\begin{rem}
The abundance of flips in the above computations is an unavoidable consequence of
putting the underlying operator space of a matrix space always `on the left'
(so for example $M_n(\Vil) = \Vil \ot \bc^n$) and the Fock space always `on the right'.
The advantages (certainly subjective) of such a convention were visible
for example when iterated QS integrals were defined.

The reader may have noticed that throughout the proof the essential work has been
done with columns of processes in question. To facilitate such operations it may be convenient
to view QS processes on operator spaces as families of linear maps (indexed
by time and exponential vectors) taking values in possibly unbounded column spaces.
This point of view is exploited in \cite{LSqsde}, where the existence and uniqueness
of solutions are proved when the stochastic generator has completely bounded
columns.
\end{rem}

\subsection*{Properties of the solution}

\begin{lemma} \label{Ctraffic}
Let $\Vil, \Wil$  be operator spaces,
$\phi \in CB\left(\Vil; M_{\khat}(\Vil)\right)$, $\theta \in CB(\Vil; \Wil)$. Then
$k^{\theta, \phi} \in \Proc^{\dagger}_{\cbsr} (\Vil ; \Wil, \Exps)$ and the following hold.
\begin{alist}
\item
If $\phi' \in CB\left(\Vil; M_{\khat}(\Vil)\right)$, $\theta' \in CB(\Vil; \Wil)$, then
$k^{\theta, \phi} = k^{\theta', \phi'}$ if and only if $ \theta = \theta'$  and 
$ (\theta^{(\khat)} \circ \phi)^{(\khat^{\ot(n-1)})} \circ \phi_{n-1}$
for all $n \in \bn$.
\item
If $\Vil$ and $\Wil$ are unital then $k^{\theta, \phi}$ is unital if and only if $\theta (1) = 1$,
$\phi (1) =0$.
\item
If $\Vil, \Wil$ are closed under adjoint operation then
$(k^{\theta, \phi} )^{\dagger} = k^{\theta^{\dagger}, \phi^{\dagger}}$.
In particular,
$k^{\theta, \phi}$ is real if and only if $\theta$ and $\theta^{(\wh{\kil})} \circ \phi$ are real.
\end{alist}
\end{lemma}
\begin{proof}
The first statement follows directly from the proof of Theorem \ref{Guniq} and the formula \eqref{Lambdadag}.
\newline
(a) Note an obvious relation
$k^{\theta, \phi} - k^{\theta', \phi} = k^{\theta-\theta', \phi}$. Further
$\wt{\up}^{\theta, \phi} = \wt{\up}^{\theta', \phi'}$ if and only if $\theta = \theta'$ and
\[ \up^{\theta, \phi}_n = (\theta^{(\khat)} \circ \phi)^{(\khat^{\ot(n-1)})} \circ \phi_{n-1}\]
 for all $n \in \bn$. The claim follows now from injectivity of the map $\Lambda$.
\newline
(b) follows again from injectivity of $\Lambda$.
\newline
(c) is an immediate consequence of \eqref{Lambdadag}.
\end{proof}

\begin{lemma} \label{compinlem}
Let $\Vil, \Wil_1, \Wil_2$  be operator spaces,
$\phi \in CB\left(\Vil; M_{\khat}(\Vil)\right)$, $\theta_1 \in CB(\Vil; \Wil_1)$,
 $\theta_2 \in CB(\Wil_1; \Wil_2)$. Then for all $t \geq 0$, $f \in \Step$
\begin{equation}  k^{\theta_2 \circ \theta_1, \phi}_{t, \ve(f)} = \theta_2 ^{(|\Fock\rangle)} \circ
k^{\theta_1, \phi}_{t, \ve(f)} \label{compin}\end{equation}
(recall that $\theta_2 ^{(|\Fock\rangle)}$ denotes the column lifting of the (completely bounded) map
$\theta_2$ to a map $C_{\Fock} (\Wil_1) \to C_{\Fock} (\Wil_2)$).   When
$k^{\theta_1, \phi}$ is completely bounded,
then so is $ k^{\theta_2 \circ \theta_1, \phi}$  and  for all $t \geq 0$
\[  k_t^{\theta_2 \circ \theta_1, \phi} = \theta_2^{(\Fock)} \circ k_t^{\theta_1, \phi}.\]
\end{lemma}
\begin{proof}
For the proof of the first statement it is enough to check that the process
$k' \in \Proc_{\sr}(\Vil; \Wil_2, \Exps)$ whose columns are given by the right hand side of \eqref{compin}
satisfies weakly the equation \eqref{genQSDE} with $\theta:=\theta_2 \circ \theta_1$,
and apply Theorem \ref{uniq}. Let $f,g \in \Step$, $ v \in \Vil$, $t \geq 0$.
\[ E^{\ve(f)} \left(k_t'(v) - \theta(v)\ot I_{\Fock}\right) E_{\ve(g)} =
   E^{\ve(f)} \left(k_{t, \ve(g)}'(v) - \theta(v)\ot E_{\ve(g)}\right) \]
\[= E^{\ve(f)} \left(\theta_2 ^{(|\Fock\rangle)} \circ
k^{\theta_1, \phi}_{t, \ve(g)}(v) - \theta(v)\ot E_{\ve(g)}\right)=
\theta_2 \left(E^{\ve(f)} \left(k_t^{\theta_1, \phi}(v) - \theta_1(v)\ot I_{\Fock} \right) E_{\ve(g)}\right)\]
\[ = \theta_2 \left(\int_{0}^t E^{\ve(f)} k_s^{\theta_1, \phi}\left(E^{\fhat(s)}\phi(v) E_{\ghat(s)}\right) E_{\ve(g)} ds \right)
  = \int_{0}^t E^{\ve(f)} k_s' \left(E^{\fhat(s)}\phi(v) E_{\ghat(s)}\right) E_{\ve(g)},\]
where the last equality is justified due to boundedness of $\theta_2$ and the fact that the integrands
above are piecewise constant. This implies that $k'$ is indeed a weak solution of the equation in question.

The second statement follows directly from the first.
\end{proof}

\begin{rem}
Observe that if $\Vil=\Wil$ and $\theta$ `commutes' with $\phi$, i.e.\
$ \phi \circ \theta =\theta^{(\khat)}\circ \phi$, then
$k^{\phi} \circ \theta = k^{\theta, \phi}$.
\end{rem}

Note the following fact:

\begin{fact}[Continuous dependence on initial conditions]
Let $\Vil,\Wil$ be operator spaces, $\phi \in CB\left(\Vil; M_{\khat}(\Vil)\right)$, $\theta \in CB(\Vil; \Wil)$.
Let $(\theta_n)_{n=1}^{\infty}$  be a sequence of maps in $CB(\Vil; \Wil)$ convergent (in $\cb$-norm) to $\theta$.
Then for each $f \in \Step$, $t \geq 0$ the sequence $(k_{t, \ve(f)}^{\theta_n, \phi})_{n=1}^{\infty}$
is convergent (in $\cb$-norm) to $k_{t, \ve(f)}^{\theta, \phi}$, and the convergence is uniform with respect to
$t$ on bounded subsets of $\br_+$.
\end{fact}

\begin{proof}
The claim follows from equality $k_{t, \ve(f)}^{\theta_n, \phi} - k_{t, \ve(f)}^{\theta, \phi} = k_{t, \ve(f)}^{\theta_n - \theta, \phi}$
and the estimate \eqref{cbs}.
\end{proof}

The result below is of the utmost importance for the next two sections. Its proof follows
easily from the uniqueness properties of the solutions of our QS differential equations.

\begin{tw}[\cite{lwjfa}]
\label{genstcoc}
Let $\Vil$ be an operator space, $\phi \in CB\left(\Vil;M_{\khat}(\Vil)\right)$. The
process $k^{\phi}$ is a weak standard QS cocycle on $\Vil$, whose associated semigroups
are norm continuous.
The generators of associated semigroups of $k^{\phi}$ are given by
$\{E^{\chat} \phi(\cdot) E_{\dhat} : c,d \in \kil\}$.
\end{tw}

The partial converse to this result may be found in \cite{LSqsde}, where it is proved,
using the methods of Section \ref{stochsection}, that each H\"older continuous standard
QS cocycle on a finite-dimensional operator space with a H\"older continuous adjoint
satisfies a QS differential equation of the form \eqref{genQSDE} with the initial condition given by the
identity mapping.

\section[Coalgebraic QSDEs on OS coalgebras]{Coalgebraic QS differential equations on OS coalgebras}
\label{coalgQSDEquat}
\setcounter{equation}{0}
In this section the general theory of QS differential equations with nontrivial initial conditions
is specified to the case of coalgebraic QS differential equations. Their solutions are shown to be
weak QS convolution cocycles.

Let $\coalg$ be an OS coalgebra, $\varphi \in CB\big(\coalg ; B(\khat)\big)$.  A
\emph{coalgebraic QS differential equation on $\coalg$}
(with the coefficient $\varphi$) is  the equation
\begin{equation} \label{CqgQSDE}
dk_t = k_t \star_{\tau} \varphi \,d\Lambda_t , \; \;
k_0 = \iota\circ\epsilon
\end{equation}
($\tau$ denoting a tensor flip exchanging the order of $\khat$ and $\Fock$,
$\iota$ indicating an ampliation).
A process  $k \in \Proc_{\wb} (\coalg ; \Exps)$ is a \emph{weak solution} of the
equation \eqref{qgQSDE}
if for all $f,g \in \Step$, $t \geq 0, a \in \coalg$
\begin{equation} \label{wCqgQSDE}
\lla \ve (f),\big( k_t (a)-\Cou (a)I_{\Fock}\big) \ve(g)\rra
=
\int^t_0  \left((\omega_{\ve(f), \ve(g)} \circ k_s) \star
(\omega_{\fhat (s), \ghat (s)} \circ \varphi )\right) (a) ds.
\end{equation}
Note that weak boundedness of $k$ is sufficient for the above formula to make sense,
as bounded linear functionals are automatically completely bounded.

To define strong solutions, it is necessary to repeat from Section \ref{QSDEOS} the construction of the
process $K$ associated with $k$.
Assume that a process $k \in \Proc(\Coalg; \Exps)$ has completely bounded columns
and let
\begin{equation} \label{Cbigphi}  \phi =R_{B(\khat)} \varphi: \coalg \to \coalg \ot B(\khat).
\end{equation}
For each $t \geq 0$, $ f\in \Step$, define
\[K_{t, \ve(f)} = \tau \circ (k_{t,\ve(f)})^{(\khat)} \circ \phi : \coalg \to C_{\Fock}
\left(B(\khat)\right)
=B(\khat; \khat \ot \Fock),\]
where $\tau$ denotes the map implemented by the tensor flip $\Fock \ot \khat\mapsto \khat \ot \Fock$.
Let the process $K\in \Proc (\coalg;B({\khat}), \Exps)$ be defined
by
\begin{equation} \label{integrproc} K_t (\zeta \ot \ve(f))(a) = K_{t, \ve(f)}(a) (\zeta)\end{equation}
($t \geq 0, \zeta \in \khat, f \in \Step$, $a \in \coalg$).

A process $k \in \Proc(\coalg;  \Exps)$ is a \emph{strong solution} of the equation
\eqref{CqgQSDE} if it is a weak solution, it has completely bounded columns and
the process $K$ introduced above is locally square
integrable. In this case the First Fundamental Formula implies that for each $a \in \coalg$
\[ k_t(a)  = \Cou(a) I_{\Fock} + \int^t_0 K_s(a) \, d\Lambda_s. \]

From the discussion above it should be clear that in fact coalgebraic QS differential
equations are a special case of QS differential equations considered in Section \ref{QSDEOS}.
This is formalised in the next lemma.

\begin{lemma}  \label{solequiv}
A process $k \in \Proc_{\wb} (\coalg;\Exps)$ is a weak (respectively, strong) solution
of the equation \eqref{CqgQSDE} if and only if it is a weak (resp., strong) solution
of the equation \eqref{genQSDE} with the initial condition $\theta:= \Cou$ and
the coefficient $\phi$ defined by \eqref{Cbigphi}.
\end{lemma}

\begin{proof}
As $\phi \in CB \left(\coalg; \coalg \ot B(\khat)\right) \subset
CB \left(\coalg; M_{\khat}(\coalg)\right)$, for the case of weak solutions it is enough to check that in the context formulated in the
lemma both sides of \eqref{wQSDE} and \eqref{wCqgQSDE} coincide. This follows from the
equalities:
\begin{align*} \Big\langle \ve(f), k_s  & \left. \left(E^{\fhat(s)}   \phi(a) E_{\ghat(s)} \right) \ve(g) \rra =
  \lla \ve(f), k_s\left(E^{\fhat(s)} \left( (\id_{\coalg} \ot \varphi)\Com(a) \right) E_{\ghat(s)} \right) \ve(g) \rra\\
 =&  \lla \ve(f), k_s \left(\left(\id_{\coalg} \ot \omega_{\fhat(s), \ghat(s)} \circ \varphi \right)
\Com(a) \right) \ve(g) \rra    \\
=&   \left((\omega_{\ve(f), \ve(g)}\circ k_s) \ot
  (\omega_{\fhat(s), \ghat(s)} \circ \varphi) \right)\left(\Com(a) \right)\\
=&\left((\omega_{\ve(f), \ve(g)} \circ k_s)
 \star(\omega_{\fhat (s), \ghat (s)} \circ \varphi)\right) (a).
\end{align*}
The equivalence of  being a strong solution is now a direct consequence of the definitions.
  \end{proof}

\begin{tw}
Let $\coalg$ be an OS coalgebra, $\varphi\in CB\left(\coalg; B(\khat)\right)$. The coalgebraic
QS differential equation \eqref{CqgQSDE} has a unique weak solution, denoted
further by $l^{\varphi}$; it is also a strong
solution.
\end{tw}

\begin{proof}
As weak solutions of the equation \eqref{CqgQSDE} are automatically weakly regular by \eqref{wc=wr},
the claims follow from Lemma \ref{solequiv}, Theorem \ref{uniq} and Theorem \ref{Guniq}.
Note that in the notation introduced after Theorem \ref{Guniq} $l^{\varphi}= k^{\Cou, \phi}$,
where, as usual, $\phi =R_{B(\khat)} \varphi$.
\end{proof}

\subsection*{Properties of the solution}

\begin{lemma} \label{Cqgtraffic}
Let $\varphi\in CB\left(\coalg;B(\khat)\right)$. The process $l^{\varphi} \in \Proc^{\dagger}_{\cbsr} (\coalg;\Exps)$ and
the  map $\varphi \mapsto l^\varphi$ is injective.
\end{lemma}

\begin{proof}
Straightforward consequence of Lemma \ref{Ctraffic} and the equalities $l^{\varphi} = k^{\Cou, \phi}$,
$\varphi= E_{B(\khat)} \phi$.
\end{proof}

\begin{lemma} \label{Rcocycl}
Let $l=l^{\varphi}$, $k=k^{\phi}$, where $\varphi \in CB\left(\coalg;B(\khat)\right)$ and $\phi = R_{B(\khat)} \varphi$.
Then for all $t \geq 0$, $ f \in \Step$
\[k_{t, \ve(f)} = R_{|\Fock\ra} l_{t,\ve(f)}.\]
\end{lemma}

\begin{proof}
Let $f,g \in \Step$, $ t \geq 0$. Note the following:
\begin{align*}  E^{\ve(f)} \left( \right.& \left. R_{|\Fock\ra} l_{t, \ve(g)} - \id_{\coalg} E_{\ve(g)} \right) =
R_{B(\khat)} \left(E^{\ve(f)}  \left( l_{t, \ve(g)} - \Cou(\cdot)  E_{\ve(g)} \right) \right) \\
=& R_{B(\khat)} \left( \int_0^t \lla \ve(f), l_s \left( E^{\fhat(s)} \phi( \cdot) E_{\ghat(s)} \right) \ve(g) \rra ds \right)\\
=&  \int_0^t E^{\ve(f)} R_{|\Fock\ra} l_{s, \ve(g)} \left( E^{\fhat(s)} \phi( \cdot) E_{\ghat(s)} \right) ds,\end{align*}
where the last equality can be justified as in the proof of Lemma \ref{compinlem}. Repeating the arguments of that proof
and using Theorem \ref{uniq} yields the desired formula.
\end{proof}

The above lemma in conjunction with Fact \ref{Rcocycl0} yields the following.

\begin{cor} \label{remRcocycl}
If either $l$ or $k$ above is completely bounded, so is the other one. In this case $k_t = R_{B(\Fock)} l_t$
for each $t \geq 0$.
\end{cor}

This circle of ideas allows us to formulate the counterpart of Theorem \ref{genstcoc}
for QS convolution cocycles.

\begin{tw}
\label{genconvcoc}
Let $\coalg$ be an OS coalgebra, $\varphi \in CB\left(\coalg;B(\khat)\right)$. The
process $l^{\varphi}$ is a weak QS convolution cocycle on $\coalg$, whose associated
convolution semigroups are norm continuous.
\end{tw}
\begin{proof}
Follows immediately from Lemma \ref{Rcocycl}, Lemma \ref{Rcocycl0} and Theorem \ref{genstcoc}.
\end{proof}

The following lemma paves the way for a converse of Theorem \ref{genconvcoc},
to be formulated in the next section.

\begin{fact}
Let $\varphi\in CB\left(\coalg;B(\khat)\right)$. The generators
of the convolution semigroups $\lambda^{c,d}$ associated with $l^{\varphi}$ are equal to
$\omega_{\chat, \dhat} \circ \varphi$ ($c,d \in \kil$). \end{fact}

\begin{proof}
Let $\phi=R_{B(\khat)} \varphi$ and note that
(in the notation introduced before Lemma \ref{Funiq})
$\up^{\Cou, \phi}_n = \varphi^{\star n}$ for all $n \in \bn$.
Therefore for all $c,d \in \kil$, $ a \in \coalg$
\[ \lambda_t^{c,d} (a) = \int_{\Gamma_t} d\sigma
\big\la \pi_{\chat} (\sigma) ,
\varphi^{\star \# \sigma} (a) \pi_{\dhat} (\sigma) \big\ra =
\sum_{n \geq 0}
\frac{t^n}{n!} (\omega_{\chat , \dhat} \circ \varphi)^{\star n}(a),\]
which concludes the proof.
\end{proof}

The above fact may be also proved using the analogous result for the associated semigroups
of a weak standard QS cocycle $k^{\phi}$ (see Theorem \ref{genstcoc})
and applying Lemma \ref{Rcocycl}.

\section{Completely positive and contractive cocycles} \label{CPCcoc}
\setcounter{equation}{0}

In this section we analyse completely positive and contractive QS convolution
cocycles on $C^*$-hyperbialgebras. It is shown that Markov-regularity is a sufficient (and necessary)
condition for such cocycles to satisfy coalgebraic QS differential equations with a completely bounded
coefficient. The precise form of the stochastic generator is given for such a case.

We begin with a straightforward consequence of
Lemma \ref{Ctraffic}.

\begin{fact}  \label{realconv}
Let $\alg$ be a  $C^*$-hyperbialgebra,
$\varphi \in CB\left(\alg;B(\khat)\right)$. The QS convolution cocycle
$l^{\varphi}$ is real if and only if $\varphi$ is real, and unital if and only if
$\varphi(1)=0$.
\end{fact}

\subsection*{Sufficient condition for the QS convolution cocycles on $C^*$-hyperbialgebras
to be stochastically generated}

The main result concerning CPC standard QS cocycles on $C^*$-algebras, including a
relevant converse of
Theorem \ref{genstcoc} is the following theorem, originating in \cite{lp} (for the
form of the generator see also \cite{Slava}).

\begin{tw} [\cite{lwptrf}, \cite{lwex}] \label{CPstandard}
Let $\alg\subset B(\hil)$ be a unital (nondegenerate) $C^*$-algebra and $k\in \Proc(\alg; \alg, \Exps)$.
Then the following are equivalent:
\begin{rlist}
\item  $k$ is Markov-regular, completely positive and contractive standard QS cocycle on  $\alg$;
\item $k=k^{\phi}$, where $\phi \in CB \big( \alg ; M_{\khat}(\alg)\big)$
satisfies $\phi (1) \leq 0$
and may be decomposed as follows:
\begin{equation}
 \phi(a) =  \Psi (a) -  a \ot \QSproj - E_{\zerohat} aJ - J^* a E^{\zerohat}
 \label{CPC0} \end{equation}
($a\in \alg$) for some maps $\Psi \in CP \big( \alg ; M_{\kilhat}(\alg'')$ and
$J \in R_{\kilhat}(\alg'')$.
\end{rlist}
\end{tw}

The convolution counterpart of Theorem \ref{CPstandard} is:

\begin{tw} \label{genstruct}
Let $\alg$ be a $C^*$-hyperbialgebra and $l\in \Proc(\alg;\Exps)$. Then the following are equivalent:
\begin{rlist}
\item $l$ is Markov-regular, completely positive and contractive QS convolution cocycle;
\item
$l=l^{\varphi}$ where $\varphi \in CB \big( \alg ; B(\kilhat) \big)$
satisfies $\varphi (1) \leq 0$
and may be decomposed as follows\tu{:}
\begin{equation}    \label{CPC1}
   \varphi (a) =  \psi (a) -  \Cou (a)
   \left(\QSproj + | e_0 \ra \la \chi | + | \chi \ra \la e_0 | \right)
   \end{equation}
($a \in \alg$)
for some map $\psi \in CP \big( \alg ; B(\kilhat) \big)$ and vector
$\chi \in \kilhat$;
\item
there is a unital $C^*$-representation $(\rho,\Kil)$ of $\alg$, a
contraction $D\in B(\kil;\Kil)$ and a vector $\xi\in\Kil$,
such that $l=l^\varphi$ where
\begin{equation} \label{alternative CPC}
\varphi (a) =
\begin{bmatrix}\la\xi | \\ D^*\end{bmatrix}
\big(\rho(a) - \Cou (a) I_{\Kil}\big)
\begin{bmatrix}|\xi \ra & D\end{bmatrix}
+ \Cou (a)\varphi (1)
\end{equation}
$(a\in\alg)$, and $\varphi (1)$ is nonpositive with block matrix of the
form
\[
\begin{bmatrix}
* & * \\
* & D^*D - I_{\kil}
\end{bmatrix}.
\]
\end{rlist}
\end{tw}

\begin{proof}
Note first that none of the notions and formulas in the theorem depends on the actual
faithful representation of $\alg$ chosen.

(ii) $\Longrightarrow$ (i)  \\
 Define $\phi = R_{B(\kilhat)} \varphi$ and assume that $\alg$ is faithfully and nondegenerately represented on some
  Hilbert space $\hil$. Then
\[ \phi (a) = \Psi (a) - a \ot \left(\QSproj +  E_{\zerohat} T + T^* E^{\zerohat} \right),\]
where $\Psi=R_{B(\kilhat)} \psi: \alg \to \alg \ot B(\khat)$ is completely positive by
Proposition \ref{Requiv}. Moreover
$\phi(1) \leq 0$.
Theorem \ref{CPstandard} shows that $k^{\phi}$ is completely positive and contractive.
Corollary \ref{remRcocycl} and again Proposition \ref{Requiv} yield (i).

(i) $\Longrightarrow$ (ii)    \\
Put for all $t \geq 0$
$k_t=  R_{B(\Fock)} l_t$. Fact \ref{Rcocycl0} and Proposition \ref{Requiv} imply that
$k$ is a Markov-regular CPC standard QS cocycle. By Theorem \ref{CPstandard} $k$ equals
$k^{\phi}$ for some $\phi \in CB\left(\alg;M_{\khat}(\alg)\right)$. Define
\[\varphi = \Cou^{(\khat)} \circ \phi \in CB\left(\alg;B(\khat)\right),\]
and let $\{\gamma_{c,d}:c,d \in \kil\}$ be generators of the associated convolution semigroups
of the cocycle $l$. By Theorem \ref{genstcoc} the generators of the associated
semigroups of $k^{\phi}$ are given by $\{E^{\chat}\phi(\cdot) E_{\dhat}: c,d \in \kil\}$.
As $ l_t = E_{B(\Fock)} k_t$ ($t \geq 0$), for each $c,d \in \kil$
\[ \gamma_{c,d} = E_{B(\kilhat)} \left( E^{\chat}\phi(\cdot) E_{\dhat} \right)= \omega_{\chat, \dhat} \circ
\varphi.\]
Therefore $l=l^{\varphi}$, as their respective associated convolution semigroups coincide.

The condition $\phi(1) \leq 0$ is clearly equivalent to $\varphi(1) \leq 0$.
It remains to prove that $\varphi$ may be decomposed as in \eqref{CPC1}.
To this end, assume that $\alg$ is faithfully and nondegenerately represented on $\hil$ in such a way that $\Cou$ extends
to a normal state on $\alg''\subset B(\hil)$. This can be always achieved by considering
a direct sum of a faithful representation and the GNS representation with respect to $\Cou$
or by considering the bidual of $\alg$ (in both cases $\Cou$ is even a vector state). 
The continuous extension of $\Cou$ to $\alg''$
will be denoted by the same letter.
Let $\phi:\alg \to \alg \ot B(\kil)$ have the form (\ref{CPC0}),
with $\Psi: \alg \to M_{\khat}(\alg'')$ being completely positive and
$J \in R_{\khat} (\alg '')$.
Put
\[ T = \Cou^{(\langle\khat|)} \circ J \in \langle \khat|, \;\; \psi = \Cou^{(\khat)}
\circ \Psi:\alg \to
B(\khat).\]
One can show that $\psi$ is completely positive
(essentially using the same techniques as ones used to
prove complete boundedness for lifted  maps, see \cite{lwex}). Further
\[ \Cou^{(\khat)} \circ (E_{\zerohat} a J) = \Cou(a) E_{\zerohat} T,\]
which can be checked by comparison of the respective matrix elements ($\chi, \zeta \in \khat$):
\[ E^{\chi} \left(\Cou^{(\khat)} \circ (E_{\zerohat} a J)\right) E_{\zeta} =
 \Cou (E^{\chi} (E_{\zerohat} a J E_{\zeta})) = \Cou(a) \langle \chi, \zerohat \rangle \Cou (J E_{\zeta}) = \]
\[= \Cou(a) E^{\chi} E_{\zerohat}  T  E_{\zeta} = E^{\chi} \Cou(a) E_{\zerohat} T E_{\zeta}.\]
(note that the normal extension $\Cou$ to $\alg''$ is necessarily multiplicative).
 Summing up, we obtain:
\begin{equation} \label{CPC}  \varphi(a) = (\Cou \ot \ida) \circ \phi (a) =
\Cou^{(\khat)} \circ \phi(a) = \psi (a) - \Cou(a)
\left(\QSproj + E_{\zerohat} T + T^* E^{\zerohat} \right). \end{equation}

Note that the above form of $\varphi$ enforces a specific form of $\phi$:
\[ \phi (a) = \Psi' (a) - a \ot \left(\QSproj +  E_{\zerohat} T + T^* E^{\zerohat} \right),\]
where $\Psi'=R_{B(\khat)}\psi: \alg \to \alg \ot B(\khat)$ is completely positive
(possibly different than $\Psi$ with which we started)
and $T \in \langle \khat |$.
The above shape of $\phi$ corresponds to $J = I_{\hil} \ot T$ in (\ref{CPC0}).

(ii) $\Longrightarrow$ (iii)  \\
 Let
\begin{equation} \label{psi Stinespring}
\begin{bmatrix} \la\xi | \\ D^* \end{bmatrix}
\rho(\cdot)
\begin{bmatrix} |\xi\ra & D \end{bmatrix}
\end{equation}
be a minimal Stinespring decomposition of $\psi$. Thus $(\rho,\Kil)$ is a
unital $C^*$-representation of $\alg$, $\xi$ is a vector in $\Kil$, $D$
is an operator in $B(\kil;\Kil)$ (and
$ \rho(\alg)\big(\bc\xi + \Ran D\big)$  is dense in $\Kil$).
Identity \eqref{alternative CPC} follows, with
\[
\varphi (1) =
\begin{bmatrix}
\|\xi\|^2 - 2\re \alpha & \la D^*\xi - c| \\
| D^*\xi - c\ra & D^*D - I_{\kil}
\end{bmatrix},
\]
where $\binom{\alpha}{c}=\chi$, so (iii) holds.

(iii) $\Longrightarrow$ (ii)  \\
Writing
\[
\begin{bmatrix}
t & \la d | \\
|d\ra & D^*D - I_{\kil}
\end{bmatrix}
\]
for the block matrix form of $\varphi (1)$, $\varphi$ has the
form \eqref{CPC1} where $\psi$ is given by~\eqref{psi Stinespring} and
\[
\chi = \binom{\frac{1}{2}(\|\xi\|^2 - t)}{D^*\xi - d}
\]
so (ii) holds.
\end{proof}

\subsection*{Precise form of the generator}

The proof of the second part of the implication (ii)
$\Longrightarrow$ (i) in Theorem \ref{genstruct} in the last
section, although short and simple, has its disadvantages. First of
all it requires working in a very specific representation, secondly
it indirectly, via Theorem \ref{CPstandard} uses deep
Christensen-Evans theorem on quasi-innerness of derivations on
represented $C^*$-algebras. Moreover to prove the existence of
suitable dilations of CPC Markov-regular QS convolution cocycles to
$^*$-homomorphic cocycles on $C^*$-bialgebras yet more explicit
description of the relevant stochastic generators is needed.  In
this section we present an alternative, more elementary approach to
this problem, following the ideas used in characterising the
structure of generators of CPC standard QS cocycles.

Various refinements and generalisations of this characterisation
were published in \cite{lwptrf} and \cite{lwex}, but the crucial analysis has been
carried out in \cite{lp} (see also \cite{Slava}).
Adapting arguments there requires some care, and again the $R$-map introduced
in Section \ref{OS coalgebras}  is an indispensable tool.
Straightforward attempt of investigating consequences of complete positivity of QS
convolution cocycles leads
to highly nontrivial considerations of a proper counterpart of the conditional
CPositivity condition. To the knowledge of the author, the $R$-map is not very helpful in that.
However, it appears that nonnegative-definite kernels with values in $C^*$-algebras behave well
under the $R$-map (in a sense to be clearly visible from the proof of the next proposition).

Further in this section $\alg$ denotes a fixed $C^*$-hyperbialgebra.
For any $\tau\in B(\alg)$ define $\partial{\tau}:\alg \times \alg \to \alg$ by
\[ \partial \tau (a_1,a_2) = \tau (a_1^*)a_2 - a_1^* \tau(a_2) - \tau(a_1^*) a_2 +
a_1^* \tau(1) a_2, \;\; a_1,a_2 \in \alg.\]
By analogy,  for any $f \in \alg^*$ define $\partial_{\Cou}f: \alg \times \alg \to \bc$
by
\[ \partial_{\Cou}f (a_1, a_2) = f(a_1^* a_2) - \Cou(a_1^*) f(a_2) - f(a_1^*) \Cou(a_2) +
\Cou(a_1^*) f(1) \Cou(a_2), \;\; a_1,a_2 \in \alg.\]
For the notion of derivations and their basic properties
that will be of use further we refer to the Appendix A.

We need to start with the finite-dimensional situation.
The key fact is the following result, corresponding to Theorem 4.1 in
\cite{lp}.

\begin{lemma}  \label{findim}
Assume that $\kil$ is finite dimensional. Let $\varphi\in CB\big(\alg; B(\kilhat)\big)$
and suppose that the
(weak) QS convolution cocycle $l:=l^{\varphi}\in
\Proc(\alg;\Exps_{\kil})$ is CPC. Then
there exist a unital representation $(\rho, \Kil)$ of $\alg$, a
($\rho,\Cou$)-derivation $\delta:\alg\to B(\bc;\Kil)$,
an operator $D \in B(\kil; \Kil)$ and a vector $d\in\kil$ such that
\begin{equation} \label{macierz00}
\varphi(a)=\begin{bmatrix}\lambda(a)&
\Cou(a)\langle d|+\delta^\dagger(a)D\\ \Cou(a)|d\ra +D^*\delta(a)&
D^*\rho(a)D-\Cou(a)I_{\kil}\end{bmatrix},
\;\;\; a\in \alg,
\end{equation}
where the functional $\lambda$ is real,
\[ \partial_{\Cou}\lambda (a_1, a_2) = \delta(a_1)^*\delta(a_2),
\;\; a_1,a_2 \in \alg,\]
and the following minimality condition holds:
\begin{equation} \label{minimality}
\Kil= \textup{cl} \Lin \{\delta(a)1 + \rho(a) D c : a \in \alg, c \in \kil\}.
\end{equation}
If $(\Kil', \rho', \delta', D')$ is another quadruple satisfying the
above conditions \textup{(}except possibly the minimality
condition\textup{)}, then there exists a unique isometry $V: \Kil \to
\Kil'$ such that
\[ \delta'(a) = V \delta(a), \;\; \rho'(a) V = V \rho(a), \;\;D' = VD, \;\; a \in \alg.\]
\end{lemma}

\begin{proof}   The proof is a modification of the argument used in the
proof of Lemma 4.5 in \cite{lp}, where $\kil$ is taken to be $\bc^d$.
Write $\varphi$ in block matrix form:
\[
\left[\begin{matrix} \lambda
 & \widetilde{\eta} \\ \eta
 & \sigma - \Cou(\cdot)I_{\kil} \end{matrix} \right].
\]
By Fact \ref{realconv} the map $\varphi$ is real,
in particular $\widetilde{\eta} = \eta^{\dag}$.
By Propositions \ref{Requiv} and \ref{Rcocycl},
$k=R_{B(\Fock_{\kil})} l$ is
a CPC standard QS cocycle and  $\phi=R_{B(\kilhat)} {\varphi}$ is real.
The map $\phi$ has block matrix form
\[
\left[ \begin{matrix}{\tau}  & \alpha^{\dag} \\ \alpha
 & \nu - \iota \end{matrix} \right],
\]
where $\tau=R\lambda$, $\alpha=R_{B(\bc;\kil)}{\eta}$ and
$\nu=R_{B(\kil)}{\sigma}$. Now Lemma 4.4 in \cite{lp}
implies that the map $\Psi$ from $\alg\times\alg$ to $\alg\ot B(\kilhat)$,
there identified with $M_{d+1}(\alg)$, defined by
\[ \Psi(a_1, a_2) = \left[ \begin{matrix} \partial\tau (a_1,a_2)  & \alpha^{\dag} (a_1^*a_2) -
a_1^* \alpha^{\dag}(a_2)\\ \alpha (a_1^* a_2) - \alpha(a_1^*)a_2 & \nu(a_1^* a_2) \end{matrix}
\right], \;\; a_1,a_2 \in \alg, \]
is nonnegative-definite. Observe that if $\psi: \alg \times \alg \to
B(\kilhat)$ is defined by the formula
 \[ \psi(a_1, a_2) = \left[ \begin{matrix} \partial_{\Cou}\lambda (a_1,a_2)  & \eta^{\dag} (a_1^*a_2) -
\Cou(a_1^*) \eta^{\dag}(a_2)\\ \eta (a_1^* a_2) - \eta(a_1^*)\Cou(a_2) & \sigma(a_1^* a_2) \end{matrix}
\right], \;\; a_1,a_2 \in \alg, \]
then $ \psi = (\Cou \ot {\rm id}_{B(\kilhat)} ) \circ \Psi$.
This in turn implies that $\psi$
is a nonnegative-definite kernel. Indeed, for any
$n\in\bn,\, a_1, \ldots ,a_n\in\alg$ and $T_1, \ldots ,T_n\in B(\kilhat)$
\[ \sum_{i,j=1}^n T_i^* \psi(a_i, a_j) T_j =
\big(\Cou \ot {\rm id}_{B(\kilhat)}\big)
\Big(\sum_{i,j=1}^n (1_{\alg} \ot T_i^*) \Psi(a_i, a_j) (1_{\alg} \ot T_j) \Big) \geq 0, \]
as $(1_{\alg} \ot T_i)^* =  (1_{\alg} \ot T_i^*) \in \alg \ot B(\kilhat)$,
$\Cou$ is CP, and $\Psi$ is nonnegative-definite.

Now let $(\Kil, \chi)$ be the minimal Kolmogorov construction associated with $\psi$. This
means that $\chi$ is a map $\alg \to B(\kilhat; \Kil)$ satisfying
\begin{align*}
& \chi(a_1)^* \chi(a_2) = \psi(a_1, a_2), \;\;\;\; a_1, a_2 \in \alg, \\
& \Kil = {\textup{cl} \Lin}\{\chi(a) \zeta: a\in \alg, \zeta \in \kilhat\}.
\end{align*}
Properties of $\psi$ imply that $\chi$ is linear and bounded. Write
$\chi = [\delta \; \gamma]$, where $\delta\in B\big(\alg;B(\bc;\Kil)\big)$
and $\gamma \in B(\alg; B(\kil;\Kil))$. Then, for any $a,b\in \alg$,
\[ \delta(a)^* \delta(b) = \partial_{\Cou}\lambda (a,b) \text{ and }
\gamma (a)^*\delta(b) = \eta(a^*b) - \eta(a^*)\Cou(b). \]
Setting $a=b=1$ shows that $\delta(1)=0$.
Now for $u \in \alg$ unitary, define
\[
\delta_u (a) = \delta(ua) - \delta(u) \Cou(a), \
\gamma_u(a) = \gamma(ua) \text{ and } \chi_u = [\delta_u \; \gamma_u],
\quad \text{ for } a\in\alg.
\]
A straightforward computation yields
\[ \chi_u (a_1)^* \chi_u (a_2) =   \chi (a_1)^* \chi (a_2).\]
The uniqueness of the minimal Kolmogorov construction implies the existence of a unique
isometry $\rho(u): \Kil \to \Kil$
given by the formula
\[ \rho(u) (\delta(a)1 + \gamma(a)c) =
\delta(ua)1 - \delta(u) \Cou(a) + \gamma(ua)c,
\;\; a \in \alg, c \in \kil.\]
It follows, by standard arguments, that
\[ \rho(a) (\delta(b)1 + \gamma(b)c) =  \delta(ab)1 - \delta(a) \Cou(b) + \gamma(ab)c,
\;\; a,b \in \alg, c \in \kil,\]
defines a bounded operator $\rho(a)$ on $\Kil$. Moreover, it is easily checked that the resulting map
$\rho: \alg \to B(\Kil)$
is indeed a representation of $\alg$. It immediately follows that $\delta$ is a
($\rho$,$\Cou$)-derivation and also, by the  minimality and the identity $\delta (1)=0$, that
$\rho$ is unital.
Put $D= \gamma(1)\in B(\kil;\Kil)$. Then $\gamma(a) = \rho(a)D$, and
furthermore
$\sigma(a) = D^* \rho(a)D$ and
$\eta(a) = \Cou(a)\eta(1) + D^* \delta(a)1$. This yields
($\ref{macierz00}$) with $d=\eta (1)1$.

The second part of the lemma follows once more from the uniqueness of the
Kolmogorov construction.
\end{proof}

The step from finite-dimensional to arbitrary noise dimension space follows
in exactly the same way as for standard cocycles.

 \begin{lemma}  \label{infdim}
Assume that $\kil$ is an arbitrary Hilbert space. Let $\varphi\in CB(\alg; B(\khat))$ and suppose that the
(weak) QS convolution cocycle $l^{\varphi}\in
\Proc(\alg;\Exps_{\kil})$ is CPC. Then the conclusions of Lemma~\ref{findim} hold.
\end{lemma}

\begin{proof}
Observe first that one can obtain, as in Theorem \ref{genstruct},
\begin{equation} \label{cont} \varphi(1)\leq 0\end{equation}
(it can be also deduced directly from the contractivity of $l^{\varphi}$ via the It\^o Formula).
Let $\{\kil_\iota:\iota \in \mathcal{I} \}$ be an indexing of the set of all
finite-dimensional subspaces of $\kil$, which is partially ordered by inclusion.
As in \cite{lwptrf} we consider finite-dimensional cut-offs of  both $l^{\varphi}$ and
$\varphi$ itself. For each $\iota \in\mathcal{I}$ denote by $\varphi_\iota$
the map $\alg \to B(\widehat{\kil_\iota})$ given by the formula
\[ \varphi_\iota(a) = P_\iota \varphi(a)  P_\iota, \;\; a \in \alg,\]
where $P_\iota\in B(\kilhat)$ is the orthogonal projection onto $\kilhat_\iota$.
Setting $l^{(\iota)} = l^{\varphi_\iota}$, $\Fock^\iota = \Fock_{\kil_\iota}$,
$\Exps^{\iota}=\Exps_{\kil_{\iota}}$ and letting $\mathbb{E}_\iota$ denote the vacuum conditional expectation
map from $B(\Fock_{\kil})$ to $B(\Fock^\iota)$,
it is easy to see that $l^{(\iota)} \in \Proc(\alg;\Exps^{\iota})$ is a
CPC QS convolution cocycle and that it satisfies
\[
l^{(\iota)}_t(a) = \mathbb{E}_\iota [l^{\varphi}_t(a)], \;\;
a\in \alg, t\in\br_+.
\]

Lemma \ref{findim} yields quadruples
$(\Kil_\iota, \rho_\iota, \delta_\iota, D_\iota)$, unique up to isometric
isomorphism, such that for all $a \in \alg$
\begin{equation*}
\varphi_\iota(a) =
\begin{bmatrix}\lambda(a)&\Cou(a)\langle d_\iota|+\delta_\iota^\dagger(a)D_\iota\\
\Cou(a)|d_\iota\ra +D_\iota^*\delta_\iota(a)&D_\iota^*\rho_\iota(a)D_\iota-\Cou(a)I_\iota\end{bmatrix},
 \end{equation*}
where $I_\iota$ denotes the identity operator on $\kil_\iota$.

Exploiting the uniqueness one can construct an inductive limit $\Kil$
of the Hilbert spaces $\Kil_\iota$. Denote by $U_\iota$ the respective isometry
$\Kil_\iota\to\Kil$. Then there is a unital representation $\rho$
of $\alg$ on $\Kil$, a ($\rho$,$\Cou$)-derivation
$\delta: \alg \to B(\bc;\Kil)$ and, for each $c \in \kil$ a vector $c_D \in \Kil$ such that
\[\rho(a)U_\iota = {U_\iota}\rho_\iota(a),
\ \
\delta(a) = U_\iota\delta_\iota(a)
\text{ and }
{c_D}  = U_\iota D_\iota c,\]
for all $\iota \in \mathcal{I}$, $a \in \alg$ and $c \in \kil_\iota$.
The map $c \mapsto c_D$ is linear; it remains to show that it is bounded.
To this end
observe that, for any $\iota\in\mathcal{I}$ such that $c\in\kil_\iota$,
\begin{align*}
\Big\la \binom{0}{c}, \varphi(1) \binom{0}{c} \Big\ra &=
\Big\la \binom{0}{c}, \varphi_\iota(1) \binom{0}{c} \Big\ra \\ &=
\big\la c, (D_\iota{^*} D_\iota - \Cou(1) I_\iota )c \big\ra =
\|D_\iota c \|^2 - \|c\|^2 = \|c_D\|^2 - \| c\|^2,
\end{align*}
and inequality \eqref{cont} implies that $\|c_D\| \leq \|c\|$. The operator
$D\in B(\kil; \Kil)$ given by $Dc = c_D$ completes the tuple whose existence
we wished to establish. Minimality holds by the construction.
\end{proof}

Automatic innerness of ($\rho,\Cou$)-derivations
leads to the following theorem.

\begin{tw}  \label{spec}
Let $\varphi\in CB(\alg;B(\khat))$, for a $C^*$-hyperbialgebra $\alg$, and
suppose that the weak QS convolution cocycle
$l^{\varphi}\in \Proc(\alg;\Exps_{\kil})$ is completely positive and contractive. Then
there exists a tuple $(\Kil,\rho,D,\xi,d,e,t)$ constisting of a
unital representation $(\rho,\Kil)$ of $\alg$, a contraction $D\in
B(\kil;\Kil)$, vectors $\xi\in\Kil$ and $d,e \in \kil$, and a real number
$t$, such that
\begin{equation}  \label{specified}
\varphi(a)=\begin{bmatrix}\lambda(a)&\Cou(a)\langle d|+\delta^\dagger(a)D\\
\Cou(a)|d\ra +D^*\delta(a)&D^*\rho(a)D-\Cou(a)I_{\kil}\end{bmatrix},
\end{equation}
$t=\lambda(1)\leq 0$, $d=(I_{\kil}-D^*D)^{1/2}e$, $\|e\|^2\leq -t$, and,
for all $a \in \alg$,
\begin{equation} \label{deltalambda}
\delta(a) = \big(\rho(a) - \Cou(a)\big)|\xi\ra, \;\;\;
\lambda(a) =
\Cou(a) (t - \|\xi\|^2)+ \langle \xi, \rho(a) \xi \rangle.
\end{equation}
\end{tw}

\begin{proof}
Lemma \ref{infdim} gives the form \eqref{specified} for some
$\rho, \Kil, \delta$ and $D$.
As all ($\rho, \Cou$)-derivations are inner (Corollary \ref{innderiv} of the appendix),
there exists $\xi \in \Kil$ such that
\[ \delta(a) = \rho(a) |\xi\ra - \Cou(a) |\xi\ra.\]

It remains to note that
\begin{equation}  \varphi(1)=
\left[ \begin{matrix} t &  \langle d| \\
|d \rangle & D^*D - I_{\kil}\end{matrix}\right],\label{contrac}\end{equation}
and the condition $\varphi(1)\leq 0$ implies contractivity of $D$, negativity of $t$ and the existence
of a vector $e\in \kil$ satisfying all the conditions above
(see the characterisation of positive matrices given in Lemma 2.1 of \cite{dilate}).
\end{proof}

\begin{cor}
Let $\varphi\in CB(\alg; B(\khat))$,  and let $l=l^{\varphi}$ be completely positive
and unital. Then  there exist a Hilbert space $\Kil$, a unital representation $\rho: \alg \to B(\Kil)$,
$\zeta \in \Kil$ and isometry $D\in B(\kil;\Kil)$ such that
\[ \varphi(a) = \left[ \begin{matrix}{\lambda(a)}
 & \langle  D^* \delta(a) | \\ | D^* \delta(a)\rangle
 & D^*\rho(a)D - \Cou(a) I_{\kil} \end{matrix}\right],\]
where (for all $a \in \alg$)
\[\delta(a) = \rho(a) \zeta - \Cou(a) \zeta,  \]
\[ \lambda(a) =  \langle \zeta, \rho(a) \zeta \rangle -  \Cou(a) \|\zeta\|^2,\]
\[ 
\lambda(1) = 0.\]
This is exactly the form corresponding to the purely algebraic case considered in
\cite{UweSch}.
\end{cor}

\begin{rem}
The characterisation in Theorem \ref{spec} yields, as announced in the beginning of this section,
an alternative proof of the second part of the implication (i) $\Longrightarrow$ (ii) of Theorem \ref{genstruct}.
Indeed, for $\varphi:\alg \to B(\kilhat)$ having a form (\ref{specified}),
define $S: \khat\to \Kil$ by $S = [ \zeta \;\; D ]$
(identifying here $\Kil$ with $|\Kil\rangle$). Then
\[ \varphi(a) = S^* \rho(a) S + \left[ \begin{matrix}{\lambda_0(a)}
 & \langle \Cou(a) \xi - D^* \zeta | \\ | \Cou(a) \xi - D^* \zeta\rangle
 & - \Cou(a) I_{\kil} \end{matrix}\right] \;\;\; (a \in \alg), \]
where $\lambda_0 (a) = \lambda(a) -\langle \zeta, \rho(a) \zeta \rangle$.
Note that as $\partial_{\Cou}\lambda_0 (a_1, a_2) = 0$
for any $a_1, a_2 \in \alg$, $\lambda_0= \lambda_0(1) \Cou$ - one can check that $\lambda_0 -  \lambda_0(1) \Cou$ is an
($\Cou,\Cou$)-derivation and use Corollary \ref{disap}.

Observe that $\psi: \alg \to B(\kilhat)$ defined by
\[ \psi(a) = S^* \rho(a) S,\;\;\;   a \in \alg,\]
is evidently completely positive.
Putting $T = [ \frac{1}{2}\lambda_0(1) \; \;\; \xi - D^* \zeta ] \in \langle\khat|$ yields
the representation (\ref{CPC1}).
\end{rem}

\begin{rem}
The proofs of the theorems in this section can be modified so that to obtain another proof
of the whole Theorem \ref{genstruct}, avoiding any references to Theorem \ref{CPstandard}.
In particular, one may prove that if $l$ is a Markov-regular weak QS convolution cocycle
satisfying weakly a coalgebraic QS differential equation \eqref{CqgQSDE} for a
linear map $\varphi:\alg \to \Op(\kilhat)$, and $l$ is completely positive and contractive,
then $\varphi$ must be completely bounded (this would first require describing what is meant
by coalgebraic QS differential equations with a possibly unbounded coefficient).
\end{rem}

\section[$^*$-homomorphic cocycles on $C^*$-bialgebras]
{$^*$-homomorphic QS convolution cocycles on $C^*$-bialgebras}
\label{multCalg} \setcounter{equation}{0}

In this section we characterise the stochastic generators of
$^*$-homomorphic (and more generally, weakly multiplicative) QS
convolution cocycles on a $C^*$-bialgebra $\alg$ in terms of
structure maps on $\alg$. Two possible definitions of quantum L\'evy
processes on $C^*$-bialgebras are also proposed and a relevant
version of the Sch\"urmann Reconstruction Theorem established.

\subsection*{Weakly multiplicative QS convolution cocycles}

As in general the processes constructed with the help of stochastic integration do not
leave the exponential domain invariant, the multiplicativity of processes on $C^*$-algebras
in principle has to be understood weakly.

\begin{deft}
Let $\alg$, $\blg$ be $C^*$-algebras, $\blg \subset B(\hil)$. A process
$k\in \Proc(\alg;\blg,\Exps)$ is called weakly multiplicative if
\begin{equation} \label{wmltp} \lla  k_t(a)^* \xi \ve(f), k_t (b) \eta \ve(g) \rra  =
\lla  \xi \ve(f), k_t (ab) \eta \ve(g) \rra \end{equation}
for all $t \geq 0$, $a,b \in \alg$, $\xi, \eta \in \hil$, $f,g \in \Step$.
\end{deft}

Compare the above definition with Definition \ref{multipalg} and
note that when $B=\bc$, they are equivalent. If $k\in
\Proc(\alg;\blg,\Exps)$ is weakly multiplicative and real, it must
be bounded. Then it (or rather its continuous extension) is
$^*$-homomorphic, so also completely bounded.

The multiplicative properties of iterated stochastic integrals in the operator-space theoretic
context are reflected by the following variant of Theorem 3.4 of \cite{lwhom}.
The notation used in its formulation is modelled on the one used in Section \ref{multalg};
as in fact we will use only a special case of this theorem, for the precise interpretation
we refer to the original paper of J.M.\,Lindsay and S.J.\,Wills.

\begin{tw}     \label{0Corolla}
Let $\alg$, $\blg\subset B(\hil)$ be unital $C^*$-algebras,
$\phi \in CB\left(\alg;M_{\khat}(\alg) \right)$, $\theta \in CB(\alg;\blg)$. Then the following
are equivalent:
\begin{alist}
\item $k^{\theta, \phi}\in \Proc(\alg;\blg,\Exps)$ is weakly multiplicative;
\item for all $n\in \bn_0$, $a,b \in \alg$
\[ \up_n^{\theta, \phi} (ab) = \sum_{\lambda \cup \mu  = \{ 1,\ldots,n\} }
\up^{\theta, \phi}_{\# \lambda} (a) (\lambda;n) \QSproj [\lambda \cap \mu ; n]
\up^{\theta, \phi}_{\# \mu} (b)(\mu ; n)\]
\end{alist}
\end{tw}

In favourable circumstances, the condition (b) of the above theorem can be simplified.
The following is Corollary 4.2 ($\alpha_2$) of \cite{lwhom}.

\begin{tw}            \label{Corolla}
Let $\alg$ be a unital $C^*$-algebra, $\phi \in CB\left(\alg;\alg \ot B(\khat) \right)$.
The weak standard QS cocycle $k^{\phi} \in \Proc(\alg; \alg,\Exps)$ is weakly multiplicative
if and only if
\begin{equation} \label{struc0} \phi(ab) =
\phi(a) (b \ot 1_{B(\khat)}) + (a \ot 1_{B(\khat)}) \phi(b)
+ \phi(a) (1_{\alg} \ot \QSproj) \phi(b), \;\;\; a,b \in \alg.\end{equation}
\end{tw}

From now on let $\alg$ be a fixed $C^*$-bialgebra. The theorem above will allow us to
characterise the stochastic generators of $*$-homomomorphic QS convolution cocycles on $\alg$.
Start with the following fact:

\begin{fact} \label{wmltequiv}
Let $\varphi \in CB\left( \alg;B(\khat) \right)$, $\phi = R_{B(\khat)} \varphi$.
Then $l:=l^{\varphi}$ is weakly multiplicative if and only if $k:=k^{\phi}$ is weakly
multiplicative.
\end{fact}

\begin{proof}
As usual $\alg$ is assumed to be faithfully and nondegenerately represented on some Hilbert space $\hil$.

Assume first that $l$ is weakly multiplicative. By Lemma \ref{realconv}
$l^{\dag}=l^{\varphi^{\dag}}$, and weak multiplicativity of $l$ implies in particular that
\[ \lla \xi \ve(f), (\id_{\alg} \ot l_t) (xy) \eta \ve(g) \rra =
   \lla (\id_{\alg} \ot l^{\dag}_t) (x^*) \xi \ve(f), (\id_{\alg} \ot l_t) (y) \eta \ve(g) \rra,\]
for all $x,y \in \alg \odot \alg$, $t\geq 0$, $\xi, \eta \in \hil$, $f,g \in \Step$.
Note that the formula above may be equivalently written as:
\begin{equation} \label{formo} \lla \xi \ve(f), \left(\id_{\alg} \ot l_{t, \ve(g)}\right) (xy) \eta  \rra =
   \lla \left(\id_{\alg} \ot l^{\dag}_{t,\ve(f)}\right) (x^*) \xi,
   \left(\id_{\alg} \ot l_{t_\ve(g)}\right) (y) \eta  \rra.\end{equation}
As both sides of the above equation are clearly continuous in $x$ and $y$ (separately),
first fixing $x$ and varying $y$ and then reverting this procedure one can deduce
that in fact the formula is valid for all $x,y \in \alg \ot \alg$. Choose  then $a,b \in \alg$
and let $x = \Com(a)$, $y = \Com(b)$. As $\Com$ is multiplicative,
$xy = \Com(ab)$ and with the help of Lemma \ref{Rcocycl}
equation \eqref{formo} takes the form
\[ \lla \xi \ve(f),  k_{t, \ve(g)} (ab) \eta  \rra =
   \lla  k^{\dag}_{t,\ve(f)} (a^*) \xi, k_{t,\ve(g)} (b) \eta  \rra,\]
which is exactly the statement of weak multiplicativity of $k$.

Assume conversely that $k$ is weakly multiplicative.  By Lemma \ref{Ctraffic}
$k^{\dag}=k^{\phi^{\dag}}$. Choose $a,b \in \alg$, $ t\geq 0$, $f,g \in \Step$. Then
\begin{align*} \lla l_t (a)^* \ve(f), l_t (b) \ve(g) \rra & = \left( l^{\dag}_{t,\ve(f)} (a^*) \right)^*
   l_{t,\ve(g)} (b)   \\
 & =  \left( \left( \Cou \ot \id_{|\Fock\ra} \right) k_{t, \ve(f)} (a^*) \right)^*
    \left( \left( \Cou \ot \id_{|\Fock\ra} \right) k_{t, \ve(g)} (b) \right) \\
  & =  \Cou \left( \left( k^{\dag}_{t, \ve(f)} (a^*)\right)^*   k_{t, \ve(g)} (b) \right)
   \end{align*}
where the second equality follows from Lemma \ref{Rcocycl} and  the third uses the equality
\[ ( \Cou \ot \id_{\la\Fock|} ) (X^*) (\Cou \ot \id_{|\Fock\ra}) (Y) = \Cou(X^*Y),\]
valid for all $X,Y \in R_{\Fock}(\alg)$. Further weak multiplicativity of $k$ implies
that
\[ \Cou \left( \left( k^{\dag}_{t, \ve(f)} (a^*)\right)^*   k_{t, \ve(g)} (b) \right)
=  \Cou \left( E^{\ve(f)} k_{t, \ve(g)} (ab)\right) \]
\[= E^{\ve(f)} \left( (\Cou \ot \id_{\la\Fock|} ) k_{t, \ve(g)} (ab)\right) =
E^{\ve(f)} l_{t, \ve(g)} (ab),\]
where in the last equality Lemma \ref{Rcocycl} was used again. Comparison of the above formulas
yields weak multiplicativity of $l$.
\end{proof}

\begin{propn}   \label{convmlt}
Let $\varphi \in CB\left( \alg;B(\khat) \right)$.
Then $l:=l^{\varphi}$ is weakly multiplicative if and only if
\begin{equation} \label{struc1} \varphi(ab) = \varphi(a) \Cou(b) + \Cou(a) \varphi(b) + \varphi(a) \QSproj \varphi(b), \;\;\;
a,b \in \alg.\end{equation}
\end{propn}

\begin{proof}
In view of Theorem \ref{Corolla} and Fact \ref{wmltequiv} it is enough to show the equivalence
of the conditions \eqref{struc0} and \eqref{struc1}, where again $\phi = R_{B(\khat)} \varphi$.
If \eqref{struc0} holds it is enough to apply the homomorphism $\Cou \ot \id_{B(\khat)}$ to it
to obtain \eqref{struc1}. Conversely,
if \eqref{struc1} holds, then for all $a_1, a_2, b_1, b_2 \in \alg$
\begin{align*}  (\id_{\alg} \ot \varphi) & \left( (a_1 \ot a_2)(b_1 \ot b_2)\right) =
 (\id_{\alg} \ot \varphi) (a_1 \ot a_2) (\id_{\alg} \ot \Cou) (b_1 \ot b_2) \\
 & + (\id_{\alg} \ot \Cou) (a_1 \ot a_2) (\id_{\alg} \ot \varphi) (b_1 \ot b_2) \\
 &+ (\id_{\alg} \ot \varphi) (a_1 \ot a_2)  (1_{\alg} \ot \QSproj)
 (\id_{\alg} \ot \varphi) (b_1 \ot b_2).\end{align*}
By linearity for all $x,y \in \alg \odot \alg$
\begin{align*} (\id_{\alg} \ot \varphi) (xy) &
 = (\id_{\alg} \ot \varphi) (x) (\id_{\alg} \ot 1_{B(\khat)}\Cou) (y) \\
 &+ (\id_{\alg} \ot 1_{B(\khat)}\Cou) (x) (\id_{\alg} \ot \varphi) (y) \\
  &+ (\id_{\alg} \ot \varphi) (x)  (1_{\alg} \ot \QSproj)
 (\id_{\alg} \ot \varphi) (y),\end{align*}
and again by (separate) continuity the formula remains valid for all $x,y \in \alg \ot \alg$.
Inserting $x= \Com(a)$, $y =\Com(b)$ (for some $a,b \in \alg$) and then using multiplicativity
of the coproduct and the counit property yields \eqref{struc0}.
\end{proof}

\subsection*{Unital $^*$-homomorphic QS convolution cocycles}

Fact \ref{realconv} and Proposition \ref{convmlt} yield the
following characterisation of the stochastic generators of unital
$^*$-homomorphic QS convolution cocycles on $C^*$-bialgebras.

\begin{tw} \label{CHH}
Let $\alg$ be a $C^*$-bialgebra, $\varphi \in CB\left(\alg; B(\khat)\right)$ and let
$l=l^{\varphi} \in \QSC(\alg;\Exps)$.
Then the following are
equivalent\tu{:}
\begin{rlist}
\item
$l$ is unital and $^*$-homomorphic;
\item
$\varphi$ vanishes at $1_{\alg}$ and satisfies
\begin{equation} \label{structuremap}
\varphi (a^*b) = \varphi (a)^* \Cou (b) + \overline{\epsilon (a)}
\varphi (b) + \varphi (a)^* \QSproj \varphi (b) ;
\end{equation}
\item
$\varphi$ has block matrix form
\begin{equation} \label{Cvarphi as block}
\begin{bmatrix}
\gamma & \delta^\dagger \\
\delta & \rho - \iota\circ\Cou
\end{bmatrix}
\end{equation}
in which $\iota$ is the ampliation $z\mapsto zI_{\kil}$\tu{;}
\begin{align}
& \rho : \alg \rightarrow B(\khat) \text{ is a unital
$^*$-homomorphism};
\label{CSch T1} \\
& \delta : \alg \rightarrow |\kil\ra \text{ is a }
(\rho, \Cou)\text{-derivation: }  \notag \\
& \hspace{3cm}
\delta (ab) = \delta (a) \epsilon (b) + \rho (a) \delta (b) ;
\label{CSch T2} \\
& \gamma : \alg \rightarrow \bc \text{ is linear and satisfies }
\notag \\
&\hspace{2cm}
\gamma(a^*b) =
\overline{\gamma (a)} \Cou (b) + \overline{\Cou (a)} \gamma (b)
+ \delta (a)^*\delta (b).
\label{CSch T3}
\end{align}
\end{rlist}
\end{tw}

\begin{deft}
Let $\alg$ be a $C^*$-bialgebra. A linear map $\varphi:\alg\to B(\kilhat)$
is called a structure map on $\alg$ if it vanishes
at $1_{\alg}$ and satisfies the formula \eqref{structuremap} for each $a,b \in \alg$.
\end{deft}

\begin{propn}   \label{cbstruct}
Every structure map $\varphi$ on a $C^*$-bialgebra $\alg$ is inner, that is there exists a
unital representation $\rho:\alg \to B(\kil)$ and a vector $\xi \in \kil$ such that
\[ \varphi(a) = \begin{bmatrix}\la\xi, (\rho(a) - \Cou(a) ) \xi \ra &  \la (\rho(a^*) - \Cou(a^*) ) \xi |  \\
|(\rho(a) - \Cou(a) ) \xi \ra & \rho(a) - \Cou(a) 1_{\kil} \end{bmatrix}, \;\;\; a \in \alg.\]
In particular,  $\varphi$ is completely bounded.
\end{propn}

\begin{proof}
It is clear that $\varphi$ must have a matrix form \eqref{Cvarphi as block},
where $\rho$ is a unital representation, $\delta$ is a $(\rho, \Cou)$-derivation, and
$\gamma$ is a linear functional on $\alg$ satisfying  \eqref{CSch T3} for all $a,b \in \alg$.
The existence of $\xi \in \kil$ such that
\[ \delta(a) = \rho(a) \xi - \Cou(a) \xi, \;\;\; a \in \alg,\]
follows from Corollary \ref{innderiv}. Then it is easy to check that the functional
$\wt{\gamma}: \alg \to \bc$ given by
\[ \wt{\gamma}(a) = \la\xi, (\rho(a) - \Cou(a) ) \xi \ra, \;\;\; a \in \alg,\]
coincides with $\lambda$ on $K:=\Ker \Cou$. As $\alg = (\Lin{1_{\alg}}) \oplus K$ as a vector space,
in fact   $\gamma =\wt{\gamma}$. The last statement follows.
\end{proof}

\begin{deft}
Every triple $(\gamma, \delta, \rho)$ such that $\rho:\alg \to B(\kil)$ is
a unital representation, $\delta: \alg \to |\kil\ra $ is a $(\rho,\Cou)$-derivation
and $\gamma:\alg \to \bc$ is a functional satisfying
\eqref{CSch T3} (for all $a,b \in \alg$) is called
a $\kil$-Sch\"{u}rmann triple on $\alg$.
\end{deft}

The above proposition, Theorem \ref{genstruct} and Theorem \ref{CHH} yield

\begin{tw} \label{homequiv}
Let $\alg$ be a $C^*$-bialgebra. There is a one-to-one correspondence between the following objects:
\begin{alist}
\item Markov-regular $^*$-homomorphic QS convolution cocycles in $\QSC(\alg; \Exps)$;
\item structure maps on $\alg$;
\item $\kil$-Sch\"urmann triples on $\alg$.
\end{alist}
The correspondence between (a) and (b) is given by $l^{\varphi} \leftrightarrow \varphi$.
\end{tw}

\begin{rem}
It is clear that the results concerning the perturbation of QS
convolution cocycles on coalgebras, presented in Section
\ref{perturb}, remain valid also for QS convolution cocycles on OS
coalgebras, if only the stochastic generators of the covolution
cocycles in question are completely bounded and the stochastic
generators of the operator cocycles implementing the perturbation
are bounded. In particular, considering the perturbation by unitary
(Weyl) cocycles, one obtains again the action of the Euclidean group
of $\kil$ on $\kil$-Sch\"urmann triples associated with unital
$^*$-homomorphic QS convolution cocycles on a $C^*$-bialgebra
$\alg$.
\end{rem}

\subsection*{Quantum L\'evy processes on $C^*$-bialgebras and their reconstruction
from generators}

Defining quantum L\'evy processes on $C^*$-bialgebras requires certain modifications
of the original, purely algebraic, definition of L.\,Accardi, M.\,Sch\"urmann and W.\,von Waldenfels
(\cite{asw}, \cite{schu}).
The problem is how to build the convolution increments of the process given that, in general,
the multiplication $\alg \odot \alg \to \alg$ need not extend
continuously to $\alg \ot \alg$. (This is a commonly met difficulty in the
theory of topological quantum groups,
see \cite{kus}). Below we outline two ways of overcoming this obstacle.

The simplest idea is to define a quantum L\'evy process using only the concept of distributions.

\begin{deft}  \label{wLp}
A \emph{weak quantum L\'{e}vy process} on a $C^*$-bialgebra $\alg$
over a unital $^*$-algebra-with-state $(\Blg , \omega)$ is a family
$\big(j_{s,t}\! :\alg \to \Blg\big)_{0 \leq s \leq t}$ of unital
$^*$-homomorphisms such that the functionals  $\lambda_{s,t}:= \omega
\circ j_{s,t}$ (which are automatically bounded, so also completely bounded) satisfy the following conditions,
for $0\leq r \leq s \leq t$:

\begin{rlist}
\item       \label{wQLPi}
$\lambda_{r,t} = \lambda_{r,s} \star \lambda_{s,t}$;
\item  \label{wQLPii}    
$\lambda_{t,t} = \Cou $;
\item \label{wQLPiii} 
$\lambda_{s,t} = \lambda_{0,t-s}$;
\item \label{wQLPiv}  \[
\omega \left( \prod^n_{i=1} j_{s_i,t_i} (a_i) \right) = \prod^n_{i=1}
\lambda_{s_i,t_i} (a_i)
\]
whenever $n \in \bn$, $a_1, \ldots, a_n \in\alg$ and
the intervals $[s_1,t_1[,\ldots ,[s_n, t_n[$ are disjoint\tu{;}
\item  \label{wQLPv}  
$\lambda_{0,t}  \to \Cou$ pointwise as $t \to 0$.
\end{rlist}
A weak quantum L\'{e}vy process on a $C^*$-bialgebra $\alg$ is called
Markov-regular if
$\lambda_{0,t} \to\Cou$ in norm, as $t \to 0$.
\end{deft}

The family $\lambda:=\big(\lambda_{0,t}\big)_{t \geq 0}$ is a pointwise continuous convolution
semigroup of functionals on $\alg$, called the \emph{one-dimensional distribution} of the
process; if the process is Markov-regular then $\lambda$ has a convolution
generator which is also referred to as the \emph{generator} of the weak
quantum L\'{e}vy process.
Two weak quantum L\'{e}vy processes on $\alg$,
$j^1$ over $(\Blg^1 ,\omega^1)$ and
$j^2$ over $(\Blg^2 ,\omega^2)$,
are said to be \emph{equivalent} if they satisfy
\[
\omega^1 \left( \prod^n_{i=1} j^1_{s_i,t_i} (a_i) \right) =
\omega^2 \left( \prod^n_{i=1} j^2_{s_i,t_i} (a_i) \right)
\]
for all $n \in \bn$, $a_1, \ldots , a_n \in \alg$ and
disjoint intervals $[s_1,t_1[, \ldots , [s_n,t_n[$.  Clearly
two weak quantum L\'{e}vy processes are equivalent if and only if their one-dimensional
distributions coincide, and if they are Markov-regular then this is
equivalent to the equality of their generators.

\begin{rem}
Note that the above definition of a weak quantum L\'evy process, in
contrast to the definition of a quantum L\'evy process on an
algebraic $^*$-bialgebra, does not yield a recipe for expressing the
joint moments of the process increments corresponding to the
overlapping time intervals, such as
\[
\omega(j_{r,t}(x) j_{s,t} (y)) \text{ where } 0 \leq r,s < t.
\]
To achieve the latter, one would have to formulate the weak convolution increment property
(wQLP\ref{wQLPi}) in greater generality and assume certain commutation relations between
the increments corresponding to disjoint time intervals.
For other examples of investigations
of the notion of independence in noncommutative probability without imposing
any particular commutation relations we refer to \cite{hkk}.
\end{rem}

As in the algebraic case, the generator of a Markov-regular weak
quantum L\'evy process vanishes on $1_{\alg}$, is real and is
\emph{conditionally positive}, that is positive on the kernel of the
counit. Observe that if $l\in \Proc(\alg;\Exps)$ is a unital
$^*$-homomorphic QS convolution cocycle then, defining $\Blg :=
B(\Fock)$, $\omega:= \omega_{\ve(0)}$, and $j_{s,t} := \sigma_s
\circ l_{t-s}$ for all $0\leq s\leq t$, we obtain a weak quantum
L\'evy process on $\alg$, called a \emph{Fock space quantum L\'evy
process}, Markov-regular if $l$ is.

The following theorem may be proved exactly along the lines of the Sch\"urmann
Reconstruction Theorem (Theorem \ref{recon});
all the necessary continuity properties follow from Theorem \ref{cbstruct}.

\begin{tw} \label{Crecon}
Let $\gamma$ be a real,
conditionally
positive linear functional on a $C^*$-bialgebra $\alg$ vanishing at $1_{\alg}$.
Then there is a (Markov-regular)
Fock space quantum L\'{e}vy process with generator $\gamma$.
\end{tw}

\begin{cor}
Every Markov-regular weak quantum L\'{e}vy process is equivalent to a Fock
space quantum L\'{e}vy process.
\end{cor}

Another notion, in a sense intermediate between weak quantum L\'{e}vy processes and Fock space quantum L\'{e}vy processes,
can be formulated in terms of product systems --- a similar idea
is mentioned in a recent paper of M.\,Skeide (\cite{Skeiden}). Recall that a
\emph{product system of Hilbert spaces} is
a `measurable' family of Hilbert spaces $E=\{E_t: t \geq 0 \}$, together with unitaries $U_{s,t}:E_s \ot E_t \to E_{s+t}$
($s,t \geq 0$) satisfying the associativity relations:
\begin{equation} \label{ass}
U_{r+s, t} (U_{r,s} \ot I_{t}) = U_{r, s+t} (I_{r} \ot U_{s,t}),
\end{equation}
($r,s,t \in \br_+$) where $I_s$ denotes the identity operator on $E_s$.
A \emph{unit} for the product system $E$ is a `measurable' family $\{u(t): t \geq0\}$ of vectors
with $u(t) \in E_t$ and $u(s+t) = U_{s,t} \big(u(s) \ot u(t)\big)$ for all
$s,t \geq 0$ (the unit is \emph{normalised} if, for all $t \geq 0$,
$\|u(t)\|=1$). For the precise
definitions we refer
to \cite{Arv}.
The unitaries $U_{s,t}$ implement the isomorphisms
$\sigma_{s,t}: B(E_s \ot E_t) \to B(E_{s+t})$.

\begin{deft}
A \emph{product system quantum L\'evy process} on $\alg$
over a product-system-with-normalised-unit $(E,u)$ is a family
$\big(j_t: \alg \rightarrow B(E_t)\big)_{ t \geq 0}$ of unital
$^*$-homomorphisms satisfying the following conditions, for
$r,s \geq 0$:
\begin{rlist}
\item
$j_{s+t} =\sigma_{s,t} \circ \left( (j_s \ot j_t) \circ \Com \right)$,
\item  \label{EuQLPii}
$j_0 = \iota_0\circ\Cou$,
\item \label{EuQLPiii}
$\omega_{u(t)}\circ j_t \to\Cou$
pointwise as $t \to 0$.
\end{rlist}
where $\iota_0$ denotes the ampliation $\bc\to B(E_0)$.
\end{deft}

The Fock space $\Fock$ corresponds to a product system $E$ by setting $E_t = \Focktot$,
and using  the obvious unitaries
whose existence is due to the exponential property \eqref{factor}
(in fact the product system mentioned above is the product system of the
$E_0$-semigroup $\{\sigma_t:t \geq 0 \}$
introduced in \eqref{shift}, see \cite{Arv}).
A normalised unit $\Omega$ is
given by $\Omega (t)= \ve(0)\in\Fock_{[0,t[}$, $t\geq 0$.
It is therefore easy to see that every Fock space quantum
L\'{e}vy process is a  product system quantum L\'{e}vy process.

\begin{propn}
Each product system quantum L\'{e}vy process on $\alg$ naturally
determines a weak quantum L\'{e}vy process on $\alg$ with the same
one-dimensional distribution.
\end{propn}

\begin{proof}
Let $j$ be a quantum L\'evy process on $\alg$ over a
product-system-with-normalised-unit $(E,u)$. We use an inductive
limit construction. Define $\wt{\Blg} := \bigcup_{t \geq 0}
(B(E_t),t)$ and introduce on $\wt{\Blg}$ the relation: $(T,r) \equiv
(S,s)$ if there is $t \geq \max\{r,s\}$ such that $\sigma_{r,t-r} (T
\ot I_{t-r} )= \sigma_{s,t-s} (S \ot I_{t-s})$, in other words we
identify operators with common ampliations. The associativity
relations \eqref{ass} imply that $\equiv$ is an equivalence
relation. Define $\Blg = \wt{\Blg}/\!\equiv$ and introduce the
structure of a unital $^*$-algebra on $\Blg$, consistent with the
pointwise operations:
\[
(T,t) + (S,t) = (T+S,t), \;\;
(S,t) \cdot (T,t) = (ST, t), \;\;
(T,t)^* = (T^*,t)
\]
($t \geq 0, S,T \in B(E_t)$).
The map $\wt{\omega}:\wt{\Blg} \to \bc$ defined by $\wt{\omega} (T,t) = \omega_{u(t)}(T)$
induces a state $\omega$ on $\Blg$.
For $s,t\in\br_+$ define
\[
j_{s,t}:\alg\to\Blg \text{ by }
x\mapsto \left[\sigma_{s, t-s} (I_s \ot j_{t-s}(x) ) \right]_{\equiv}.
\]
It is easy to see that the family
$\big(j_{s,t}\big)_{0 \leq s \leq t}$ is a weak quantum
L\'evy process on $\alg$ over $(\Blg, \omega)$.
\end{proof}

The construction in the above proof, informed by the case of QS
convolution cocycles, is a special case of the familiar construction of
$C^*$-algebraic inductive limits. The completion of $\Alg$ with respect to
the norm induced from $\wt{\Blg}$ is a unital $C^*$-algebra that may be
called the \emph{ $C^*$-algebra of finite range operators on the product
system $E$}.

\begin{rem}
A version of the reconstruction theorem also holds for unital,
completely positive, QS convolution cocycles on
$C^*$-hyperbialgebras. It is easily seen that if
$l\in\Proc(\alg;\Exps)$ is a Markov-regular, unital, completely
positive QS convolution cocycle on a $C^*$-hyperbialgebra $\alg$,
then the generator of its Markov convolution semigroup is real,
vanishes at $1_{\alg}$ and is conditionally positive. The GNS-type
construction from the proof of Theorem \ref{Crecon} yields a
completely bounded map $\varphi:\alg \to B(\kilhat)$ for which the
cocycle $l^{\varphi}$ is unital and completely positive according to
Proposition \ref{realconv} and Theorem \ref{genstruct} (of course
there is no reason why it should be $^*$-homomorphic, if $\alg$ is
not a $C^*$-bialgebra). Clearly the Markov convolution semigroup of
$l^{\varphi}$ coincides with that of $l$.
\end{rem}

\section{Dilations} \label{Dilations}
\setcounter{equation}{0}

In this section $\alg$ is a fixed $C^*$-bialgebra, and we are
concerned with the possibility of dilating completely positive and
contractive QS convolution cocycles to $^*$-homomorphic ones,
possibly by extending the dimension of the noise space. We begin
with casting the characterisation of generators of $^*$-homomorphic
cocycles obtained in the previous section in the form given in
Theorem \ref{spec}.

\begin{propn} \label{homprecise}
Let $(\Kil, \rho, D, \xi, d, t)$ be a tuple as in Theorem \ref{spec}
and let $\varphi$ be the map in $CB\big(\alg;B(\khat)\big)$ given by
the formulas \eqref{specified} and \eqref{deltalambda}. Then the
(weak) QS convolution cocycle $l^{\varphi}\in
\Proc(\alg;\Exps_{\kil})$ is $^*$-homomorphic if and only if the
following conditions hold\textup{:} \\
{\rm (i)} $D$ is a partial isometry,\\
{\rm (ii)} $Dd = 0$,\\
{\rm (iii)} $DD^* \in \rho(\alg)'$,\\
{\rm (iv)}  $t = - \| d \|^2,$\\
{\rm(v)} $DD^* \delta = \delta$,\\
where $\delta$ is the $(\rho,\Cou)$-derivation
$a\mapsto\big(\rho(a)-\Cou(a)I_{\kil}\big)|\xi\ra$.
\end{propn}

\begin{proof}
Fact \ref{realconv} and Proposition \ref{convmlt} imply that $l$ is
$^*$-homomorphic if and only if $\varphi$ is real and
\begin{equation} \varphi(ab) = \varphi(a)\QSproj \varphi(b) + \Cou(a) \varphi(b) +
\Cou(b) \varphi(a),\;\; \; a,b \in \alg.\label{homold}\end{equation}
In the language of Theorem $\ref{spec}$, the structure
relations \eqref{homold} translate into the following identities:
\begin{align*}
&D^* \rho(a) DD^* \rho(b) D = D^* \rho(ab) D,  \\
&D^* \delta(ab) + \Cou(ab)|d\ra  =
D^* \rho(a) D \big( D^* \delta(b) + \Cou(b)|d\ra\big) + D^* \delta(a) \Cou(b) + \Cou(a)\Cou(b)|d\ra,  \\
&\lambda(a^*b) = \big\langle D^* \delta(a)1 + \Cou(a) d, D^* \delta(b)1 + \Cou(b) d \big\rangle +
\lambda(a^*)\Cou(b)
+ \Cou(a^*) \lambda(b),
\end{align*}
for all $a,b \in \alg$. As in Proposition 3.3 of \cite{dilate},
this in turn may be shown to be equivalent to the conditions (i)-(v).
\end{proof}

Additionally, $l$ is unital and $^*$-homomorphic if and only if
(iii), (v) are satisfied, $D$ is an isometry, $\xi =0$, and $t =0$.

\subsection*{Stochastic dilations of CPC QS convolution cocycles}

This subsection is patterned
on \cite{dilate}, with all necessary modifications. Whenever the proofs in the convolution context
are straightforward adaptations of ones for standard QS cocycles, only the reference and the general ideas
behind the reasoning are indicated.

\begin{deft}
Let $\kil_0$ be a closed subspace of a standard noise Hilbert space $\kil$.
A QS convolution cocycle $j \in \Proc(\alg;\Exps_{\kil})$ is said to be a stochastic
dilation of a QS convolution cocycle $l\in \Proc(\alg;\Exps_{\kil_0})$ if
\[ l_t = \be_0 \circ j_t, \;\;\; t \geq0,\]
where $\condexp_0$ denotes the $\kil_0$-vacuum
conditional expectation (see Section \ref{Focksection}).
\end{deft}

The following result follows in exactly the same way as its counterpart
for standard cocycles (\cite{dilate}, Lemma 1.2).

\begin{propn} \label{dilgen}
Let $\varphi\in CB(\alg;B(\khat))$ and $\psi\in CB(\alg; B(\khat_0))$, and
let  $j=l^{\varphi} \in \Proc(\alg;\Exps_{\kil})$ and $l=l^{\psi}\in
\Proc(\alg;\Exps_{\kil_0})$
be the respective QS convolution cocycles. Then $j$ is
a stochastic dilation of $l$ if and only if
$\psi (\cdot) = P_0 \varphi(\cdot) P_0$, where $P_0\in B(\khat)$ denotes
the orthogonal projection onto $\khat_0$.
\end{propn}

\begin{rem}
Observe that the above characterisation excludes the possibility of
obtaining the exchange free dilations,  --- it can be seen directly
from ($\ref{homold}$) that if a Markov-regular $^*$-homomorphic QS
convolution cocycle is generated by a map having the form
\[
\begin{bmatrix}*&*\\ *&0\end{bmatrix}
\]
then it is identically 0. This uses Corollary \ref{disap}.
As to the creation/annihilation free dilations they are possible only for those CPC QS convolution cocycles, whose generators have the form
\[\begin{bmatrix}0&0\\ 0&*\end{bmatrix}.\]
\end{rem}

\begin{tw}
Every Markov-regular completely positive and contractive QS convolution cocycle on
a $C^*$-bialgebra $\alg$
admits a Markov-regular $^*$-homomorphic stochastic dilation.
\end{tw}

\begin{proof}
Let $l\in \Proc(\alg;\Exps_{\kil_0})$ be a Markov-regular CPC QS convolution cocycle.
Then $l=l^{\varphi}$
for some  $\varphi\in CB\big(\alg; B(\khat_0)\big)$ and we can assume that
$\varphi$ has matrix form
(\ref{specified}) for a tuple $(\Kil, \rho, D, \xi, d, e, t)$ with the
properties described
in Theorem \ref{spec}. Let $\kil_1, \kil_2$ be Hilbert spaces,
suppose that $d_1 \in \kil_1$, $d_2 \in \kil_2$, $D_1 \in B(\kil_1; \Kil)$
(all as yet unspecified) and consider the
map $\psi:\alg \to B(\kilhat)$, where $\kil :=
\kil_0 \oplus \kil_1 \oplus \kil_2$, given by ($a \in \alg$)

\begin{equation}  \psi(a) = \left[ \begin{matrix}{\lambda(a)}
 & \Cou(a)\langle d|+\delta^\dagger(a)D
 & \Cou(a)\langle d_1|+\delta^\dagger(a)D_1
 & \Cou(a)\langle d_2| \\
\Cou(a)|d\ra +D^*\delta(a) &
D^*\rho(a)D - \Cou(a)I_0
 & D^* \rho(a) D_1 & 0 \\
\Cou(a)|d_1\ra +D_1^*\delta(a) &
D_1^* \rho(a) D &
  D_1^*\rho(a)D_1 - \Cou(a)I_1 & 0 \\
\Cou(a)|d_2\ra  &
0 & 0& - \Cou(a)I_2
 \end{matrix}\right],\label{dilated} \end{equation}
with $I_i$ denoting $I_{\kil_i}$, $i=0,1,2$. Now
observe that $\psi$ can also be written in the form
\begin{equation}  \psi(a) = \left[ \begin{matrix}
{\lambda(a)} & \Cou(a)\langle \wt{d}|+\delta^\dagger(a)\wt{D}\\
\Cou(a)|\wt{d}\ra +\wt{D}^*\delta(a)&
\wt{D}^*\rho(a)\wt{D} - \Cou(a)I_{\kil} \end{matrix}\right],\label{specified2} \end{equation}
where
\[
\wt{d}=  \left( \begin{matrix} d \\ d_1 \\ d_2 \end{matrix}
\right) \in \kil \text{ and }
\wt{D} = \left[ \begin{matrix} D & D_1 & 0 \end{matrix} \right]
\in B(\kil;\Kil).
\]
As $\psi$ is clearly completely bounded, it generates a weak QS
convolution cocycle $l^{\psi}\in \Proc(\alg;\Exps_{\kil})$. It
follows from Proposition \ref{dilgen} that $l^{\psi}$ is a
stochastic dilation of $l^{\varphi}$; it remains to show that we can
choose the parameters $\kil_1, \kil_2$, $d_1, d_2$ and $D_1$ so that
$l^{\psi}$ is $^*$-homomorphic.

To this end, it suffices to put $\kil_1 = \Kil$, $\kil_2 = \bc$,
\[
D_1 =  \big( I_1 - DD^*   \big)^{\frac{1}{2}}, \; \;
d_1 = D e, \;\; d_2 = \sqrt{-(t + \|e\|^2)}.
\]
The above definitions make sense as $\|e \|\leq -t$ and $D$ is a
contraction.
It remains then
to check properties (i)-(v) of Proposition \ref{homprecise}.
First note that
\[ \wt{D}\wt{D}^* = DD^* + I_1 - DD^* = I_1,\]
which implies that the conditions (i), (iii) and (v)  are satisfied (one can easily check
that $\wt{D}^* \wt{D}$ is a selfadjoint projection). Further
we obtain (ii):
\[
\wt{D} \wt{d} =
D (I_0-D^*D)^{\frac{1}{2}} e +
\left( I_1 - DD^*\right)^{\frac{1}{2}}D e = 0.
\]
Finally (iv) follows since
\[
\|\wt{d}\|^2 =
\|(I_0-D^*D)^{1/2} e\|^2 + \| D e\|^2 - \big(t+\|e\|^2\big) = -t.
\]
This completes the proof.
\end{proof}

If $l$ in the above theorem is unital and $\dim \Kil = \dim
\Ran(I_{\Kil} - DD^*)$, then it is possible to obtain a unital
$^*$-homomorphic dilation $j \in \Proc(\alg;\Exps_{\kil})$ of $l$
 (with the noise dimension space $\kil=\kil_0 \oplus \Kil$).

\subsection*{Stinespring Theorem for QS convolution cocycles}

As the previous section was a variation on the theme of
\cite{dilate}, this one addresses the convolution counterpart of the
problem considered in \cite{Stine} for standard QS cocycles. We
shall show (in Theorem \ref{Stin}) that each Markov-regular,
completely positive, contractive QS convolution cocycle has a
Stinespring-like decomposition in terms of a $^*$-homomorphic
cocycle perturbed by a contractive process.

First we need some remarks on QS differential equations of the type:
\begin{equation}
{\rm d}W_t = F_t \big(I_{\khat} \ot W_t\big)\, {\rm d}\Lambda_t, \;\;\; W_0 =
I_{\Fock},\label{eqgen}
\end{equation}
where $F\in \Proc(\kilhat;\Exps)$ is a bounded process.
We say that $W$ is a weak solution of the above equation if for all
$f,g \in \Step$ and $t \geq0$
\[ \la \ve(f), (W_t - I_{\Fock}) \ve(g) \ra  =
\int_0^t {\la \fhat(s) \ot \ve(f),
F_s  (I_{\kilhat} \ot W_s) (\ghat(s) \ot \ve(g) \ra ds}.\]

The solution of the above equation is given by the iteration procedure:
\begin{align*}
& X_t^0 = I_{\Fock}, \;
X_t^1 = \int_0^t F_s (I_{\kilhat} \ot X_s^0) d\Lambda_s,\;\cdots\;,
X_t^{n+1} = \int_0^t F_s (I_{\kilhat} \ot X_s^n)
d\Lambda_s,\;\cdots \\
& W_t \ve(f) := \sum_{n=0}^{\infty} X_t^n \ve (f).
\end{align*}
Sufficient conditions for the above heuristics to be justified are that
$F$ is strongly measurable and has locally uniform bounds;
this is also sufficient for the uniqueness of strongly regular strong
solutions of the equation (\cite{Stine}, Proposition 3.1).
These conditions are clearly satisfied when
\[F_s = (\id_{B(\kilhat)} \ot l_s)(T), \;\;\; s\geq 0,\]
where $l$ is a Markov-regular, CPC QS convolution cocycle and $T\in
B(\khat) \ot \alg$.

Now let $j$ be the $^*$-homomorphic QS convolution cocycle
$l^{\varphi}$
 ($\varphi\in CB (\alg; B(\khat))$) and let $T\in B(\khat)\ot \alg$.
Assume that $W\in \Proc(\Exps)$
is a bounded solution to the equation
\begin{equation} \label{geneq} {\rm d}W_t =
(\id_{B(\kilhat)} \ot
j_t) (T) \big(I_{\khat} \ot W_t\big) {\rm d} \Lambda_t,\;\;\; W_0 = I_{\Fock}. \end{equation}
We shall identify sufficient conditions for $W$ to be a contractive
process later.
The next question to be addressed is: when can we expect a process $k\in
\Proc(\alg;\Exps)$ defined by
\[k_t(a) = j_t(a) W_t,\;\;\; a\in \alg, \; t\geq 0,\]
to be a Markov-regular QS convolution cocycle?

The quantum It\^o formula yields
\begin{align*}
\big\langle \ve(f), k_t (a) \ve(g) \big\rangle
= &
\big\langle j_t(a^*) \ve(f), W_t \ve (g) \big\rangle \\
= &
\Cou(a) \la \ve(f),\ve(g) \ra + \\
& \int_0^t \textrm{ds}
\Big( \big\la \tilde{j}_s (I_{\kilhat} \ot a^*)(\fhat(s) \ot
\ve(f)), \tilde{j}_s(T)   \wt{W}_s (\ghat(s) \ot \ve(g)) \big\ra + \\
& \quad \quad
\big\la \tilde{j}_s (\phi(a^*)) (\fhat(s) \ot \ve(f)),  \wt{W}_s (\ghat(s)
\ot \ve(g)) \big\ra  + \\
& \quad \quad
\big\la \tilde{j}_s (\phi(a^*)) (\fhat(s) \ot \ve(f)),
(\QSproj\ot I_{\Fock} )\tilde{j}_s(T)\wt{W}_s (\ghat(s)\ot\ve(g))\big\ra\Big)
\end{align*}
($f,g \in \Step, t\geq0$), where $\phi = (\varphi \ot \id_{\alg}) \circ \Com$,
$\tilde{j}_s = (\idB \ot j_s)$ and $\widetilde{W}_s = I_{\kilhat} \ot W_s$.
Defining analogously $\tilde{k}_s = (\idB \ot k_s)$ we see that the above
equation may be written as
\begin{align*}
&\big\langle \ve(f), k_t (a)  \ve(g) \big\rangle =
\Cou(a) \langle \ve(f), \ve(g) \rangle + \\
&
\int_0^t \textrm{ds} \left( \Big\la \fhat(s) \ot \ve(f),
\tilde{k}_s \left((I_{\kilhat} \ot a)T + \phi(a) +
\phi(a)(\QSproj\ot\ida)T\right) (\ghat(s) \ot \ve(g)) \Big\ra \right).
\end{align*}

 The process  $k$ is equal to $l^{\psi}$ for some  $\psi \in CB (\alg;
B(\khat))$ if and only if $\wt{\psi} := (\psi \ot \id_{\alg})\circ\Com$ is given by
\begin{equation} \label{eq1}
a\mapsto
(I_{\kilhat} \ot a)T + \phi(a) + \phi(a) (\QSproj\ot \ida)T.
\end{equation}
Note that we need to work with the left version of the $R$-map
introduced in Section \ref{OS coalgebras} because of the tensor flip
in the definition of the coalgebraic QS differential equation (\ref{coalgQSDE}).
Let $\tau = (\idB \ot \Cou ) ( T) \in B(\khat)$.
Then $(\ref{eq1})$ implies that
\begin{equation} \label{sigm} \psi(a) = \Cou(a) \tau + \varphi(a) (1 +
\QSproj\tau),\end{equation}
and so
\begin{equation} \label{eq2} \wt{\psi}(a) = \tau \ot a +  \phi(a) +
\phi(a) (\QSproj\tau \ot\ida).  \end{equation}
Comparing $(\ref{eq1})$ with $(\ref{eq2})$ yields
\begin{equation} (I_{\kilhat} \ot a)T + \phi(a) (\QSproj\ot \ida)T =  \tau
\ot a + \phi(a) (\QSproj\tau \ot\ida). \label{compar}\end{equation}
If $T = \tau \ot \ida$ then this condition is automatically satisfied. If
$j$ is unital, then $T= \tau \ot \ida$ is also necessary for \eqref{compar}
to hold: put $a=\ida$ and use $\phi(\ida)=0$.

Observe that when $T = \tau \ot \ida$ the equation $(\ref{geneq})$ takes
the simple form
\begin{equation} \label{pert} {\rm d}W_t = (\tau \otimes U_t W_t ){\rm d}
\Lambda_t,\;\;\; W_0 = I_{\Fock}, \end{equation}
with $U_t = j_t(1)$.  In this case the condition on $\tau$ assuring
contractivity of $W$ is also particularly simple.

\begin{tw} \label{gener}
Let $j=l^{\varphi}$ where $\varphi\in CB (\alg; B(\khat))$ and
$\alg$ is a $C^*$-bialgebra. Suppose that $j$ is $^*$-homomorphic
and $\tau\in B(\khat)$ satisfies the condition
\begin{equation} \label{con} \tau + \tau^* +\tau^*\QSproj\tau \leq 0. \end{equation}
Then the equation (\ref{pert}), with $U_t:= j_t(1)$, has a
unique contractive strong solution $W\in \Proc(\Exps)$.
Moreover the process $W_t^* j_t(\cdot) W_t$ is equal to $l^{\theta}$,
where
\[ \theta(a) = \Cou(a) \left( \tau^*+\tau + \tau^*\QSproj \tau \right) +
(I_{\kilhat} + \tau^* \QSproj)\varphi(a) (I_{\kilhat} + \QSproj\tau),
\;\;\; a\in \alg.\]
\end{tw}

\begin{proof}
The discussion before the theorem shows that the equation (\ref{pert}) has a unique
strongly regular strong solution $W\in \Proc(\Exps)$. The It\^o formula
yields, for
$u = \sum_{i=1}^k {\lambda_i \ve(f_i)}$, $k\in \bn$,
$\lambda_1, \ldots, \lambda_k\in \bc$,
$f_1,\ldots, f_k\in \Step$,
\begin{align*}
\la W_t u, W_t u \ra - \la u, u \ra =
&
\sum_{i,j=1}^k \overline{\lambda_i} \lambda_j
\int_0^t \textrm{ds} \Big( \big\la \fhat_i(s) \ot \ve(f_i),
\tau\fhat_j(s) \ot U_s\ve(f_j) \big\ra \\
& \qquad \qquad +
\big\la \tau\fhat_i(s) \ot U_s\ve(f_i),  \fhat_j(s) \ot \ve(f_j) \big\ra  \\
& \qquad \qquad +
\big\la \tau\fhat_i(s) \ot U_s\ve(f_i),  \QSproj\tau\fhat_j(s) \ot
U_s\ve(f_j) \big\ra\Big).
\end{align*}
As $U_s=j_s(1)$ and $j$ is $^*$-homomorphic, each $U_s$ is a
projection. Therefore putting
\[x(s) = \sum_{i=1}^k {\lambda_i \fhat_i(s) \ot U_s \ve(f_i)}, \;\;\;s \in [0,t],\]
 yields
\[ \la W_t u, W_t u \ra - \la u, u \ra=  \int_0^t \textrm{ds} \left\la x(s), \left((\tau+\tau^*+
\tau^*\QSproj\tau) \ot \idf\right) x(s) \right\ra \leq 0.\]
It follows that $W$ is contractive.

The proof of the second part of the theorem is a combination of the
considerations before its formulation and one more application of the It\^o
formula. Again let $f,g\in \Step$, $t \geq 0$, $a \in \alg$ and
$T=\tau \ot 1_{\alg}$,
let $\tilde{j}$, $\tilde{k}$, $\wt{W}$ and $\psi$ be defined as in the
discussion before the theorem and set
$\wt{\psi} = (\psi\ot\id_{\alg})\circ\Com$. Then
\begin{align*}
\langle \ve(f), W_t^* j_t(a) W_t\ve (g) \rangle
= & \,
\langle W_t \ve(f), j_t(a) W_t \ve (g) \rangle \\
= &
 \, \Cou(a) \la \ve(f),\ve(g) \ra + \\
& 
\int_0^t \textrm{ds}
\Big( \big\la \widetilde{W}_s (\fhat(s) \ot \ve(f)),
\tilde{k}_s(\wt{\psi}(a))   (\ghat(s) \ot \ve(g)) \big\ra + \\
& 
\big\la \tilde{j}_s (T) \widetilde{W}_s  (\fhat(s) \ot \ve(f)),
\tilde{j}_s(I_{\kilhat} \ot a)
   \widetilde{W}_s (\ghat(s) \ot \ve(g)) \big\ra + \\
& \big\la \tilde{j}_s (T) \widetilde{W}_s (\fhat(s) \ot \ve(f)),
(\QSproj\ot \idf) \tilde{k}_s(\wt{\psi}(a))(\ghat(s) \ot \ve(g))
\big\ra \Big).
\end{align*}
Finally,  (\ref{sigm}) yields
\begin{align*}
\big\la \ve(f), W_t^* j_t(a) W_t\ve(g) \big\ra
=&
\Cou(a) \big\la \ve(f),\ve(g) \big\ra \\
& \qquad
+ \int_0^t \textrm{ds}
\lla \fhat(s) \ot \ve(f), \widetilde{W}^*_s
\tilde{j}_s(\wt{\theta}(a)) \widetilde{W}_s (\ghat(s) \ot \ve(g)) \rra.
\end{align*}
where $\wt{\theta} = \big(\theta\ot\id_{\alg}\big)\circ\Com$.
This completes the proof.
\end{proof}

For each $t\geq 0$ denote the orthogonal projection from $\Fock$ onto
$\Fock_{[t, \infty[}$ by $P_{\kil,[t,\infty[}$. The following result
may be proved by differentiation, as with its predecessor for standard QS
cocycles, Lemma 4.2 of \cite{Stine}.

\begin{fact} \label{Stingen}
Let $\kil$ be an orthogonal direct sum of Hilbert spaces:
$\kil_0\oplus\kil_1$, let
$ \varphi\in CB\big(\alg;B(\wh{\kil}_0)\big)$ and
$\psi\in CB\big(\alg; B(\wh{\kil})\big)$, and
let  $k^0=l^{\varphi} \in \Proc(\alg;\Exps_{\kil_0})$ and
$k=l^{\psi}\in \Proc(\alg;\Exps_{\kil})$
be the respective weak QS convolution cocycles. Then
\[k_t(a) = k^0_t(a) \ot P_{\kil_1,  [t,\infty[},
\;\;a\in \alg, t\geq 0,\]
if and only if
\[
\psi (a) = \left[\begin{matrix}\varphi(a) & 0 \\
 0 & - \Cou(a)I_1 \end{matrix} \right],
\;\;a\in \alg,
\]
where $I_1 = I_{\kil_1}$.
\end{fact}

We are ready for the main theorem of this section.

\begin{tw} \label{Stin}
Let $k \in  \Proc(\alg;\Exps_{\kil_0})$ be a Markov-regular CPC QS
convolution cocycle on a $C^*$-bialgebra $\alg$. Then there exists
another Hilbert space $\kil_1$, a Markov-regular, $^*$-homomorphic
QS convolution cocycle $j \in \Proc(\alg;\Exps_{\kil})$, where $\kil
:= \kil_0\oplus\kil_1$, and a contractive process $W\in
\Proc(\Exps_{\kil})$, such that
\[ \wt{k}_t(a) = W_t^*j_t(a)W_t, \;\; t\geq 0, a \in \alg,\]
where $\wt{k}_t(a) := k_t(a) \ot P_{\kil_1,  [t, \infty[} $.
A process $W$ may be chosen so that
it satisfies the QS differential equation
\begin{equation} \label{pert1} {\rm d}W_t = (l \otimes U_t W_t ){\rm d}
\Lambda_t,\;\;\; W_0 = I_{\Fock_\kil} \end{equation}
for some $l \in B(\hat{\kil})$ in which $U \in \Proc(\Exps_{\kil})$
is the  projection-valued process given by
$U_t = j_t(1)$, $t \geq 0$.
\end{tw}

\begin{proof}
Let $\varphi\in CB\big(\alg; B(\kilhat)\big)$ be the stochastic  generator of $k$
(so that $k=l^{\varphi}$) and let
$(\Kil, \rho, D, \xi, d, e, t)$ be an associated tuple, as in
Theorem \ref{spec}. Set $\kil_1=\Kil$ and define
$\theta: \alg \to B(\widehat{\kil})$ by
\[ \theta(a) = \left[ \begin{matrix}
\lambda(a) - t \Cou(a) & 0 & \delta^\dagger(a)\\
 0 & - \Cou(a) I_0 & 0 \\
\delta(a)  & 0 & \rho(a) - \Cou(a) I_1
\end{matrix}
 \right], \;\;\; a\in \alg, \]
where $I_i$ denotes $I_{\kil_i}$, $i=0,1$ and $\delta$ is the
$(\rho,\Cou)$-derivation:
$a\mapsto\big(\rho(a)-\Cou(a)I_{\kil}\big)|\xi\ra$. The map $\theta$
is completely bounded and as such generates a Markov-regular weak QS
convolution cocycle $j=l^{\theta} \in \Proc(\alg;\Exps_{\kil})$. It
is easily checked that $\theta$ satisfies the structure relations of
Theorem \ref{CHH}, so $j$ is $^*$-homomorphic. Now choose any
contraction $B\in B(\kil_1; \kil_0)$  and define
\[ \tau =   \left[ \begin{matrix} \frac{1}{2}t & \la \xi| & 0 \\
 0 & -I_0 & B  \\ 0 & D & -I_1     \end{matrix}\right]
\in B(\kilhat). \]
Then
\[ \tau^* + \tau + \tau^*\QSproj \tau = \left[ \begin{matrix} t & \la\xi| & 0 \\
 |\xi\ra & D^*D-I_0 & 0  \\ 0 & 0 & B^*B-I_1
\end{matrix}\right]\leq 0 ,\]
 as $B$ is  a contraction, and $\varphi(1)\leq 0$ (see (\ref{contrac})).

Theorem \ref{gener} yields the existence of  a contractive process $W \in \Proc(\Exps_{\kil})$
satisfying the QS differential equation (\ref{pert1})
and  shows that the process $l \in \Proc(\alg;\Exps_{\kil})$ given by
\[ l_t(a) = W_t^*j_t(a)W_t, \;\; t\geq 0, a \in \alg,\]
is equal to $l^{\psi}$ where  $\psi:\alg \to B(\hat{\kil})$ is defined by
\begin{eqnarray*}
\psi(a) &=& \Cou(a) \left( \tau^*+\tau + \tau^*\QSproj \tau \right) + (1 + \tau^* \QSproj ) \theta(a) (1 + \QSproj \tau) \\
 &=& \Cou(a) \left[ \begin{matrix} t & \la\xi| & 0 \\
 |\xi\ra & D^*D-I_0 & 0  \\ 0 & 0 & B^*B-I_1     \end{matrix}\right] \\
  &&+ \left[ \begin{matrix} 1 & 0 & 0 \\
 0 & 0& D^* \\ 0 & B^* & 0     \end{matrix}\right] \cdot  \left[ \begin{matrix}   \lambda(a) - t \Cou(a) & 0 &
\delta^\dagger(a)\\
 0 & - \Cou(a) I_0 & 0 \\ \delta(a) & 0 & \rho(a) - \Cou(a) I_1
\end{matrix}\right]  \cdot  \left[ \begin{matrix} 1 & 0 & 0 \\
 0 &  0 & B \\ 0 & D & 0   \end{matrix}\right] \\
 &=&  \Cou(a) \left[ \begin{matrix} t & \la\xi| & 0 \\
 |\xi\ra & D^*D-I_0 & 0  \\ 0 & 0 & B^*B-I_1     \end{matrix}\right] \\
 && + \left[ \begin{matrix}   \lambda(a) - t \Cou(a) & \delta^\dagger(a)D & 0 \\
 D^* \delta(a) & D^* \rho(a) D - \Cou(a) D^*D & 0 \\ 0 & 0 &  - \Cou(a) B^*B \end{matrix} \right] \\
&=& \left[ \begin{matrix}   \lambda(a) &  \delta^\dagger(a)D + \Cou(a)\la \xi | & 0 \\
 D^* \delta(a)+ \Cou(a)|\xi  \ra & D^* \rho(a) D - \Cou(a)I_0 & 0 \\ 0 & 0
&  - \Cou(a) I_1 \end{matrix}\right] =
\left[ \begin{matrix}   \varphi(a) & 0 \\  0 &  - \Cou(a) I_{\kil}
\end{matrix}\right]. \end{eqnarray*}
Application of Fact \ref{Stingen}  completes the proof.
\end{proof}

\section{Examples} \label{Exam}
\setcounter{equation}{0}

In this section we present several features of $^*$-homomorphic
convolution cocycles on three types of example of $C^*$-bialgebras -
algebras of continuous functions on compact semigroups, universal
$C^*$-algebras of discrete groups and full compact quantum groups.
We focus on connections between the results obtained in this chapter
and the case of purely algebraic QS convolution cocycles analysed in
Chapter 3.

Whenever $\Alg$ is a (purely algebraic) $^*$-bialgebra, $D$ is a
dense subspace of $\kil$ and $\varphi \in L(\Alg; \Opdag(\Dhat))$ we
will follow the convention of Chapter 3 and write $l^{\varphi}$ for
a QS convolution cocycle on $\Alg$ satisfying the coalgebraic QS
differential equation \eqref{qgQSDE}. This hopefully will not cause
any confusion with $l^{\varphi} \in \QSC(\alg; \Exps)$ being a
solution of \eqref{CqgQSDE} for a $C^*$-bialgebra $\alg$ and a map
$\varphi\in CB \left(\alg; B(\khat)\right)$; it will be always clear
from the context which notion is understood.

\subsection*{Commutative case - continuous functions on a semigroup }

Let $G$ be a compact semigroup with identity $e$ and $\alg = C(G)$ denote the
algebra of all continuous functions on $G$. The algebra $\alg$ has the structure of a $C^*$-bialgebra
with coproduct given (with the help of the natural identification
$C(G)\ot C(G) \cong C(G\times G)$) by
\begin{equation} \Com(f) (s,t) = f(st), \;\;\; s,t \in G, \; f\in \alg, \label{comcop}\end{equation}
and counit given by
\[ \Cou (f) = f(e), \;\;\; f \in \alg.\]

Following the standard in quantum probability idea
(going back to \cite{AFL} and beyond),
any $G$-valued stochastic process $\{X_t: t \in \br_+\}$
on the probability space $(\Omega, \mathfrak{F}, \mu)$ can be dually described by
a family of unital $^*$-homomorphisms $\{l_t: t \in \br_+\}$ mapping $\alg$ into
$L_{\infty}(\Omega, \mathfrak{F},  \mu)$. These homomorphisms are given by
\[ l_t (f) := f \circ X_t, \;\; f \in C(G), t \geq 0,\]
and determine uniquely the original process.

Recall that a process $\{X_t: t \in \br_+\}$ on a semigroup with identity is called a
\emph{L\'evy process} if
it has identically distributed, independent increments, $\mu(\{X_0=e\})=1$ and
the distributions of $X_t$ weakly converge to the Dirac measure $\delta_{\{e\}}$
(the distribution of $X_0$) as $t$ tends to $0$. As well known, and often pointed in
the literature, not all L\'evy processes have
generators defined on the whole of $C(G)$. In our language, this fact corresponds to
the fact that not all $^*$-homomorphisms on $\alg$ are Markov-regular. Now Markov-regularity
of the process
corresponds to the norm continuity of the convolution semigroup given by
\[ \lambda_t (f) = \int_{\Omega} f \circ X_t d\mu , \;\; f \in C(G), t \geq 0 \]
(note that
standard weak continuity of this semigroup corresponds in the algebraic formulation
to the pointwise continuity of the Markov semigroup) . This yields the following:

\begin{propn}
Let $\{X_t: t\geq 0\}$ be a L\'evy processes on a normal compact
semigroup $G$. It is equivalent (in the sense of identical
finite-dimensional distributions, see Section \ref{QLevy}) to a
Markov-regular $^*$-homomorphic QS convolution cocycle on $\alg$ if
and only if it satisfies the following condition:
\begin{equation} \mu \left(\{X_t = e \} \right)\stackrel{t \to 0^+}{\longrightarrow} 1. \label{nor}\end{equation}
\end{propn}

\begin{proof}
It is easily seen that the condition \eqref{nor} implies the existence of the bounded
generator $\gamma: \alg \to \bc$ from which the process can be reconstructed. The other
direction can be seen by considering the Markov semigroup of a given QS convolution cocycle and
judiciously choosing continuous functions on $G$ with values in $[0,1]$ and
equal to $0$ outside of some neighbourhood of $e$ (using the assumption that $G$ is a normal
topological space).
\end{proof}

Processes satisfying (\ref{nor}) were investigated for example in \cite{Gren}. They are called
there homogenous processes of discontinuous type and their laws are shown to be compound Poisson
distributions (Theorem 2.3.5 of \cite{Gren}). Note that the results of Appendix A   on innerness of $\kil$-Sch\"urmann triples in conjunction with the definition of
quantum Poisson processes in \cite{franz} may be interpreted as the noncommutative counterpart of the above statement.

In general all L\'evy processes on a semigroup may be equivalently
realised (again in the sense of equal finite-dimensional
distributions) as quantum L\'evy processes on $^*$-bialgebras
(\cite{schu}, \cite{FranzSchott}).

\subsection*{Cocommutative case - group algebras}

Let $\Gamma$ be a discrete group. Denote by $\alg=C^*(\Gamma)$ the enveloping $C^*$-algebra
of the Banach algebra $l^1(\Gamma)$ (\cite{Pedersen}), called a universal
(or full) $C^*$-algebra of
$\Gamma$. By construction (as the algebra of functions on $\Gamma$
with finite support is dense in $\alg$), there exists a (so-called universal) unitary representation
$L:\Gamma \to \alg$ such that $\Alg:= \text{Lin}\{L_g : g \in \Gamma\}$ is dense in $\alg$.
Due to universality  the mappings $\Com$ and $\Cou$ defined on the image of $L$ by:
\[ \Com(L_g) = L_g \ot L_g, \;\; \Cou(L_g) =1, \;\; g \in \Gamma\]
extend to $^*$-homomorphisms on $\alg$. It is easy to check that $\alg$ equipped with the
comultiplication $\Com$ and the counit $\Cou$ becomes a cocommutative
$C^*$-bialgebra.

\begin{tw}
Let $\alg = C^*(\Gamma)$ for a discrete group $\Gamma$. Then the formula
\begin{equation}\label{W and l}
W(t,g) = l_t(L_g) \quad (g\in\Gamma, t\geq 0)
\end{equation}
defines a bijective correspondence   between unital $^*$-homomorphic
QS convolution cocycles $l$ on $\alg$ and maps $W:\br_+ \times\Gamma \to
B(\Fock)$ satisfying the following conditions:
\begin{rlist}
\item for each $g \in \Gamma$ the family $\{W(t,g): t \geq 0\}$ is a QS operator  cocycle;
\item for each $t \geq 0$ the family $\{W(t,g) :g \in \Gamma\}$ is a unitary representation of
$\Gamma$ on $\Fock$.
\end{rlist}
\label{corp}
\end{tw}

\begin{proof}
Assume that $l \in \Proc(\alg;\Fock)$ is a $^*$-homomorphic QS
convolution cocycle and define a map $W:\br_+\times\Gamma\to
B(\Fock)$ by \eqref{W and l}. Then, for all $t,s \geq 0,$ $g,h \in
\Gamma$
\begin{eqnarray*}
 W(t+s,g) &=& l_{t+s} (L_g) = (l_t \ot (\sigma_t \circ l_s)) \Com (L_g)  = l_t(L_g) \ot \sigma_t
(l_s (L_g)) \\ &=&  W(t,g) \ot \sigma_t (W(s,g)),\\
 W(t,g) W(t,h) &=& l_t(L_g) l_t(L_h) = l_t (L_g L_h) = l_t (L_{gh}) = W(t,gh),\\
 W(t,g)^* &=& l_t(L_g)^* = l_t(L_g^*) = l_t (L_{g^{-1}}) = W(t, g^{-1}),\\
 W(t,e) &=& l_t(L_e) = l_t (\ida) = I_{\Fock}, \;\;\; W(0,g) = j_0(L_g) = I_{\Fock}.
\end{eqnarray*}

Conversely, suppose that $W:\br_+ \times\Gamma \to B(\Fock)$ is a map satisfying the conditions
{\rm (i) and (ii)}. For each $t \geq 0$ define the linear map
$l_t: \Alg \to B(\Fock)$ by
\[ l_t(L_g) = W(t,g), \;\;\; g\in \Gamma.\]
Due to universality, the map $l_t$ continuously extends to a $^*$-homomorphism on
the whole $\alg$.
Properties of $W$ guarantee that for $a \in \Alg$ ($t, s \geq 0$),
\[ l_0 (a) = \Cou (a) I_{\Fock}, \;\; l_{t+s} (a) = (l_t \otimes \sigma_t (l_s) ) \Com (a),\]
and continuity of the extension assures that the above equations remain valid for any
$a \in \alg$.
\end{proof}

\noindent On the level of the stochastic generators the above correspondence takes the following form.

\begin{propn} \label{corg}
Let $\Alg:= \Lin\{L_g : g \in \Gamma\}$ for a discrete group
$\Gamma$. Then
\[
\psi_g = \varphi (L_g), \quad g\in\Gamma,
\]
determines a bijective correspondence between
maps $\varphi \in L( \Alg; B(\kilhat))$ satisfying
\begin{equation}
\varphi (ab) = \varphi (a) \Cou(b) + \Cou(a) \varphi(b) + \varphi(a)
       \QSproj \varphi(b), \;\; \varphi (a)^* = \varphi(a^*),
       \;\;\varphi(1) = 0, \label{pier}
\end{equation}
and maps $\psi : \Gamma \to B(\kilhat)$ satisfying
\begin{equation}    \label{drug}
 \psi_{gh} = \psi_g + \psi_h + \psi_g \QSproj \psi_h,  \;\;
(\psi_g)^* = \psi_{g^{-1}},\;\; \psi_e = 0;
\end{equation}
\end{propn}

\begin{proof} Elementary calculation.
 \end{proof}

\begin{rem}
Identities (\ref{drug}) may be considered as a special (time-independent)
case of formulae (4.2-4) in \cite{flows}. They are equivalent to $\psi$
having the block matrix form
\begin{equation}\label{lambda xi U}
\psi_g =
\begin{bmatrix}
i \lambda_g - \frac{1}{2} \|\xi_g\|^2 & - \langle \xi_g |U_g \cr
 |\xi_g\rangle & U_g - I_{\kil} \cr
\end{bmatrix},
\end{equation}
for a unitary representation $U$ of $\Gamma$ on $\kilhat$
and maps $\lambda: \Gamma \to \br$ and $\xi: \Gamma \to \kil$
satisfying
\[
\xi_{gh} = \xi_g + U_g \xi_h
\text{ and }
\lambda_{gh} =
\lambda_g + \lambda_h - {\rm Im} \langle \xi_g, U_g \xi_h \rangle.
\]
\end{rem}

Observe that according to Theorem \ref{HH} each map  $\varphi\in L( \Alg; B(\khat))$
satisfying (\ref{pier})
generates a unital, real and weakly multiplicative QS convolution cocycle
 $l^{\varphi}$ on $\Alg$. Further $l^{\varphi}$ continuously
extends to a $^*$-homomorphic QS convolution cocycle on $\alg$ (see
Lemma \ref{extend} below). On the other hand, given a map $\psi$
such as in (\ref{drug}), for each fixed $g \in \Gamma$ the QS
differential equation of the form
\[ W_0(g) = I_{\Fock}, \;\;\; dW_t(g) = \psi(g) W_t(g) d\Lambda_t\]
yields a unitary cocycle $\{W_t(g): t \geq 0\}$ (\cite{lwjfa}). The map $W:\br_+ \times \Gamma$ given by
$W(t,g)= W_t(g)$ satisfies the conditions of theorem \ref{corp}. One can easily see that
the correspondences described in Theorems \ref{corp} and \ref{corg} are consistent with this
construction.

\begin{propn}
A unital $^*$-homomorphic QSCC $l$ on $\Alg$ is equal to
$l^{\varphi}$ for some $\varphi\in L( \Alg; B(\khat))$ if and only
if it is weakly measurable.
\end{propn}
\begin{proof}
One direction is trivial. For the other consider the unitary cocycles associated with $l$ by
Theorem \ref{corp}. Theorem 6.7 of \cite{lwjfa} implies that each of these cocycles is
stochastically generated (as it is weakly measurable). Denoting respective generators by
$\psi_g$ one can see that the so obtained map $\psi:\Gamma \to B(\kil)$ satisfies the conditions
(\ref{drug}). Theorem \ref{corg} (and discussion above) imply the desired conclusion.
\end{proof}

If a $^*$-homomorphic QS convolution cocycle $l^{\varphi}$ on $\alg$
is Markov-regular, the automatic innerness of its stochastic
generator (Corollary \ref{innderiv}) implies in particular that the
triple $(\lambda, \xi, U)$ corresponding to $\varphi$ by Theorem
\ref{corg} must also be inner, in the following sense: there exists
a vector $\eta \in \kil$ such that
\[ \xi(g) =   U(g) \eta - \eta, \;\;\;\; \lambda(g) = {\rm Im} \langle \eta, U_g \eta \rangle.\]

Elements of a $C^*$-bialgebra $\blg$ are called \emph{group-like}
when they satisfy $\Com b = b \ot b$, as the elements $L_g$ do. On such elements
the solution $\big(k_t(b)\big)_{t\geq 0}$,
of the mapping QS differential equation \eqref{CqgQSDE},
is given by the solution of the operator QS differential equation
\[
dX_t = X_t d \Lambda_L (t), \;\;\; X_0 = I_{\Fock},
\]
where $L= \varphi(b) \in B(\kilhat)$. For more on this we refer to Section
4.1 of \cite{schu}.

\subsection*{Compact quantum groups}

The concept of compact quantum groups was introduced by S.L.\,Woronowicz, in \cite{CMP}. For our
purposes it is most convenient to adopt the following definition:

\begin{deft}[\cite{wor2}] \label{CQG}
A compact quantum group is a pair $(\alg, \Com)$, where $\alg$ is a
$C^*$-algebra, and $\Com:\alg \to \alg \ot \alg$ is a unital,
$^*$-homomorphic map which is coassociative and satisfies the
quantum cancellation properties:
\[ \overline{\Lin}((1\ot \alg)\Com(\alg) ) = \overline{\Lin}((\alg \ot 1)\Com(\alg) )
= \alg \ot \alg. \]
\end{deft}

For the concept of Hopf $^*$-algebras and their unitary
corepresentations, as well as unitary corepresentations of compact
quantum groups, we refer to Appendix B (see also \cite{Schm}). Here
it is sufficient to note the facts contained in the following
theorem.

\begin{tw} [\cite{CMP}]
Let $\alg$ be a compact quantum group and let $\Alg$ denote the
linear span of the matrix coefficients of irreducible unitary
corepresentations of $\alg$. Then $\Alg$ is a dense $^*$-subalgebra
of $\alg$, the coproduct of $\alg$ restricts to an algebraic
coproduct $\Com_0$ on $\Alg$ and there is a natural counit $\Cou$
and coinverse $\mathcal{S}$ on $\Alg$ which makes it a Hopf
$^*$-algebra.
\end{tw}

\begin{rem}[\cite{coamen}]
In the above theorem $(\Alg,\Com_0, \Cou, \mathcal{S})$ is the
unique dense Hopf $^*$-subalgebra of $\alg$, in the following sense:
if $(\Alg',\Com'_0, \Cou', \mathcal{S}')$ is a Hopf $^*$-algebra, in
which $\Alg'$ is a dense $^*$-subalgebra of $\alg$ and the coproduct
of $\alg$ restricts to an algebraic coproduct $\Com'_0$ on $\Alg'$,
then $(\Alg',\Com'_0, \Cou', \mathcal{S}')$ equals $(\Alg,\Com_0,
\Cou, \mathcal{S})$.
\end{rem}

The Hopf  $^*$-algebra arising here is called the \emph{associated
Hopf  $^*$-algebra} of $(\alg, \Com)$. When $\alg = C(G)$ for a
compact group $G$, $\Alg$ is the algebra of all matrix coefficients
of unitary representations of $G$; when $\alg$ is a universal
$C^*$-algebra of a discrete group $\Gamma$, $\Alg=\Lin \{L_g: g\in
\Gamma\}$ (see the beginning of the previous subsection).
M.\,Dijkhuizen and T.\,Koornwinder observed that the Hopf
$^*$-algebras arising in this way have intrinsic algebraic
structure.

\begin{deft}
A Hopf  $^*$-algebra $\Alg$ is called a CQG (compact quantum group)
algebra if it is the linear span of all matrix elements of its
finite dimensional unitary corepresentations.
\end{deft}

\begin{tw}[\cite{Koor}]
Each Hopf  $^*$-algebra associated with a compact quantum group is a
CQG algebra. Conversely, if  $\Alg$ is a CQG algebra then
\begin{equation}
\label{norm on a CQG}
\|a\|:=
\sup \big\{\|\pi(a)\|:
\pi \text{ is a $^*$-representation of } \Alg
\text{ on a Hilbert space}\big\}
\end{equation}
defines a $C^*$-norm on $\Alg$ and
the completion of $\Alg$ with respect to this norm is a compact
quantum group whose comultiplication extends that of
$\Alg$.
\end{tw}

The compact quantum group obtained in this theorem is called the
\emph{universal compact quantum group of $\Alg$} and is
denoted $\Alg_{\rm u}$.

For further use note the following extension of Lemma 11.31 in \cite{Schm}:

\begin{lemma}              \label{extend}
Let $E$ be a dense subset of a Hilbert space $H$ and let $\Alg$ be a CQG algebra.
Suppose that $\pi:\Alg \to \Opdag(\wh{E})$ is a real, unital and weakly
multiplicative (where the latter, as usual, means multiplicative
with respect to the product $\schur$ introduced in \eqref{bullet}). Then $\pi$ admits a
continuous extension to a unital $^*$-homomorphism from $\Alg_u$ to $B(H)$.
\end{lemma}

\begin{proof}
Let $(v_{i,j})_{i,j=1}^n$ be any finite dimensional unitary corepresentation of $\Alg$.
For any $i,j \in \{1, \ldots,n\}, \xi \in E$,
\[ \| \pi (v_{i,j}) \xi \|^2 \leq \sum_{k=1}^n \| \pi (v_{k,j}) \xi \|^2 =
\sum_{k=1}^n \la \pi(v_{k,j}) \xi , \pi(v_{k,j}) \xi \ra \]
\[ = \sum_{k=1}^n \la \xi , \pi (v_{k,j})^{\dagger} \pi(v_{k,j}) \xi \ra
=\left\la \xi , \pi \left( \sum_{k=1}^n  v_{k,j}^* v_{k,j} \right)\xi \right \ra =
\|\xi \|^2, \]
as $\sum_{k=1}^n  v_{k,j}^* v_{k,j} = 1_{\Alg}$. This implies that for each $a\in \Alg$
the operator $\pi(a)$ extends to a bounded map on $H$. A moment of reflection suffices to see that
so extended $\pi$ is a unital $^*$-homomorphism on $\Alg$. The definition of the canonical
norm on $\Alg$ implies that $\pi$ is contractive on $\Alg$ with respect to this norm
and as such may be continuously extended to a unital  representation of $\Alg_u$.
\end{proof}

\begin{deft}
A compact quantum group $(\alg, \Com)$ is called full
 if the $C^*$-norm
it induces on its associated CQG algebra $\Alg$
coincides with its canonical norm defined in \eqref{norm on a CQG}
(equivalently, if $\alg$ is $^*$-isomorphic to $\Alg_{\rm u}$).
\end{deft}

The notion of full compact quantum groups was introduced in \cite{coamen} and in \cite{Baaj}
(in the first paper they were called universal compact quantum groups).  It is very relevant for our
 context, as the above facts imply the following

\begin{propn}
Each full compact quantum group $\alg$ is a $C^*$-bialgebra (with
the counit being a continuous extension of the counit of its
associated Hopf $^*$-algebra $\Alg$). There is a bijective
correspondence between unital $^*$-homomorphic QS convolution
cocycles on $\alg$ and unital, real and weakly multiplicative QS
convolution cocycles on $\Alg$.
\end{propn}

Both families of examples described in previous two subsections - algebras of continuous
functions on compact groups and full $C^*$-algebras of discrete groups -  are full compact
quantum groups. Moreover most of the genuinely quantum (i.e.\ neither commutative nor
cocommutative) compact quantum groups considered in the literature also fall into this category
- among them the queen of all examples, quantum $SU_q(2)$.

Reconnecting with the contents of Chapter 3, we obtain

\begin{tw}
Let $k\in \Proc_{\mathrm{cb}}(\alg;\Exps)$, where $\alg$ is a full
compact quantum group with the associated Hopf $^*$-algebra $\Alg$.
Then the following are equivalent:
\begin{rlist}
\item  $k$ and $k^{\dagger}$ are H\"older continuous QS convolution cocycles;
\item $k|_{\Alg}=l^{\varphi}$ for some map $\varphi \in L(\Alg; B(\khat))$.
\end{rlist}
\end{tw}

\begin{proof}
One direction follows from the fact that $\Alg$ is an (algebraic) coalgebra and Theorem
\ref{Q}. The other is trivial.
\end{proof}

Specialising to $^*$-homomorphic cocycles yields a much stronger
result:

\begin{tw}
Let $k\in \Proc(\alg;\Exps_D)$, where $\alg$ is a full
compact quantum group with the associated Hopf $^*$-algebra $\Alg$
and $D$ is a dense subspace of $\kil$. Then the following are
equivalent:
\begin{rlist}
\item $k$ is H\"older continuous, unital and $^*$-homomorphic
      QS convolution cocycle;
\item $k$ is bounded and $k|_{\Alg}=l^{\varphi}$ for some
 $\varphi \in L(\Alg; \Opdag(\Dhat))$ satisfying the structure relations
 \eqref{alghomstruct} and vanishing at $1_{\Alg}$.
\end{rlist}
\end{tw}

\begin{proof}

(i)$\Rightarrow$(ii) \\
Follows from the previous theorem, Theorem \ref{Q}, and the implication (i)$\Rightarrow$(ii) of Theorem
\ref{HH}. Note that Theorem \ref{Q} yields even $\varphi \in L\left(\Alg; \Opdag(\khat)\right) =
L\left(\Alg; B(\khat)\right)$.

(ii)$\Rightarrow$(i) \\  Theorem \ref{HH} guarantees that
$l=k|_{\Alg}$ is real, unital, and weakly multiplicative. Lemma
\ref{extend} shows that $l$ admits a continuous extension to a
$^*$-homomorphic unital process in $\Proc(\Alg;\Exps_D)$ which must
coincide with $k$. Application of the previous theorem ends the
proof.
\end{proof}

 The above theorem could be equivalently formulated in the following way:

\begin{tw}
Let $k\in \Proc(\Alg;\Exps_D)$, where $\alg$ is a full
compact quantum group with associated Hopf $^*$-algebra $\Alg$ and
$D$ is a dense subspace of $\kil$. Then the following are
equivalent:
\begin{rlist}
\item  $k$ extends to a  H\"older continuous unital and $^*$-homomorphic
QS convolution cocycle on $\alg$;
\item  $k=l^{\varphi}$ for some
 $\varphi \in L(\Alg; \Opdag(\Dhat))$ satisfying the structure relations
 \eqref{alghomstruct} and vanishing at $1_{\Alg}$.
\end{rlist}
\end{tw}

\begin{rem}
In the course of the proof of the last theorem it was established that each map $\varphi$ defined on
a CQG algebra $\Alg$ with values in $\mathcal{O}^{\dagger}(\Dhat)$ satisfying the conditions
\eqref{alghomstruct} and vanishing at $1_{\Alg}$
must be bounded-operator-valued. However, $\varphi$ need not extend to $\alg$
(see examples in \cite{Skeide}). If it is continuous,
then it is necessarily completely bounded.
\end{rem}

\subsection*{Construction of $C^*$-hyperbialgebras by conditional expectations
   and QS convolution cocycles}

In previous subsections the variety of examples of $C^*$-bialgebras
were presented and $^*$-homomorphic QS convolution cocycles on these
were described and given alternative interpretation. Here the way of
obtaining new $C^*$-hyperbialgebras as certain subalgebras of other
$C^*$-hyperbialgebras is recalled. The construction was explicitly
described (in the context of compact quantum groups) in the papers
\cite{Kal1} and \cite{Kal2}, but its origins go back much further
(see \cite{ChaV} and references therein).
 In fact all known examples of noncommutative
$C^*$-hyperbialgebras arise in this way from $C^*$-bialgebras.

\begin{fact}   \label{expecthyper}
Let $(\alg, \Com, \Cou)$ be a $C^*$-hyperbialgebra. Assume that $\wt{\alg}$ is a unital
$C^*$-subalgebra of $\alg$ and that there exists a conditional expectation (that is norm-one
projection) $P$ from $\alg$  onto $\wt{\alg}$ satisfying the following identities:
\[ (P \ot \id_{\alg}) \Com P = (P \ot P) \Com = (\id_{\alg} \ot P) \Com.\]
Then $\wt{\alg}$ equipped with the coproduct $\wt{\Com}$ and the counit $\wt{\Cou}$, where
\[ \wt{\Com} = (P \ot P) \Com|_{\wt{\alg}}, \;\;\; \wt{\Cou} = \Cou|_{\wt{\alg}},\]
is a $C^*$-hyperbialgebra.
\end{fact}

Two particular cases of this construction are \emph{double coset bialgebras}
and \emph{Delsarte $C^*$-hyperbialgebras}; they are described below.

Let $(\alg_1,\Com_1,\Cou_1)$ and $(\alg_2,\Com_2,\Cou_2)$ be
$C^*$-bialgebras and assume that the latter is a \emph{quantum
subsemigroup} of the former.
 This means that there exists a unital $^*$-homomorphism
$\pi:\alg_1 \to \alg_2$ which is surjective and intertwines the coalgebraic structure:
\[ (\pi \ot \pi)\circ \Com_1 = \Com_2 \circ \pi, \;\; \Cou_2 \circ \pi = \Cou_1.\]
Assume additionally that $\alg_2$ admits a Haar state $h$ (recall formula \eqref{Haar}). Define the following $C^*$-subalgebras
of $\alg_1$ (\emph{the algebras of left and right cosets of $\alg_2$}):
\[\alg_1 /\alg_2 = \{ a \in \alg_1: (\id_{\alg_1 } \ot \pi ) \circ \Com_1 (a) = a \ot 1\},\]
\[ \alg_1 \backslash \alg_2 = \{ a \in \alg_1: (\pi \ot \id_{\alg_1 }  ) \circ \Com_1 (a) = 1 \ot a\}\]
and the \emph{double coset bialgebra}
\[ \alg_2 \backslash \alg_1 / \alg_2 = \alg_1 /\alg_2 \cap \alg_1 \backslash \alg_2.\]
It can be checked that the map $P:\alg_1 \to \wt{\alg}:=\alg_2 \backslash \alg_1 / \alg_2$ defined by
\[ P(a) = \left((h \circ \pi) \ot \id_{\alg_1} \ot (h \circ \pi) \right) \left(\Com_1 \ot \id_{\alg_1} \right)
 \Com_1 (a), \;\;\;  a \in \alg_1,\]
satisfies the conditions given in Fact \ref{expecthyper}. Its action may be understood as averaging (twice)
over the quantum subsemigroup; the construction is common in theory of classical hypergroups (\cite{hyp}).

Let $(\alg,\Com,\Cou)$ be a $C^*$-bialgebra and assume that a compact group $\Gamma$ acts
(continuously with respect to the topology of pointwise convergence)
on $\alg$ by $C^*$-algebra automorphisms satisfying
\[ (\gamma \ot \gamma) \circ \Com = \Com \circ \gamma,\;\; \Cou \circ \gamma = \Cou, \;\;\; \gamma \in \Gamma.\]
Let $\wt{\alg}$ be the fixed point subalgebra,
$\wt{\alg}= \{ a \in \alg: \forall_{\gamma \in \Gamma}\; \gamma(a) = a\}$.
It is easily checked that the map $P:\alg \to \wt{\alg}$ given by
\[ P(a) = \int_{\Gamma} \gamma(a) d\gamma, \;\;\; a \in \alg\]
(where $d\gamma$ denotes the normalised Haar measure on $\Gamma$), satisfies the assumptions of Fact \ref{expecthyper}.
The resulting $C^*$-hyperbialgebra is called a \emph{Delsarte $C^*$-hyperbialgebra}.

The connection between QS convolution cocycles on $\alg$ and $\wt{\alg}$ is given in the following fact:

\begin{fact} \label{corresp}
Assume that $(\wt{\alg}, \wt{\Com}, \wt{\Cou})$ is a $C^*$-hyperbialgebra
arising from a $C^*$-hyperbialgebra  $(\alg, \Com, \Cou)$ via the construction presented
in Fact \ref{expecthyper}, with the associated
conditional expectation $P$. Then there is a bijective correspondence between QS convolution
cocycles on $\wt{\alg}$ and $P$-invariant convolution increment processes on $\alg$ with
initial condition given by the idempotent functional $\Cou \circ P$.
\end{fact}

\begin{proof}
Assume first that  $\wt{l}\in \QSC(\wt{\alg};\Exps)$ and define $l \in \Proc(\alg;\Exps)$ by
\[ l_t = \wt{l}_t \circ P, \;\; t \geq 0.\]
Then clearly $l_0(a) = \Cou (P(a))$ for all $a \in \alg$, and $l$ is $P$-invariant. It remains to
check it is a convolution increment process. Choose $t,s \geq 0$ and compute:
\begin{eqnarray*} l_{s+t} &=& \wt{l}_{s+t} \circ P =
\left( \wt{l}_s \ot (\sigma_s \circ \wt{l}_t) \right) \wt{\Com} P  \\
&=& \left( \wt{l}_s \ot (\sigma_s \circ \wt{l}_t )\right) (P \ot P) \Com P =
\left( \wt{l}_s \ot (\sigma_s \circ \wt{l}_t) \right) (P \ot P) \Com \\
&=&  \left( l_s \ot (\sigma_s \circ l_t) \right) \Com.\end{eqnarray*}

Conversely, if $l \in \Proc(\alg;\Exps)$ is a $P$-invariant convolution increment process, with  initial
condition given by $\Cou \circ P$, then the process $\wt{l}\in \Proc(\wt{\alg};\Exps)$, defined simply by
the restriction of $l$ is a QS convolution cocycle on $\wt{\alg}$ -- again the only thing to be checked
is the convolution increment property: for all $s,t \geq 0$, $a \in \wt{\alg}$,
\begin{eqnarray*}\wt{l}_{s+t} (a) &=& l_{s+t} (a) = \left( l_s \ot (\sigma_s \circ l_t) \right) \Com (a) =
\left( l_s \ot (\sigma_s \circ l_t )\right) \Com P (a) \\
&=& \left( l_s \ot (\sigma_s \circ l_t) \right) (P \ot P)\Com P (a) = \left( \wt{l}_s \ot
(\sigma_s \circ \wt{l}_t) \right) \wt{\Com} (a).\end{eqnarray*}
\end{proof}

\begin{rem}
If $\Cou = \Cou \circ P$ (as is in the case of Delsarte $C^*$-hyperbialgebras, but usually not for double coset bialgebras), the processes
$l$ in the proof of the above theorem are obviously QS convolution cocycles. Assuming this is the case, and noting that Markov-regularity is clearly
preserved under the correspondence in the above fact, it is easily checked that if $\varphi\in CB(\wt{\alg}; B(\kilhat))$ then
$\wt{l}:=l^{\varphi} \in \QSC(\wt{\alg}; \Exps)$ corresponds to the process $l^{\psi} \in \QSC (\alg;\Exps)$ generated by
$\psi:= \varphi \circ P \in CB\left(\alg;B(\kilhat)\right)$. There is an analogous correspondence on the level of weak QS convolution cocycles.
\end{rem}

\section [Cocycles on locally compact quantum groups]  {Towards QS convolution cocycles on locally compact quantum groups} \label{multiplier}
\setcounter{equation}{0}

In this short section we discuss a possible approach to QS convolution cocycles on multiplier $C^*$-bialgebras.
No satisfactory results concerning the existence and characterisation of cocycles are known in this generality, and the whole section
should be considered as an announcement of a problem rather than the formulation of solutions.
The motivation for considering this question lies in the recently developed theory
of locally compact quantum groups.

\subsection*{Multiplier algebras and multiplier $C^*$-bialgebras}

Recall that if $\alg$ is a $C^*$-algebra then a closed (twosided, selfadjoint) ideal $\mathsf{I}$ of $\alg$ is called \emph{essential}
if for all $a \in \alg$ the equality $a \mathsf{I} = \{0\}$ implies $a=0$.

\begin{deft} Let $\alg$ be a $C^*$-algebra. The multiplier algebra
of $\alg$, denoted by $\bM(\alg)$, is the biggest $C^*$-algebra containing $\alg$ as
an essential ideal.
\end{deft}

The definition above requires certain comments. The biggest is understood in the following sense:
whenever $\blg$ is another $C^*$-algebra
containing  $\alg$ as an essential ideal, the identity map on $\alg$ extends to an injective $^*$-homomorphism from $\blg$ to $\bM(\alg)$.
Therefore the uniqueness of $\bM(\alg)$ up to an isomorphism follows directly from the definition. The easiest way to show the existence
is to exhibit a concrete model for $\bM(\alg)$. Below we describe a model of so-called \emph{double centralizers} (known also as \emph{multipliers}).

Let $\alg$ be a $C^*$-algebra. By a double centralizer on $\alg$ we understand a pair of maps $S,T:\alg\to\alg$ satisfying the following condition:
\[  T(a) b =a S(b), \;\;\; a,b \in \alg.\]
It is easy to see that both $S$ and $T$ must be linear module maps ($aT(b) = T(ab)$,
$S(a)b = S(ab)$ for all $a,b \in \alg$) and the Closed Graph Theorem implies they are continuous.
The vector space $\mathsf{DC}(\alg)$ of all double centralisers equipped with the norm $\| (S,T)\| := \|S\|(= \|T\|)$ is a Banach space.
Define the multiplication by $(S_1, T_1) \cdot (S_2, T_2):= (S_2 S_1, T_1 T_2) $ and the adjoint operation by
$(S,T)^* = (S',T')$, where $S'(a) = (T (a^*))^*$, $T'(a) = (S(a^*))^*$ for all $a \in \alg$ (note that
$(S,T)^* = (T^{\dag}, S^{\dag})$ in the earlier nomenclature). The norm introduced before is submultiplicative
and one may check it satisfies the $C^*$-condition, so that $\mathsf{DC}(\alg)$ is a $C^*$-algebra. The algebra $\alg$ may be embedded in
$\mathsf{DC}(\alg)$ via the identification $a \mapsto (L_a, R_a)$ (left and right multiplication by $a$). It may be shown that the image
of $\alg$ in this embedding is an essential ideal of $\mathsf{DC}(\alg)$, and that $\mathsf{DC}(\alg)$ satisfies the maximality condition
required in the definition of the multiplier algebra.

Note that the multiplier algebra $\bM(\alg)$ is always unital, and if $\alg$ is unital then $\bM(\alg) \cong \alg$.
 If $X$ is a locally compact
topological space and $\alg=C_0(X)$, then $\bM(\alg) \cong C_b(X) \cong C(\beta X)$, where
$C_0(X)$ denotes the algebra of continuous functions on $X$ vanishing at infinity, $C_b(X)$ the algebra
of continuous bounded functions and $\beta X$  the Stone-\v{C}ech compactification
of $X$.

Apart from the norm topology, there is another useful locally convex topology on $\bM(\alg)$, the so-called \emph{strict topology}. It is induced by the
family of seminorms $\{p_a:a\in\alg\}\cup \{p'_a:a\in\alg\}$ defined by $p_a (b) = \|ab\|$ and
 $p'_a(b) = \|ba\|$ for all $a \in \alg$ and
$b \in \bM(\alg)$. The unit ball of $\alg$ is strictly dense in the unit ball of $\bM(\alg)$, and \emph{strict} maps (that is bounded
linear maps which are strictly continuous on bounded subsets)
from $\alg$ to $\bM(\blg)$ extend uniquely to strict maps  on $\bM(\alg)$. The most useful examples of strict maps are
 \emph{slice maps} and \emph{nondegenerate} $^*$-homomorphisms. By a slice map is meant a map of the form
$\omega \ot \id_{\alg}:\alg \ot \alg \to \alg$, where  $\omega\in \alg^*$; a $^*$-homomorphism $\psi:\alg \to \bM(\blg)$
(where $\blg$ is another $C^*$-algebra) is
called nondegenerate if $\psi(\alg) \blg$ is dense in $\blg$. Note that nondegenerate $^*$-homomorphisms extend to
unital $^*$-homomorphisms from $\bM(\alg)$ to $\bM(\blg)$.
Elementary proofs of the well-known facts above may be found for example in the paper \cite{Johan}; for another presentation
of multiplier algebras, geared towards (and using the language of) the theory of Hilbert $C^*$-modules, see the book of \cite{Lance}.
In the latter book it is shown that a completely positive map is strict if and only if the image in this map of any (equivalently, every) approximate unit
in $\alg$ is Cauchy with respect to the strict topology in $\bM(\blg)$.

\begin{deft}
By a multiplier $C^*$-bialgebra is meant a $C^*$-algebra $\alg$ together with a nondegenerate $^*$-homomorphism
$\Com: \alg \to \bM(\alg \ot \alg)$ and a character $\Cou:\alg\to \bc$, satisfying the
standard coassociativity and counit identities:
\[ (\Com \ot \id_{\alg})\circ \Com = (\id_{\alg} \ot \Com) \circ \Com \;\;\; \text{ and } \;\;\; (\Cou \ot \id_{\alg}) \circ \Com =
(\id_{\alg} \ot \Cou) \circ \Com  = \id_{\alg}.\]
\end{deft}

The equalities above have to be understood in the only possible way, that is whenever it is necessary one in fact considers
continuous extensions of the strict maps in question (e.g.\ $\Com \ot \id_{\alg}$ is treated as a map from $\bM(\alg \ot\alg)$ to
$\bM(\bM(\alg \ot \alg) \ot \alg) \hookrightarrow \bM(\alg \ot \alg \ot \alg)$). The necessity of introducing such a definition stems
from the fact that already for classical, locally compact but noncompact group the comultiplication
may be defined on $C_0(G)$ by the formula \eqref{comcop}, but there is no guarantee that it will take values in
$C_0 (G \times G)$ - it rather maps into $C_b(G \times G)\cong \bM(C_0 (G) \ot C_0(G))$.
When a multiplier $C^*$-bialgebra $\alg$ is unital, it is a $C^*$-bialgebra.

\subsection*{QS convolution cocycles on multiplier $C^*$-bialgebras}

Let $\alg$ be a multiplier $C^*$-bialgebra.  The convolution product $\phi_1\star\phi_2$ retains
its meaning as long as the map $\phi_1\ot \phi_2$ has a natural extension to the multiplier algebra
$\bM(\alg \ot \alg)$. This is the case for example if both $\phi_1, \phi_2$
are bounded functionals ($\phi_1\ot \phi_2$ may be then written as a composition of slice maps). Therefore the notion of convolution semigroups of functionals may be
used in the context of multiplier $C^*$-bialgebras without any changes.

\begin{deft}
A process $l \in \Proc_{\wb} (\alg ; \Exps)$ is a weak convolution
increment process on a multiplier $C^*$-bialgebra $\alg$ if it
 satisfies formula \eqref{Csemigroup decomp} (for the associated
convolution semigroups given by \eqref{Cassoc semigroups}); it is
a weak QS convolution cocycle if in addition $l_0 (a) = \Cou(a) I_{\Fock}$ for all
$a \in \alg$.
\end{deft}

A satisfactory `strong' definition of a QS convolution cocycle is at
the moment available only for nondegenerate $^*$-homomorphic
cocycles.

\begin{deft}
A process $l\in \Proc_{\cb}(\alg; \Exps)$ is called a nondegenerate
$^*$-homomorphic QS convolution cocycle if the following conditions
are satisfied:
\begin{rlist}
\item for each $t \geq 0$ the map $l_t$ is a $^*$-homomorphism which is nondegenerate as a map from $\alg$ to $B(\Focktot)$;
\item for all $s,t \geq 0$
\begin{equation}
\label{MCconvinc}
l_{s+t} = l_s \star (\sigma_s \circ l_t);
\end{equation}
\item  for all $a\in \alg$
\[l_0 (a) = \Cou(a) I_{\Fock}.\]
\end{rlist}
\end{deft}

Observe that the formula \eqref{MCconvinc} makes sense;
the (completely bounded) map $l_s \ot (\sigma_s \circ l_t):\alg \ot\alg \to B(\Focktos) \ot B(\Fockstot)\subset B(\Focktot)$
is a nondegenerate $^*$-homomorphism, and so
extends in a natural way to a (unital) map from $\bM(\alg \ot \alg)$ to $B(\Focktot)$.

\begin{rem}
If $\alg$ is unital, nondegeneracy means simply that $l$ must be a
unital process; in any case, each nondegenerate $^*$-homomorphic QS
convolution cocycle on $\alg$ has a natural extension to a unital
$^*$-homomorphic process on $\bM(\alg)$. This may at first suggest
that one may reduce the considerations to QS convolution cocycles on
$\bM(\alg)$ -- note however that $\bM(\alg)$ need not be a
$C^*$-bialgebra in the sense of Definition \ref{Cbialg}. In general
the algebra $\bM(\alg) \ot \bM(\alg)$ is smaller than $\bM(\alg
\ot\alg)$.
\end{rem}

By a coalgebraic QS differential equation on $\alg$ (with
coefficient $\varphi \in CB\left(\alg;B(\khat)\right)$) is
understood again the equation of the form \eqref{CqgQSDE}. There is
no problem with understanding the notion of a weak solution: the
formula \eqref{wCqgQSDE} remains meaningful with the same
assumptions on the process as before. The strong solutions are more
problematic: recall that the first step of the corresponding
definition in Section \ref{coalgQSDEquat} required introducing the
map $\phi$ via the $R$-map. To do it in the context of multiplier
$C^*$-bialgebras one needs to insert additional assumption on
$\varphi$: it should be `stably strict', that is the map $\id_{\alg}
\ot \varphi:\alg \ot \alg \to \alg \ot B(\khat)$ should extend in a
natural way to  a map from $\bM(\alg \ot\alg)$ (in general expected
to have values in $\bM(\alg \ot B(\khat)$). Note that although
nondegenerate $^*$-homomorphisms are `stably strict', it is not
sufficient for our purposes. One would like at least to be able to
work with the maps $\varphi$ satisfying the structure relations
\eqref{structuremap}. The next step, a definition of a process $K$
which is to be actually QS integrated (see \eqref{integrproc})
causes further domain-related problems; despite certain efforts we
were not able to identify any natural sufficient conditions on
$\varphi$ allowing for completing this procedure.

As at the moment there is no satisfactory definition of strong
solutions of coalgebraic QS differential equation on multiplier
$C^*$-bialgebras, there are also no results on existence of
solutions of such equations (uniqueness of weak solutions can be
dealt with via methods similar to those of Section
\ref{coalgQSDEquat}). One of the basic questions would be: what are
the minimal assumptions on the coefficient $\varphi:\alg \to
B(\kilhat)$ so that some form of Picard iteration, as in Sections
\ref{QSDEOS} and \ref{coalgQSDEquat} is possible? The main problem,
as should be apparent from the above considerations, lies in an
insufficient understanding of the interplay between operator-space
theoretic notions
and the theory of multiplier $C^*$-algebras.

\renewcommand{\theequation}{\Alph{chapter}.\arabic{equation}}

\renewcommand{\chaptername}{Appendix}

\renewcommand{\thechapter}{\Alph{chapter}}
\setcounter{chapter}{0}  
\setcounter{propn}{0}

\renewcommand{\thepropn}{\Alph{chapter}.\arabic{propn}}

\chapter{$(\pi_1,\pi_2)$-derivations}

In this appendix we give an extension of the
innerness theorem of E.\,Christensen, for completely bounded derivations on a
$C^*$-algebra, to $(\pi_1,\pi_2)$-derivations, and prove automatic complete
boundedness for $(\pi,\chi)$-derivations, when $\chi$ is a character.

\begin{deft}
Let $\alg$ be a $C^*$-algebra and $\Xlg$ be a Banach $\alg$-bimodule.
A map
$\delta:\alg \to \Xlg$ is called a derivation if for all $a,b \in \alg$
\[ \delta(ab) = \delta(a) b + a \delta(b).\]
A derivation $\delta$ is inner if there exists an element
$x\in \Xlg$ such that for all $a \in \alg$
\[ \delta (a) = ax - xa.\]
\end{deft}

The following theorem is a modification (due to J.R.\,Ringrose) of the celebrated
Sakai theorem on boundedness of $C^*$-algebra derivations.

\begin{tw} [\cite{Sakai}, \cite{Ringrose}] \label{Ringrose}
Let $\alg$ be a $C^*$-algebra and $\Xlg$ be a Banach $\alg$-bimodule.
 Every derivation
$\delta:\alg \to \Xlg$ is bounded.
\end{tw}

For certain $\alg$-bimodules inner derivations are characterised by their complete boundedness. For a
simple proof of the next theorem, and connections with not necessarily real
homomorphisms between $C^*$-algebras we refer to \cite{Pisiersim}.

\begin{tw} [\cite{Der}]    \label{Christens}
Let $\alg\subset B(\hil)$ be a $C^*$-algebra, and let $\delta: \alg \to B(\hil)$
be a derivation. Then $\delta$ is completely bounded if and only if it is inner.
\end{tw}

For the purpose of this work
we are interested in the particular class of Banach $\alg$-bimodule-valued
derivations captured by the following definition.

\begin{deft}
Let $\alg$ be a $C^*$-algebra with representations $(\pi_1, \hil_1)$, $(\pi_2, \hil_2)$.
A map $\delta: \alg \to B(\hil_2;\hil_1)$ is called a $(\pi_1,\pi_2)$-derivation
if for all $a,b \in \alg$
\[ \delta(ab)  = \pi_1(a) \delta(b) + \delta(a) \pi_2(b).\]
A $(\pi_1,\pi_2)$-derivation $\delta$ is inner if it is implemented by an operator
$T \in B(\hil_2;\hil_1)$, in the sense that for all $a \in \alg$
\[ \delta (a)  = \pi_1(a) T - T \pi_2 (a).\]
\end{deft}

\begin{tw} \label{RingChristen}
Let $\alg$ be a $C^*$-algebra with representations $(\pi_1, \hil_1)$, $(\pi_2, \hil_2)$
and let $\delta: \alg \to B(\hil_1;\hil_2)$ be a $(\pi_1,\pi_2)$-derivation. Then
$\delta$ is completely bounded if and only if it is inner.
\end{tw}

\begin{proof}
Let $(\rho, \hil)$ be a faithful representation of $\alg$ and set
$\Hil = \hil_2 \oplus \hil_1 \oplus \hil$,
$\wt{\alg} = (\pi_2 \oplus \pi_1 \oplus \rho) (\alg)$. Then $\wt{\alg}$ is a $C^*$-subalgebra
of $B(\Hil)$, and a map $\wt{\delta} : \alg \to B(\Hil)$ defined
by ($a \in \alg$)
\[ \wt{\delta} \left(\left( \begin{array}{ccc} \pi_2(a) & & \\
  & \pi_1(a) & \\ & & \rho(a) \end{array} \right)\right) =
  \left( \begin{array}{ccc}   0 & & \\ \delta(a) & 0 & \\ & & 0 \end{array} \right).
\]
For any $a \in \alg$ write $\wt{a}$ for a relevant element of $\wt{\alg}$. Then
\begin{align*}  \wt{\delta} (\wt{a}) \wt{b} + \wt{a} \wt{\delta} (\wt{b}) & =
 \left( \begin{array}{ccc}   0 & & \\ \delta(a) & 0 & \\ & & 0 \end{array} \right)
 \left( \begin{array}{ccc} \pi_2(b) & & \\
  & \pi_1(b) & \\ & & \rho(b) \end{array} \right)\\
& + \left( \begin{array}{ccc} \pi_2(a) & & \\  & \pi_1(a) & \\ & & \rho(a) \end{array} \right)
  \left( \begin{array}{ccc}   0 & & \\ \delta(b) & 0 & \\ & & 0 \end{array} \right)  \\
& =  \left( \begin{array}{ccc}   0 & & \\ \delta(a)\pi_2(b) + \pi_1(a) \delta(b)  & 0 & \\ & & 0 \end{array} \right)
= \left( \begin{array}{ccc}   0 & & \\ \delta(ab)  & 0 & \\ & & 0 \end{array} \right)\\
& = \wt{\delta} (\wt{ab}) = \wt{\delta} (\wt{a} \wt{b}) ,     \end{align*}
so $\wt{\delta}$ is a derivation. Observe that  $\delta$ is inner if and only if $\wt{\delta}$
is; if $\wt{\delta}$ is implemented by $T\in B(\Hil)$, $\delta$ is implemented by
$P_{\hil_1} T P_{\hil_2}$ (a (2,1) `matrix coefficient' of $T$). The claim now follows from
Theorem \ref{Christens}.
\end{proof}

\begin{tw} \label{cbderiv}
Let $\alg$ be a  $C^*$-algebra with a  representation $(\pi, \hil)$ and
a character (nonzero multiplicative functional)
$\chi$. Then every $(\pi, \chi)$-derivation is completely bounded.
\end{tw}

\begin{proof}
Let $\delta: \alg \to B(\bc;\hil)=|\hil\ra$ be a $(\pi, \chi)$-derivation.
By Theorem  \ref{Ringrose} $\delta$ is bounded.
Without loss of generality we may suppose that the $C^*$-algebra $\alg$ and
representation $\pi$ are both unital; if necessary by extending $\pi$,
$\chi$ and $\delta$ to the unitisation of $\alg$ in the natural
way:
\[
(a,z)\mapsto \pi (a) + zI_{\hil}, \ \
(a,z)\mapsto \chi (a) + z \ \text{ and } \
(a,z)\mapsto \delta (a).
\]
Denote by $K$ a kernel of $\chi$, and let $P:\alg \to K$ denote
the canonical (bounded) projection onto $K$:
\[ P(a) = a - \chi (a) 1, \;\; a \in \alg. \]
Then $K$ is a maximal ideal, and as such also
a (nonunital) $C^*$ -subalgebra of $\alg$. Let $\{u_i:i \in I\}$ be an approximate unit of
$K$ (contained in the unit ball). Its existence guarantees that the set
$K_2 = {\rm lin}\{ab: a, b \in K\}$ is dense in $K$. For any $n \in \bn$, $a_1, \ldots, a_n,
b_1, \ldots, b_n \in K$
\begin{align*} \sum_{k=1}^n  \delta(a_k^*)^* \delta (b_k)  & = \lim_{i \in I}
\left(\sum_{k=1}^n \delta (a_k^* u_i)^* \delta(b_k u_i)  \right)
 = \lim_{i \in I}  \left(\sum_{k=1}^n (\pi(a_k^*) \delta (u_i))^* \pi(b_k)
\delta (u_i)  \right) \\ & =
\lim_{i \in I}   \delta(u_i)^* \pi (\sum_{k=1}^n a_k b_k) \delta(u_i),\end{align*}
so
\[ \left| \sum_{k=1}^n  \delta(a_k^*)^* \delta (b_k)  \right| \leq
 \| \delta\|^2 \| \sum_{k=1}^n a_k b_k \|.\]
This allows us to define a bounded functional $\wt{\gamma}: K_2 \to \bc$ by
\[\wt{\gamma} ( \sum_{k=1}^n a_k b_k) = \sum_{k=1}^n \delta(a_k^*)^* \delta (b_k) \]
(its continuous extension to $K$ will be extended by the same letter). Put then
\[ \gamma = \wt{\gamma} \circ P.\]
For any $a, b \in \alg$
\[ \delta(a)^* \delta(b)  =  \delta (a - \chi(a) 1)^*
\delta(b - \chi(b) 1)  = \wt{\gamma}((a^* - \chi(a^*) 1)(b - \chi(b)1)) \]
\begin{equation} \label{expres} = \gamma ((a^* - \chi(a^*) 1)(b - \chi(b)1)) =  \gamma(a^* b) - \gamma(a^*)
\chi(b) - \chi(a^*) \gamma(b),  \end{equation}
where we used the fact that $\delta(1)=0$. The last formula has a natural matricial equivalent. Fix $ n \in \bn$ and let
$\delta_n: M_n (\alg) \to M_n (B(\bc; H)) = B(\bc^n; H^n)$, $\gamma_n: M_n(\alg) \to M_n(\bc)$ and
$\chi_n: M_n(\alg) \to M_n(\bc)$ be respective tensorizations of $\delta$, $\gamma$ and
$\chi$.
Then for any  $\wt{a}=[a_{i,j}] \in M_n(A)$, $\wt{b}=[b_{i,j}] \in M_n(A)$
\begin{equation} \delta_n(\wt{a})^* \delta_n(\wt{b}) = \gamma_n (\wt{a}^* \wt{b}) +
\gamma_n (\wt{a}^*) \chi_n (\wt{b}) + \chi_n (\wt{a}^*) \gamma_n (\wt{b})
\label{matexpres} \end{equation}
Indeed, for any $i,j \in \{1, \ldots,n\}$  equation
\eqref{expres} implies
\begin{eqnarray*} (\delta_n (\wt{a})^* \delta_n (\wt{b}))_{i,j} &=& \sum_{k=1}^n  \delta(a_{k,i})^*
\delta(b_{k,j})= \\
 && \sum_{k=1}^n  \left(\gamma(a_{k,i}^* b_{k,j}) - \gamma(a_{k,i}^*)
\chi(b_{k,j}) - \chi(a_{k,i}^*) \gamma(b_{k,j}) \right).\end{eqnarray*}
Complete boundedness of $\delta$ follows easily from the estimate:
\[ \| \delta_n\|^2 \leq \|\gamma_n\| + 2 \|\gamma_n\|\|\chi_n\|=
  3 \|\gamma\|\]
 and arbitrariness of $n \in \bn$.
\end{proof}

Theorems \ref{RingChristen} and \ref{cbderiv} yield the following corollaries:

\begin{cor} \label{innderiv}
Let $\alg$ be a  $C^*$-algebra with  with a
representation $(\pi, \hil)$ and a character $\chi$.
Then every $(\pi,\chi)$-derivation is inner.
\end{cor}

\begin{cor} \label{disap}
If $\alg$ is a $C^*$-algebra with a character $\chi$ then every $(\chi,\chi)$-derivation
on $\alg$ vanishes.
\end{cor}

\chapter{Hopf $^*$-algebras and their corepresentations}

\setcounter{propn}{0}

In this appendix we present the basics of the theory of Hopf $^*$-algebras and their
corepresentations, following the presentation in \cite{Schm}.
In the second part we give the relevant definitions for compact quantum groups.

\begin{deft}
A unital bialgebra $\Alg$ (with multiplication $m$)
 is called a Hopf algebra if there exists a linear mapping $\mathcal{S}: \Alg \to \Alg$, called the antipode or the coinverse
of $\Alg$, such that for all $a \in \Alg$
\begin{equation} m (\mathcal{S} \ot \id_{\Alg} ) \Com (a) = m (\id_{\Alg} \ot \mathcal{S} ) \Com (a)= \Cou(a)1.
\label{antipode}\end{equation}
\end{deft}

Note that the antipode may be viewed as a convolution inverse of the identity mapping; in the convolution notation \eqref{antipode} takes form
\[ (\Sant \star \id_{\Alg}) (a) = (\id_{\Alg} \star \Sant) (a)= \Cou(a), \;\;\; a \in \Alg.\]

\begin{propn}
The antipode of a Hopf algebra is necessarily a unital antihomomorphism and a counital anti-coalgebra morphism, that is
\[ \Sant(1)=1, \;\;\; \Sant(ab) = \Sant(b) \Sant(a), \;\; a,b \in \Alg, \]
and
\[ \Cou \circ \Sant = \Cou, \;\;\; \Com \circ \Sant = \tau \circ (\Sant \ot \Sant) \circ \Com,\]
where $\tau$ denotes the tensor flip on $\Alg \odot \Alg$.
\end{propn}

\begin{deft}
A unital $^*$-bialgebra is called a Hopf $^*$-algebra if it is a Hopf algebra.
\end{deft}

Although the definition of a Hopf $^*$-algebra does not specify the behaviour of
the antipode under the involution, it follows from the definitions that for all $a \in \Alg$
\[ \Sant(\Sant(a^*)^*) = a.\]
In particular, the antipode of a Hopf $^*$-algebra is always invertible.

\begin{deft}
Let $\Coalg$ be a coalgebra. A corepresentation of $\Coalg$ on a vector space $V$ is a linear map
$\varphi:V \to V \odot \Alg$ such that
\[ (\id_{V} \ot \Com) \circ \varphi = (\varphi \ot \id_{\Alg} ) \circ \varphi, \;\;\; (\id_V \ot \Cou) \circ \varphi = \id_V.\]
A linear subspace $W \subset V$ is called $\varphi$-invariant if $\varphi(W) \subset W \odot \Alg$. Then the restriction of $\varphi$
to $W$ is also a corepresentation of $\Coalg$, called a subcorepresentation of $\varphi$.
A corepresentation $\varphi$  is called irreducible if it does not have any  nontrivial
(that is different from itself and from zero)
subcorepresentations.
\end{deft}

It is obvious that the coproduct of $\Coalg$ is itself a corepresentation of $\Coalg$ (on $\Coalg$). The following fact may therefore be viewed as
a generalisation of the Fundamental Theorem on Coalgebras:

\begin{fact}
If $\varphi$ is a corepresentation of a coalgebra $\Coalg$ on a vector space $V$, then any element of $V$ is contained in a finite-dimensional
$\varphi$-invariant subspace.
\end{fact}

Assume now that $V$ is a finite-dimensional Hilbert space and let $(e_1, \ldots, e_n)$ be an orthonormal basis of $V$. Whenever
$\varphi$ is a representation of a Hopf $^*$-algebra $\Alg$ on $V$, there exist uniquely determined
elements $v_{ij}\in \Alg$ ($i,j=1, \ldots, n$), called the \emph{matrix coefficients of $\varphi$ with respect to the basis
$(e_1, \ldots,e_n)$}, such that for each $j=1, \ldots, n$
\[ \varphi(e_j) = \sum_{i=1}^n e_i \ot v_{ij}.\]

\begin{deft}
A corepresentation of a Hopf $^*$-algebra $\Alg$ on a finite-dimensional Hilbert space $V$ is called unitary if for some (equivalently, for any)
orthonormal basis of $V$ the matrix coefficients of $\varphi$ with respect to this basis satisfy the conditions:
\begin{equation} \label{unitcorep} \Sant (v_{ij}) = (v_{ji})^*, \;\;\; i,j = 1, \ldots, n.\end{equation}
\end{deft}

The conditions \eqref{unitcorep} are equivalent to the equality $v^*v =v v^* = I_{M_n(\Alg)}$,
where $v=(v_{ij})_{i,j=1}^n \in M_n(\Alg)$.
This explains the motivation
behind the above definition.

Recall the definition of a compact quantum group (Definition \ref{CQG}). The notion of unitary corepresentations in the context of compact
quantum groups needs to be slightly different from the one given below, as in principle we do not want to exclude infinite-dimensional corepresentations.
Let $K(\Hil)$ denote the $C^*$-algebra of all compact operators on a Hilbert space $\Hil$
and recall that for any $C^*$-algebra $\blg$ its multiplier algebra (see Section \ref{multiplier}) is denoted by $\bM(\blg)$.

\begin{deft}
Let $(\alg, \Com)$ be a compact quantum group, and $\Hil$ be a Hilbert space. By a unitary corepresentation of $\alg$ acting on
$\Hil$ is understood a unitary $U \in \bM(K(\Hil) \ot \alg)$ satisfying the condition
\[ (\id_{\Hil} \ot \Com) U = U_{12} U_{13}\]
(in the above formula operators $U_{12}, U_{13} \in \bM(K(\Hil) \ot \alg \ot \alg)$ are defined via the so-called leg notation, and certain
natural extensions of maps in question are implicitly understood). A corepresentation $U$ is called irreducible
if for any projection $P \in B(\Hil)$ the commutation relation
$(P \ot 1_{\alg} ) U = U (P\ot 1_{\alg})$ implies $P=0$ or $P=1_{B(\Hil)}$.
\end{deft}

S.\,L.\,Woronowicz proved that
every irreducible corepresentation of a compact quantum group must be finite-dimensional (that is, $\Hil$ is finite-dimensional). Note
that in such a case $K(\Hil)=B(\Hil)$ is unital, there is no need to consider multiplier algebras and  unitaries
$U_{12}$ and $U_{13}$ in $B(\Hil) \ot \alg \ot \alg$ are defined simply by
\[ U_{12} = U \ot 1_{\alg}, \; \; \; U_{13} = (\id_{\Hil} \ot \tau) (U \ot 1_{\alg}) (\id_{\Hil} \ot \tau),\] where
$\tau:\alg \ot\alg \to \alg \ot \alg$ is the tensor flip.

To complete the list of notions used in Section \ref{Exam}  one more definition is needed.

\begin{deft}
Let $U$ be a finite-dimensional unitary corepresentation of a compact quantum group $\alg$ acting on a Hilbert space
$\Hil$. By matrix coefficients of $U$ with respect to an orthonormal basis
$(e_1, \ldots, e_n)$ of $\Hil$ are understood  elements $v_{ij} \in \alg$ defined by
\[v_{i,j} = (\la e_i| \ot 1_{\alg}) U (|e_j \ra \ot 1_{\alg}), \;\; i,j =1, \ldots n.\]
\end{deft}

One may check that the notions of corepresentations, their irreducibility and matrix coefficients introduced above are consistent
for finite-dimensional compact quantum groups (which are also Hopf $^*$-algebras). For more information on the vast topic
of topological quantum groups we refer to \cite{kus} and references therein.

\addcontentsline{toc}{chapter}
         {\protect\numberline{}{Bibliography}}

\small
\singlespacing


\begin{thebibliography}{GLSW}

\bibitem [Acc] {Luigi} L.\,Accardi, On the quantum Feynman-Kac formula,
\emph{Rend.\,Sem. Mat.\,Fis.\,Milano} \textbf{48} (1978), 135--180.

\bibitem [AFL]{AFL}
L.\,Accardi, A.\,Frigerio and J.\,Lewis, Quantum stochastic processes,
\emph{Publ. \!Res. \!Inst. \!Math. \!Sci.} \textbf{18} (1982) no. 1,
97--133.

\bibitem[Arv]{Arv}
W.\,Arveson, ``Noncommutative dynamics and $E$-semigroups," Springer-Verlag, New York 2003.


\bibitem[ASW]{asw}
L.\,Accardi, M.\,Sch\"{u}rmann and W.\,von Waldenfels,
Quantum independent increment processes on superalgebras,
\emph{Math.\,Z.} \textbf{198} (1988) no.\,4, 451--477.


\bibitem[BaS] {Baaj}
S.Baaj and G.Skandalis, Unitaires multiplicatifs et dualit\'e pour les produits crois\'es de $C^*$-alg\`ebres,
 \emph{Ann.\,Sci.\,\'Ecole Norm.\,Sup.} \textbf{26}  (1993)  no.\,4, 425--488.


\bibitem[Bel] {Slava}
V.\,Belavkin, Quantum stochastic positive evolutions: characterization, construction, dilation,
\emph{Comm.\,Math.\,Phys.}  \textbf{184}  (1997)  no.\,3, 533--566.


\bibitem[Bha] {Raja} B.V.R.\,Bhat,  Dilations, cocycles and product systems,
\emph{in} \cite{qiip}, Volume I.


\bibitem[Bia] {Biane}
P.\,Biane,  ``Calcul stochastique non-commutatif,'' \emph{in}
Lectures on probability theory (Saint-Flour, 1993), Lecture Notes in
Mathematics, \textbf{1608}, Springer, Berlin, 1995, pp.\,1--96.


\bibitem[BlH] {hyp}
W.R.\,Bloom and H.\,Heyer, ``Harmonic analysis of probability measures on hypergroups,''
de Gruyter Studies in Mathematics, 20. Walter de Gruyter \& Co., Berlin, 1995.


\bibitem[BMT]{coamen}
E.\,Bedos, G.\,Murphy and L.\,Tuset,
Co-amenability for compact quantum groups,
\emph{J.\,Geom.\,Phys.} \textbf{40} (2001) no.\,2, 130--153.

\bibitem[ChV]{ChaV}
Yu.\,Chapovsky and L.\,Vainerman,
Compact quantum hypergroups,
\emph{J.\,Operator Theory} \textbf{41} (1999) no.\,2, 261--289.

\bibitem[Chr]{Der}
E.\,Christensen,
Extensions of derivations II,
\emph{Math.\,Scand.}
\textbf{50} (1982), 111-122.

\bibitem[Con]{nongeo} A.\,Connes, ``Noncommutative geometry,'' Academic Press, Inc., San
Diego, CA, 1994.

\bibitem[Dav]{Davies}
E.B.\,Davies,
``One-parameter Semigroups,''
Academic Press, London 1980.

\bibitem[DiK] {Koor}
M.\,Dijkhuizen and T.\,Koornwinder,
CQG algebras -- a direct algebraic approach to compact quantum groups,
\emph{Lett.\,Math.\,Phys.}
\textbf{32} (1994) no.\,4, 315--330.

\bibitem[EfR] {ERuan}
E.G.\,Effros  and Z.J.\,Ruan, ``Operator Spaces,'' Oxford University Press,
Oxford 2000.



\bibitem[Fra]{franz}
U.\,Franz,
L\'evy processes on quantum groups and dual groups,
\emph{in}~\cite{qiip}, Volume II.

\bibitem[FSc]{FranzSchott}
U.\,Franz and R.\,Schott,
``Stochastic Processes and Operator Calculus on Quantum Groups,''
Mathematics and its Applications \textbf{490}, Kluwer, Dordrecht 1999.

\bibitem[FrS] {UweSch}
U.\,Franz and M.\,Sch\"{u}rmann,
L\'evy processes on quantum
hypergroups, \emph{in}
``Infinite Dimensional
Harmonic Analysis,''
\emph{eds.\ H. \!Heyer, T. \!Hirai \& N. \!Obata}
Gr\"abner, Altendorff 2000,
pp. \!93--114.

\bibitem
[Glo] {glo} P.\,Glockner, Quantum stochastic differential equations
on $^*$-bialgebras, \emph{Math.\,Proc. Camb.\,Phil.\,Soc.}
\textbf{109} (1991) no.\,3, 571--595.

\bibitem[GLW] {Stine}
D.\,Goswami, J.M.\,Lindsay and S.J.\,Wills, A stochastic Stinespring theorem,
\emph{Math.\,Ann.}
\textbf{319} (2001) no.\,4,  647--673.

\bibitem[GLSW] {dilate}
D.\,Goswami, J.M.\,Lindsay, K.B.\,Sinha  and S.J.\,Wills, Dilation of Markovian cocycles
on a von Neumann algebra, \emph{Pacific.\,J.\, Math.}
\textbf{211} (2003), 221--247.

\bibitem[Gre]{Gren}
U.\,Grenander,
``Probabilities on Algebraic Structures,''
John Wiley \& Sons,
New-York-London 1963.


\bibitem[Gui] {Guichardet}
A.\,Guichardet, ``Symmetric Hilbert Spaces and Related Topics,'' Lecture Notes in Mathematics
\textbf{267}, Springer, Heidelberg 1970.

\bibitem [HKK] {hkk}
J.\,Hellmich, C.\,K\"{o}stler and B.\,K\"{u}mmerer,
Noncommutative continuous Bernoulli shifts,
\emph{preprint, Queen's University, Kingston}, \#math.OA/0411565.

\bibitem[He1] {Heyercomp}
H.\,Heyer,
``Infinitely divisible probability measures on compact groups,'' in
Lectures on operator algebras (dedicated to the memory of David M. Topping; Tulane Univ. Ring and Operator Theory Year, 1970-1971, Vol. II),
Lecture Notes in Mathematics {\bf 247}, Springer-Verlag, Berlin 1972.

\bibitem[He2] {Heyer}
H.\,Heyer,
``Probability Measures on Locally Compact Groups,''
Springer-Verlag, Berlin 1977.

\bibitem[He3] {Heyer2}
H.\,Heyer, ``Structural aspects in the theory of probability.
A primer in probabilities on algebraic-topological structures," Series on Multivariate Analysis \textbf{7},
World Scientific Publishing Co., River Edge, NJ, 2004.

\bibitem[Hud] {Hud}
R.\,L.\,Hudson,
An introduction to quantum stochastic calculus and some of its
applications,
\emph{in}
``Quantum Probability Communications XI,''
\emph{eds.\ S. \!Attal \& J.M. \!Lindsay},
World Scientific, Singapore 2003,
pp. \!221-271.

\bibitem[HLP]{flows}
R.L.\,Hudson, J.M.\,Lindsay and K.R.\, Parthasarathy, Flows of quantum noise,
 \emph{J.\,Appl.\,Anal.}  \textbf{4}  (1998)  no.\,2, 143--160.


\bibitem[Hun]{Hunt}
G.A.\,Hunt, Semi-groups of measures on Lie groups,  \emph{Trans.\,Amer.\,Math.\,Soc.} \textbf{81} (1956), 264--293.

\bibitem[Kal] {Kal1}
A.A.\,Kalyuzhnyi,
Conditional expectations on compact quantum groups and new examples of quantum hypergroups, \emph{Methods Funct.\,Anal.\,Topology}
\textbf{7}  (2001)  no.\,4, 49--68.


\bibitem[KaCh] {Kal2}
A.A.\,Kalyuzhnyi and Yu.A.\,Chapovsky,
A factorization of conditional expectations on Kac algebras and quantum double coset hypergroups,
\emph{Ukrain.\,Mat.\,Zh.}  \textbf{55}  (2003)  no.\,12, 1669--1677;  translation in  \emph{Ukrainian Math.\, J.} \textbf{55}
(2003)  no. 12, 1994--2005.

\bibitem[KlS] {Schm}
A.\,Klimyk and K.\,Schm\"udgen,
``Quantum Groups and their Representations,''
Texts and Monographs in Physics,
Springer-Verlag, Berlin 1997.


\bibitem[Ku1]{Johan}
J.\,Kustermans,
One-parameter representations on $C^*$-algebras,
\emph{Preprint Odense Universitet}, \#funct-an/9707010.



\bibitem[Ku2]{kus}
J.\,Kustermans,
Locally compact quantum groups,
\emph{in} \cite{qiip}, Volume I.

\bibitem[Lan]{Lance}
E.C.\,Lance, ``Hilbert $C^*$-modules. A toolkit for operator algebraists.''
London Mathematical Society Lecture Note Series \textbf{210},
Cambridge University Press, Cambridge 1995.


\bibitem[Lin]{lect}
J.M.\,Lindsay,
Quantum stochastic analysis -- an introduction,
\emph{in} \cite{qiip}, Volume I.


\bibitem[LiP]{lp}
J.M.\,Lindsay and K.R.\,Parthasarathy,
On the generators of quantum stochastic flows,
\emph{J.\,Funct.\,Anal.}
\textbf{158} (1998) no.\,2, 521--549.

\bibitem[LiS]{MartinKalyan}
J.M.\,Lindsay and K.B.\,Sinha,
Feynman-Kac representation of some noncommutative elliptic operators,
\emph{J.\,Funct.\,Anal.}
\textbf{147} (1997) no.\,2, 400--419.

\bibitem[$\text{LS}_1$]
{LSqscc1}
J.M.\,Lindsay and A.G.\,Skalski,
Quantum stochastic convolution cocycles I,
\emph{Ann. \!Inst. \!H. \!Poincar\'{e}, Probab. \!Statist.}
\textbf{41} (2005), no.\,3  (En hommage a Paul Andr\'e Meyer),
581--604.

\bibitem[$\text{LS}_2$]
{LSbedlewo}
J.M.\,Lindsay and A.G.\,Skalski,
Quantum stochastic convolution cocycles---algebraic and $C^*$-algebraic, ``Quantum Probability and Related Topics,''
\emph{eds.\
M. \!Bo\.zejko, R. \!Lenczewski, W. \!M\l otkowski \& J. \!Wysocza\'nski},
Banach Center Publications, Polish Academy of Sciences,
Warsaw 2006, \emph{in press}.

\bibitem[$\text{LS}_3$]
{LSqsde}
J.M.\,Lindsay and A.G.\,Skalski,
On quantum stochastic differential equations, \emph{J. Math. Anal. Appl.} (2006), doi:10.1016/j.jmaa.2006.07.105


\bibitem[$\text{LS}_4$]
{LSqscc2}
J.M.\,Lindsay and A.G.\,Skalski,
Quantum stochastic convolution cocycles II, \emph{preprint}.


\bibitem[LW$_1$]{lwjfa}
J.M.\,Lindsay and S.J.\,Wills,
Markovian cocycles on
operator algebras, adapted to a Fock filtration,
\emph{J.\,Funct.\,Anal.}
\textbf{178} (2000) no.\,2, 269--305.

\bibitem[LW$_2$]{lwptrf}
J.M.\,Lindsay and S.J.\,Wills,
Existence, positivity, and contractivity
for quantum stochastic flows with infinite dimensional noise,
\emph{Probab.\,Theory Related Fields} \textbf{116} (2000), 505--543.



\bibitem[LW$_3$]{lwex}
J.M.\,Lindsay and S.J.\,Wills,
Existence of Feller cocycles on a $C^*$-algebra,
\emph{Bull.\,London Math.\,Soc.}
\textbf{33} (2001) no.\,5, 613--621.

\bibitem[LW$_4$]{lwhom}
J.M.\,Lindsay and S.J.\,Wills,
Homomorphic Feller cocycles on a $C^*$-algebra,
\emph{J.\,London Math.\,Soc.\ (2)}
\textbf{68} (2003) no.\,1, 255--272.




\bibitem[Mey]{mey}
P.-A.\,Meyer,
``Quantum Probability for Probabilists,'' 2nd Edition,
Lecture Notes in Mathematics
\textbf{1538},
Springer-Verlag,
Berlin 1995.

\bibitem[Pal]{Pal}
A.\,Pal, A counterexample on idempotent states on a compact quantum group,
 \emph{Lett.\,Math.\,Phys.}  \textbf{37}  (1996) no.\,1, 75--77.

\bibitem[Par]{par}
K.R.\,Parthasarathy,
``Introduction to Quantum Stochastic Calculus,''
Birkh\"{a}user, Basel 1992.

\bibitem[Ped]{Pedersen}
G.K.\, Pedersen,
``$C^*$-algebras and their automorphism groups", London Mathematical Society Monographs \textbf{14},
Academic Press Inc., London-New York 1979.


\bibitem[Pis$_1$]{Pisiersim}
G.\,Pisier,
``Similarity problems and completely bounded maps,"  Lecture Notes in Mathematics  \textbf{1618},
Springer-Verlag, Berlin 2001.


\bibitem[Pis$_2$]{Pisier}
G.\,Pisier,
``Introduction to Operator Space Theory,'' London Mathematical Society Lecture Note Series
\textbf{294},
CUP, Cambridge 2003.


\bibitem[QIIP]{qiip}
``Quantum Independent Increment Processes,
Vol. \!I: From Classical Probability to Quantum Stochastics;
Vol. \!II: Structure of Quantum L\'evy Processes,
Classical Probability and Physics,''
\emph{eds.\ U. \!Franz \& M. \!Sch\"urmann},
Lecture Notes in Mathematics \textbf{1865}; \textbf{1866},
Springer-Verlag, Heidelberg 2005.


\bibitem[Rin]{Ringrose}
 J.R.\,Ringrose, Automatic continuity of derivations of operator algebras.
\emph{J. London Math. Soc.} \textbf{5} (1972), 432--438.

\bibitem [Sak]{Sakai}
S.\,Sakai, On a conjecture of Kaplansky.
\emph{T\^ohoku Math. J.} \textbf{12} (1960), 31--33.

\bibitem[Sch]{schu}
M.\,Sch\"{u}rmann,
``White Noise on Bialgebras,''
Lecture Notes in Mathematics \textbf{1544},
Springer, Heidelberg 1993.

\bibitem[ScS] {Skeide}
M.\,Sch\"urmann and M.\,Skeide,
Infinitesimal generators on the quantum group $SU_q(2)$,
\emph{Infin.\,Dimens.\,Anal.\,Quantum Prob.\, Relat.\,Top.}
\textbf{1} (1998) no.\,4, 573-598.


\bibitem[Ska] {CPQSCC}
 A.\,Skalski, Completely positive quantum stochastic convolution cocycles and their dilations,
 \emph{Math.\,Proc.\,Camb.\,Phil.\,Soc.} (to appear).


\bibitem[Ske] {Skeiden}
 M.\,Skeide, L\'evy processes and tensor product systems of Hilbert modules,
  \emph{in}
``Quantum Probability and Infinite Dimensional
 Analysis, From Foundations to Applications''
\emph{eds.\ M.\,Sch\"urmann \& U.\,Franz,}
World Scientific Publishing, Singapore 2004,
pp. 492--503.

\bibitem[Swe] {Sweedler}
M.E.\,Sweedler, ``Hopf Algebras,''
Benjamin, New York 1969.

\bibitem[$\text{Wor}_1$]{CMP}
S.L.\,Woronowicz,
Compact matrix pseudogroups,
\emph{Comm. \!Math. \!Phys.}
\textbf{111} (1987) no. \!4, 613--665.

\bibitem [$\text{Wor}_2$]
{wor2}
S.L.\,Woronowicz,
Compact quantum groups,
\emph{in} ``Sym\'etries Quantiques,'' Proceedings, Les Houches 1995,
\emph{eds.\ A. \!Connes, K. \!Gawedzki \& J. \!Zinn-Justin},
North-Holland, Amsterdam 1998, pp.\,845--884.

\end{thebibliography}
\end{document}